\documentclass{aims}
\usepackage{listings}
\usepackage{color,xcolor}
\usepackage{tikz}
\usetikzlibrary{decorations.pathreplacing,calligraphy}
\usepackage{amsmath,bbm}
\usepackage{amssymb}
\usetikzlibrary{shapes.geometric, arrows}
\tikzstyle{startstop} = [rectangle, rounded corners, minimum width=3cm, minimum height=2cm,text centered, draw=black, fill=red!30]
\tikzstyle{arrow} = [thick,->,>=stealth]
\usetikzlibrary{positioning,chains,fit,shapes,calc}
\definecolor{grey1}{rgb}{0.0,0.0,0.0}
\definecolor{grey2}{rgb}{0.2,0.2,0.2}
\definecolor{grey3}{rgb}{0.4,0.4,0.4}
\definecolor{grey4}{rgb}{0.5,0.5,0.5}

\definecolor{myblue}{RGB}{100,100,160}
\definecolor{mygreen}{RGB}{80,160,80}

\lstdefinestyle{highlight-fonts}{
	basicstyle=\ttfamily\mdseries\footnotesize,
	commentstyle=\ttfamily\mdseries\itshape\footnotesize,
	moredelim=[l][\color{\hiddenA}\ttfamily\mdseries\itshape\footnotesize]{\%\#},
	moredelim=[l][\color{\hiddenB}\ttfamily\mdseries\itshape\footnotesize]{\%!},
	keywordstyle=\ttfamily\mdseries\footnotesize,
	keywordstyle={[2]\ttfamily\bfseries\footnotesize},
}

\lstdefinestyle{highlight-colors}{
	basicstyle=\color{grey1}\ttfamily\mdseries\footnotesize,
	commentstyle=\color{grey2}\ttfamily\mdseries\itshape\footnotesize,
	moredelim=[l][\color{\hiddenA}\ttfamily\mdseries\itshape\footnotesize]{\%\#},
	moredelim=[l][\color{\hiddenB}\ttfamily\mdseries\itshape\footnotesize]{\%!},
	keywordstyle=\color{grey2}\ttfamily\mdseries\footnotesize,
	keywordstyle={[2]\color{grey2}\ttfamily\bfseries\footnotesize},
}

\lstdefinestyle{highlight-small}{
	basicstyle=\color{grey1}\ttfamily\mdseries\scriptsize,
	commentstyle=\color{grey2}\ttfamily\mdseries\itshape\scriptsize,
	moredelim=[l][\color{\hiddenA}\ttfamily\mdseries\itshape\scriptsize]{\%\#},
	moredelim=[l][\color{\hiddenB}\ttfamily\mdseries\itshape\scriptsize]{\%!},
	keywordstyle=\color{grey2}\ttfamily\mdseries\scriptsize,
	keywordstyle={[2]\color{grey2}\ttfamily\bfseries\scriptsize},
}

\lstdefinestyle{gobble1}{
	xleftmargin=0.3em,
}

\lstdefinestyle{gobble2}{
	xleftmargin=-0.9em,
}

\lstdefinestyle{gobble3}{
	xleftmargin=-2.1em,
}

\lstdefinelanguage{coco}[]{Matlab}{
	keywords={},
	morekeywords={},
	morekeywords=[2]{classdef,properties,private,protected,public,%
Access,Static,methods,if,function,end,for,while,else,elseif,switch,%
case,otherwise,do,repeat,until},
	otherkeywords={{,end},:end,(end,end),\{end,end\},[end,end],/private,private/},
	morecomment=[l]{},
	basicstyle=\ttfamily\mdseries\footnotesize,
	commentstyle=\ttfamily\mdseries\footnotesize,
	moredelim=[l][\color{\hiddenA}\ttfamily\mdseries\footnotesize]{\%\#},
	moredelim=[l][\color{\hiddenB}\ttfamily\mdseries\footnotesize]{\%!},
	keywordstyle=\ttfamily\mdseries\footnotesize,
	keywordstyle={[2]\ttfamily\mdseries\footnotesize},
	numbers=none,
	numberstyle=\scriptsize,
	numbersep=1em,
	breaklines=false,
	breakatwhitespace=true,
	breakindent=2.5em,
	showlines=false,
	lineskip=-0.2ex,
	frame=none,
	fontadjust=true,
	columns=[c]fixed,
	basewidth={0.575em,0.45em},
	fontadjust=true,
	tabsize=3,
	showstringspaces=false,
	aboveskip=1.5\medskipamount,
	belowskip=1.5\medskipamount,
	xleftmargin=0.925em,
	xrightmargin=0em,
	rangeprefix={\%!},
	includerangemarker=false,
	belowcaptionskip=\bigskipamount
}

\lstdefinelanguage{coco-highlight-fonts}[]{coco}{
	style=highlight-fonts,
}

\lstdefinelanguage{coco-highlight-colors}[]{coco}{
	style=highlight-colors,
}

\lstdefinelanguage{coco-small}[]{coco}{
	style=highlight-small,
}

\lstdefinelanguage{coco-highlight}[]{coco-highlight-colors}{}

\newcommand{\mcode}[1]{{\lstinline[language=coco-highlight,basewidth={0.6em,0.45em}]|#1|}}

\usepackage{graphicx}
\usepackage[colorlinks=true]{hyperref}
\hypersetup{urlcolor=blue, citecolor=red}
\def\currentvolume{X}
\def\currentissue{X}
\def\currentyear{200X}
\def\currentmonth{XX}
\def\ppages{X--XX}

\theoremstyle{definition}

\DeclareMathOperator{\rg}{rg}

\newcommand{\jrem}[1]{#1}
\newcommand{\transp}{\mathsf{T}}
\title[Sensitivity Analysis Using the Adjoint Method] %Use the shortened version of the full title
      {Sensitivity Analysis for Periodic Orbits and Quasiperiodic Invariant Tori Using the Adjoint Method}

\author[Harry Dankowicz and Jan Sieber]{}


\email{danko@illinois.edu}
\email{j.sieber@exeter.ac.uk}

\thanks{$^*$ Corresponding author: danko@illinois.edu}

\begin{document}
\maketitle

% Enter the first author's name and address:
\centerline{\scshape Harry Dankowicz$^*$}
\medskip
{\footnotesize
% please put the address of the first author
 \centerline{Department of Mechanical Science and Engineering}
   \centerline{University of Illinois at Urbana-Champaign}
   \centerline{Urbana, IL 68101, USA}
}

\medskip

\centerline{\scshape Jan Sieber}
\medskip
{\footnotesize
 % please put the address of the second  and third author
 \centerline{College of Engineering, Mathematics
and Physical Sciences}
   \centerline{University of Exeter}
   \centerline{Exeter EX4 4QF, United Kingdom}
}

\bigskip

% The name of the associate editor will be entered by an editorial staff
% "Communicated by the associate editor name" is not needed for special issue.
 \centerline{(Communicated by the associate editor name)}

%The abstract of your paper
\begin{abstract}
This paper presents a rigorous framework for the continuation of solutions to nonlinear constraints and the simultaneous analysis of the sensitivities of test functions to constraint violations at each solution point using an adjoint-based approach. By the linearity of a problem Lagrangian in the associated Lagrange multipliers, the formalism is shown to be directly amenable to analysis using the \textsc{coco} software package, specifically its paradigm for staged problem construction. The general theory is illustrated in the context of algebraic equations and boundary-value problems, with emphasis on periodic orbits in smooth and hybrid dynamical systems, and quasiperiodic invariant tori of flows. In the latter case, normal hyperbolicity is used to prove the existence of continuous solutions to the adjoint conditions associated with the sensitivities of the orbital periods to parameter perturbations and constraint violations, even though the linearization of the governing boundary-value problem lacks a bounded inverse, as required by the general theory. An assumption of transversal stability then implies that these solutions predict the asymptotic phases of trajectories based at initial conditions perturbed away from the torus. Example \textsc{coco} code is used to illustrate the minimal additional investment in setup costs required to append sensitivity analysis to regular parameter continuation.
\textbf{200} words.
\end{abstract}

\section{Introduction}

In the search for optimal designs of engineering structures, the adjoint method is a commonplace tool for computing the sensitivities of response functions to variations in the design parameters given a set of constraints that respect physical laws and discrete design decisions (see the review \cite{Tortorelli94} and extensive references cited therein, as well as the more recent literature, e.g., \cite{rubino2018adjoint}). As observed in this literature, compared to direct approaches that necessitate the inversion of a problem linearization, the adjoint approach is advantageous in cases where the number of response functions is small relative to the number of design variables.

In cases where the constraints are in the form of initial-value problems, the adjoint method results in adjoint differential equations that must be solved in backward time in order to determine the desired sensitivities (see, e.g., \cite{TraversoMagri2019}, where such an adjoint approach is used to perform real-time data assimilation for a predictive reduced-order model of a problem in nonlinear thermoacoustics). Modifications to such a formalism may be derived to handle so-called hybrid systems, in which continuously differentiable time-dependence is interrupted by event-driven state-space jumps and vector field discontinuities associated, e.g., with collisions in a multibody mechanical system~\cite{Sandu20} and switches in a relay.

A reduced form of the adjoint formalism occurs in the computation of the asymptotic phase of a limit cycle~\cite{Govaerts06}. Generalizations to systems with delay~\cite{Ahsan21,novivcenko2012phase} or event-driven discontinuities~\cite{Park18,Shirisaka17} address the necessary modifications to the general theory and the associated adjoint conditions. A further generalization of the reduced formalism occurs in the analysis of asymptotic phases associated with perturbations off invariant tori in~\cite{demir2010}.

With its origin in constrained-design optimization, the adjoint approach derives from the analysis of a Lagrangian (functional) in terms of the design variables, the response variables, and an auxiliary set of Lagrange multipliers. The adjoint conditions are then obtained by imposing vanishing infinitesimal variations of the Lagrangian with respect to the design and response variables. Due to the linearity of the Lagrangian in the Lagrange multipliers, the adjoint conditions are also linear in these variables, albeit with coefficients that depend nonlinearly on the design and response variables. This linearity lends itself to a staged approach to problem construction in which constraints and the associated terms in the adjoint conditions are added successively according to an ordering that is sensible to the designer. Such a staged approach also naturally respects a modular approach to constraints, in which constraints are grouped and composed through the application of glue.

The \textsc{Matlab}-based software package \textsc{coco} \cite{dankowicz2013recipes,COCO} supports such a staged approach of problem construction of both constraints and adjoint conditions. It is, therefore, able to integrate sensitivity analysis with parameter continuation along implicitly-defined solution manifolds, almost out of the box. In recent papers~\cite{ahsan2020optimization,li2017staged,li2020optimization}, this functionality was used to demonstrate a successive continuation approach to locating extrema of a single objective function along families of solutions to nonlinear boundary-value problems. In the present paper\footnote{\jrem{The code included in this paper constitutes fully executable scripts. Complete code, including that used to generate the results in Fig.~\ref{fig:invc}, is available at \url{https://github.com/jansieber/adjoint-sensitivity2022-supp}.}}, we demonstrate this functionality in the context of the more general problem of sensitivity analysis.

With this past work in mind, the purpose of this paper is to collect, and present original rigorous derivations of, earlier results, specifically those pertaining to phase reduction analysis near limit cycles and quasiperiodic invariant tori. Here, and in some contrast to the above-cited references, the emphasis is on a formalism that respects its origin in a problem Lagrangian and that, thereby, provides new interpretations of the results obtained in the existing literature. Particular care is taken in investigating the solvability of the adjoint conditions, especially in the case of quasiperiodic invariant tori. Several example low-level encodings in \textsc{coco} demonstrate the utility of its core functionality, as well as the ways in which the analysis may be generalized to more complicated examples.

The remainder of the paper is organized as follows. In Section~\ref{sec:adjoint-based sensitivity analysis}, we develop the theoretical foundation and illustrate its application to an algebraic problem using explicit analysis and with an implementation in \textsc{coco}. Section~\ref{sec:Examples} contrasts the direct approach with the adjoint method for deriving well-known results regarding sensitivities associated with an initial-value problem, a final-value problem, and a periodic boundary-value problem in smooth dynamical systems. It further considers the implications of event-driven discontinuities of state and vector field and presents results for a periodic orbit with a single state-space discontinuity. A more demanding analysis of two-dimensional quasiperiodic invariant tori is considered in Section~\ref{sec:quasiperiodic invariant tori}. Here, normal hyperbolicity is relied upon to establish the solvability of the adjoint conditions and transversal stability is assumed only to derive the asymptotic phase dynamics. A discussion of the regularity of the problem linearization in Section~\ref{sec:regularity} shows that while the linear problem lacks a bounded inverse, the sensitivities computed using the adjoint method represent components of the problem Jacobian with bounded inverses. The theoretical results are illustrated in the context of continuation along a family of invariant curves of a two-dimensional map in Section~\ref{sec:invc}. Section~\ref{sec:construction} reviews the construction paradigm of \textsc{coco} with particular emphasis on the stages introduction of constraints and contributions to the adjoint conditions. These principles are illustrated in Section~\ref{sec:invc:coco} using an explicit \textsc{coco} encoding of the invariant curve problem from Section~\ref{sec:invc}. The paper concludes in Section~\ref{sec:conclusions} with a brief outlook on ideas worth exploring.

\section{Adjoint-based sensitivity analysis}
\label{sec:adjoint-based sensitivity analysis}
We precede our treatment of periodic orbits and quasiperiodic invariant tori with a general theoretical formalism and an algebraic example that illustrates the promise of the adjoint-based analysis and its implementation in the \textsc{coco} software package. The section assumes some familiarity with \textsc{coco}'s design principles and command structure, as can be gleaned from extensive tutorial documentation included with the \textsc{coco} release and the monograph \cite{dankowicz2013recipes}.

\subsection{Theoretical preliminaries}
\label{sec:preliminaries}
Consider the Lagrangian
\begin{equation}
\label{eq:genlag}
    L(u,\mu,\lambda,\eta):=\langle\lambda,\Phi(u)\rangle_{\mathcal{R}_\Phi}+\langle\eta,\Psi(u)-\mu\rangle_{\mathbb{R}^{n_\Psi}}
\end{equation}
in terms of the \textit{continuation variables} $u\in\mathcal{U}_\Phi$, \textit{zero functions} $\Phi:\mathcal{U}_\Phi\rightarrow\mathcal{R}_\Phi$, \textit{monitor functions} $\Psi:\mathcal{U}_\Phi\rightarrow\mathbb{R}^{n_\Psi}$, \textit{continuation parameters} $\mu\in\mathbb{R}^{n_\Psi\times 1}$, and \textit{adjoint variables}  $\lambda\in\mathcal{R}^\ast_\Phi$ and $\eta\in\mathbb{R}^{1\times n_\Psi}$, where $\mathcal{R}_\Phi^*$ denotes the dual space of linear functionals on $\mathcal{R}_\Phi$. The variations of $L$ with respect to $u$, $\lambda$, and $\eta$ vanish at a point $(u,\mu,\lambda,\eta)$ provided that
\begin{equation}\label{eq:gensys}
    \Phi(u)=0,\,\Psi(u)-\mu=0,\, (D\Phi(u))^\ast\lambda+ (D\Psi(u))^\ast\eta=0.
\end{equation}
Consistent with the \textsc{coco} syntax, we refer to the first two of these equations as an \textit{extended continuation problem}. The final equation is the corresponding set of \textit{adjoint conditions}, which are linear in the adjoint variables.

Suppose that the \textit{zero problem} $\Phi(u)=0$ is regular\footnote{The equation $\Phi=0$ on $\mathcal{U}$ is said to be \textit{regular with dimensional deficit} $d$ at a solution point $\tilde{u}$ if there exists a function $\Psi:\mathcal{U}\rightarrow\mathbb{R}^d$ such that the map
$F:u\mapsto(\Phi(u),\Psi(u))$ is continuously Frech\'{e}t differentiable on a neighborhood of $\tilde{u}$ and $DF(\tilde{u})$
has a bounded inverse.} with dimensional deficit $d<n_\Psi$ at a solution $\tilde{u}$ and let $\tilde{\mu}:=\Psi(\tilde{u})$. Let $\mathbb{I}$, $\mathbb{J}_1$, and $\mathbb{J}_2$ be three disjoint subsets of $\{1,\ldots,n_\Psi\}$ such that $\|\mathbb{I}\cup\mathbb{J}_1\|=d$, $\mathbb{J}_1\cup\mathbb{J}_2=\{1,\ldots,n_\Psi\}\setminus\mathbb{I}$, and the reduced continuation problem
\begin{equation}
    \Phi(u)=0,\,\Psi_{\mathbb{I}\,\cup\,\mathbb{J}_1}(u)-\tilde{\mu}_{\mathbb{I}\,\cup\,\mathbb{J}_1}=0
\end{equation}
is regular at $\tilde{u}$ with zero dimensional deficit. Finally, denote by  $(\tilde{\lambda},\tilde{\eta})$ the solution to the adjoint conditions
\begin{equation}\label{gen:adjoint}
   (D\Phi(\tilde{u}))^\ast\lambda+(D\Psi(\tilde{u}))^\ast\eta=0
\end{equation}
with $\tilde{\eta}_{\mathbb{J}_2\setminus k}=0$ and $\tilde{\eta}_k=1$ for some $k\in\mathbb{J}_2$. Then, we show below that the components of $(\tilde{\lambda},\tilde{\eta}_{\mathbb{I}\,\cup\,\mathbb{J}_1})$ describe the sensitivities of the monitor function $\Psi_k(u)$ to violations of the zero problem and perturbations in $\mu_{\mathbb{I}\,\cup\,\mathbb{J}_1}$, respectively, that perturb $u$ away from $\tilde{u}$. 

Indeed, for small $(\delta\Phi,\delta\mu_{\mathbb{I}\,\cup\,\mathbb{J}_1})\in\mathcal{R}_\Phi\times\mathbb{R}^{n_\Psi}$, the perturbed problem
\begin{equation}
    \Phi(u)=\delta\Phi,\,\Psi_{\mathbb{I}\,\cup\,\mathbb{J}_1}(u)-\tilde{\mu}_{\mathbb{I}\cup\mathbb{J}_1}=\delta\mu_{\mathbb{I}\,\cup\,\mathbb{J}_1}
\end{equation}
has a locally unique solution $u=\tilde{u}+\delta u$, where
\begin{equation}\label{gen:implicit:diff}
    \delta\Phi=D\Phi(\tilde{u})\delta u+\mathcal{O}(\|\delta u\|^2),\quad \delta\mu_{\mathbb{I}\,\cup\,\mathbb{J}_1}=D\Psi_{\mathbb{I}\,\cup\,\mathbb{J}_1}(\tilde{u})\delta u+\mathcal{O}(\|\delta u\|^2)\mbox{.}
\end{equation}
The perturbation $\delta u$ results in a perturbation to the value of the monitor function  $\Psi_k$ given by
\begin{align}\label{gen:diff:insert}
    \delta\Psi_k&=D\Psi_k(\tilde{u})\delta u+\mathcal{O}(\|\delta u\|^2)=\langle(D\Psi_k(\tilde{u}))^\ast,\delta u\rangle_{\mathcal{U}_\Phi}+\mathcal{O}(\|\delta u\|^2),
\end{align}
where the second equality follows from the formal definition of the adjoint. We may determine the sensitivity of $\Psi_k$ to violations of the zero problem and perturbations in $\mu_{\mathbb{I}\,\cup\,\mathbb{J}_1}$ by solving \eqref{gen:implicit:diff} for $\delta u$ in terms of $\delta\Phi$ and $\delta\mu_{\mathbb{I}\,\cup\,\mathbb{J}_1}$, substituting the result into the middle expression of \eqref{gen:diff:insert} and identifying the coefficients in front of $\delta\Phi$ and $\delta\mu_{\mathbb{I}\,\cup\,\mathbb{J}_1}$, respectively. We call this the \emph{direct differentiation} approach.

Alternatively, upon inserting the solution $(\tilde{\lambda},\tilde{\eta})$ from \eqref{gen:adjoint} with $\tilde{\eta}_{\mathbb{J}_2\setminus k}=0$ and $\tilde{\eta}_k=1$ in the rightmost expression of \eqref{gen:diff:insert}, again applying the formal definition of adjoints, and using \eqref{gen:implicit:diff}, we arrive at the well-known adjoint sensitivity formula
\begin{align}\label{gen:adj:sensitivity}
  \delta\Psi_k&=-\langle\tilde{\lambda},\delta\Phi\rangle_{\mathcal{R}_\Phi}-\langle\tilde{\eta}_{\mathbb{I}\,\cup\,\mathbb{J}_1},\delta\mu_{\mathbb{I}\,\cup\,\mathbb{J}_1}\rangle_{\mathbb{R}^{n_\Psi}}+\mathcal{O}(\|(\delta\Phi,\delta \mu_{\mathbb{I}\,\cup\,\mathbb{J}_1})\|^2).
\end{align}
This confirms the claimed relationship between the components of $(\tilde{\lambda},\tilde{\eta}_{\mathbb{I}\,\cup\,\mathbb{J}_1})$ and the sensitivities of $\Psi_k(u)$ to violations of the zero problem and perturbations in $\mu_{\mathbb{I}\,\cup\,\mathbb{J}_1}$, respectively, that perturb $u$ away from $\tilde{u}$. The \emph{adjoint-based} analysis thus produces the sensitivities directly from the solution of \eqref{gen:adjoint} without the need to first invert \eqref{gen:implicit:diff}.

In the context of parameter continuation using the \textsc{coco} package, elements of $\mathbb{I}$ index \textit{inactive} continuation parameters that impose permanent constraints on the continuation variables during a given continuation run. In contrast, elements of $\mathbb{J}_1\cup\mathbb{J}_2$ index \textit{active} continuation parameters that track the values of the corresponding monitor functions during continuation. Different choices of $\mathbb{J}_2$ and $k$ in the assignment of values to elements of $\tilde{\eta}$ then yield different combinations of sensitivities.

The \textsc{coco} construction paradigm does not commit the user to a particular choice of values for the adjoint variables $\eta$. Instead, much like with the dimension of the solution manifold, this decision can be made at runtime (e.g., using the function \mcode{coco_set_parival} applied to the corresponding complementary continuation parameters), after problem construction. This is a consequence of the general formalism which produces terms in the adjoint conditions associated with each individual zero or monitor function but omits the conditions that would follow from also performing variations with respect to elements of $\mu$.

Given $\tilde{u}$, once the adjoint conditions have been derived, a particular choice of $\mathbb{J}_2$ and $k$ has been made, and the components of $\eta_{\mathbb{J}_2}$ have been assigned accordingly to $0$ or $1$, a solution for $\tilde{\lambda}$ and $\tilde{\eta}_{\mathbb{I}\,\cup\,\mathbb{J}_1}$ may be obtained by solving a linear problem. If a solution has been found for one such $\tilde{u}$, we may continue such a solution under variations in $\mu_{\mathbb{J}_1\cup\,\mathbb{J}_2}$. 

Since \textsc{coco} treats the adjoint conditions as part of an augmented continuation problem, which is assumed to be nonlinear, a solution to the adjoint conditions is typically computed with \textsc{coco} using continuation under variations in $\eta_k$ from $0$ to $1$. For the same reason, \textsc{coco} does not currently take advantage of the linear form in order to solve a matrix version of the adjoint conditions with columns of $\eta_{\mathbb{J}_2}$ representing different choices of $k$. Simultaneous analysis in the current release of \textsc{coco} for different choices of $\mathbb{J}_2$ and $k$ instead requires duplicate copies of the adjoint conditions in the augmented continuation problem. 

\subsection{An algebraic example}
\label{sec:illustration}
We illustrate the ideas of the previous section by considering the Lagrangian
\begin{align}
    L&=\lambda_{\mathrm{am,1}}\left( c_1^2-a_1^2-b_1^2\right)+\lambda_{\mathrm{am,2}}\left(c_2^2-a_2^2-b_2^2\right)+\lambda_\mathrm{fr}\left(\omega_1-\omega_2-\epsilon\right)\nonumber\\
    &\qquad+\lambda_{\mathrm{de},1,\mathrm{re}}\left((1-\omega_1^2)a_1+2\zeta\omega_1b_1-1\right)+\lambda_{\mathrm{de},1,\mathrm{im}}\left((1-\omega_1^2)b_1-2\zeta\omega_1a_1\right)\nonumber\\
    &\qquad+\lambda_{\mathrm{de},2,\mathrm{re}}\left( (1-\omega_2^2)a_2+2\zeta\omega_2b_2-1\right)+\lambda_{\mathrm{de},2,\mathrm{im}}\left((1-\omega_2^2)b_2-2\zeta\omega_2a_2\right)\nonumber\\
    &\qquad+\eta_\mathrm{da}\left(c_1-c_2-\Delta\right)+\eta_\mathrm{av}\left(\frac{\omega_1+\omega_2}{2}-\bar{\omega}\right)+\eta_\epsilon\left(\epsilon-\epsilon_0\right)+\eta_\zeta\left(\zeta-\zeta_0\right)
\end{align}
in terms of the continuation variables $a_1$, $b_1$, $c_1$, $a_2$, $b_2$, $c_2$, $\omega_1$, $\omega_2$, $\zeta$, and $\epsilon$, continuation parameters $\Delta$, $\bar{\omega}$, $\epsilon_0$, and $\zeta_0$, and adjoint variables $\lambda_{\mathrm{am,1}}$, $\lambda_{\mathrm{am,2}}$, $\lambda_\mathrm{fr}$, $\lambda_{\mathrm{de},1,\mathrm{re}}$, $\lambda_{\mathrm{de},1,\mathrm{im}}$, $\lambda_{\mathrm{de},2,\mathrm{re}}$, $\lambda_{\mathrm{de},2,\mathrm{im}}$, $\eta_\mathrm{da}$, $\eta_\mathrm{av}$, $\eta_\epsilon$, and $\eta_\zeta$. The corresponding nonlinear algebraic zero problem
\begin{gather}
    c_1^2-a_1^2-b_1^2=0,\,c_2^2-a_2^2-b_2^2=0,\,\omega_1-\omega_2-\epsilon=0,\label{eq:const1}\\
    (1-\omega_1^2)a_1+2\zeta\omega_1b_1-1=0,\,(1-\omega_1^2)b_1-2\zeta\omega_1a_1=0,\\
    (1-\omega_2^2)a_2+2\zeta\omega_2b_2-1=0,\,(1-\omega_2^2)b_2-2\zeta\omega_2a_2=0.\label{eq:const3}
\end{gather}
corresponds to the search for harmonic solutions of the form 
\begin{align}
    x_1(t)=c_1\cos(\omega_1 t+\phi_1) =a_1\cos\omega_1 t+b_1\sin\omega_1 t,\\x_2(t)=c_2\cos(\omega_2 t+\phi_2)=a_2\cos\omega_2 t+b_2\sin\omega_2 t
\end{align}
of the linear differential equations
\begin{equation}
\label{eq:linosc}
    \ddot{x}_1+2\zeta\dot{x}_1+x_1=\cos\omega_1 t,\,\ddot{x}_2+2\zeta\dot{x}_2+x_2=\cos\omega_2 t,
\end{equation}
for excitation frequencies $\omega_1$ and $\omega_2$ that differ by $\epsilon$. 

Restricting attention to $c_1,c_2>0$, straightforward analysis of \eqref{eq:const1}-\eqref{eq:const3} yields the relationship
\begin{align}
    \Delta&=    \frac{1}{\sqrt{(1-(\bar{\omega}+\epsilon_0/2)^2)^2+4\zeta_0^2(\bar{\omega}+\epsilon_0/2)^2}}\nonumber\\
    &\qquad-\frac{1}{\sqrt{(1-(\bar{\omega}-\epsilon_0/2)^2)^2+4\zeta_0^2(\bar{\omega}-\epsilon_0/2)^2}}
    \label{eq:Deltaexp}
\end{align}
between the four continuation parameters. We may obtain the sensitivity of $\Delta$ with respect to $\bar{\omega}$, $\zeta_0$, or $\epsilon_0$ by direct differentiation of the right-hand side. For example, the sensitivity of $\Delta$ with respect to $\bar{\omega}$ behaves as $o(\epsilon_0)$ in the limit as $\epsilon_0\rightarrow 0$ provided that
\begin{equation}
    3\bar{\omega}^6+5(2\zeta_0^2-1)\bar{\omega}^4+(16\zeta_0^4-16\zeta_0^2+1)\bar{\omega}^2+1-2\zeta_0^2=0,
\end{equation}
corresponding to an inflection point in the frequency response curve for the differential equation
\begin{equation}
    \ddot{x}+2\zeta_0\dot{x}+x=\cos\bar{\omega}t.
\end{equation}
For $\zeta_0=0.1$, two such inflection points are located at $\bar{\omega}\approx 0.92$ and $1.06$.

Alternatively, in the absence of an explicit solution, we may consider the linearization
\begin{gather}
    2c_1\delta_{c_1}-2a_1\delta_{a_1}-2b_1\delta_{b_1}=\delta_{\mathrm{am},1},\label{eq:sec2lin1}\\
    2c_2\delta_{c_2}-2a_2\delta_{a_2}-2b_2\delta_{b_2}=\delta_{\mathrm{am},2},\\
    \delta_{\omega_1}-\delta_{\omega_2}-\delta_\epsilon=\delta_\mathrm{fr},\\
    -2\omega_1a_1\delta_{\omega_1}+(1-\omega_1^2)\delta_{a_1}+2\omega_1b_1\delta_\zeta+2\zeta b_1\delta_{\omega_1}+2\zeta\omega_1\delta_{b_1}=\delta_{\mathrm{de},1,\mathrm{re}},\\
    -2\omega_1b_1\delta_{\omega_1}+(1-\omega_1^2)\delta_{b_1}-2\omega_1a_1\delta_\zeta-2\zeta a_1\delta_{\omega_1}-2\zeta\omega_1\delta_{a_1}=\delta_{\mathrm{de},1,\mathrm{im}},\\
    -2\omega_2a_2\delta_{\omega_2}+(1-\omega_2^2)\delta_{a_2}+2\omega_2b_2\delta_\zeta+2\zeta b_2\delta_{\omega_2}+2\zeta\omega_2\delta_{b_2}=\delta_{\mathrm{de},2,\mathrm{re}},\\
    -2\omega_2b_2\delta_{\omega_2}+(1-\omega_2^2)\delta_{b_2}-2\omega_2a_2\delta_\zeta-2\zeta a_2\delta_{\omega_2}-2\zeta\omega_2\delta_{a_2}=\delta_{\mathrm{de},2,\mathrm{im}},\\
    \delta_{c_1}-\delta_{c_2}=\delta_\Delta,\,
    \frac{\delta_{\omega_1}+\delta_{\omega_2}}{2}=\delta_{\bar{\omega}},\,
    \delta_\epsilon=\delta_{\epsilon_0},\,\delta_\zeta=\delta_{\zeta_0}\label{eq:sec2lin8}
\end{gather}
around some solution to the extended continuation problem. The sensitivities of $\Delta$ with respect to $\bar{\omega}$, $\epsilon_0$, $\zeta_0$, and the constraint violations $\delta_{\mathrm{am},1}$, $\delta_{\mathrm{am},2}$, $\delta_\mathrm{fr}$, $\delta_{\mathrm{de},1,\mathrm{re}}$, $\delta_{\mathrm{de},1,\mathrm{im}}$, $\delta_{\mathrm{de},2,\mathrm{re}}$, and $\delta_{\mathrm{de},2,\mathrm{im}}$ may then be obtained by solving for $\delta_\Delta$ and inspecting the coefficients in front of $\delta_{\bar{\omega}}$, $\delta_{\epsilon_0}$, $\delta_{\zeta_0}$, $\delta_{\mathrm{am},1}$, $\delta_{\mathrm{am},2}$, $\delta_\mathrm{fr}$, $\delta_{\mathrm{de},1,\mathrm{re}}$, $\delta_{\mathrm{de},1,\mathrm{im}}$, $\delta_{\mathrm{de},2,\mathrm{re}}$, and $\delta_{\mathrm{de},2,\mathrm{im}}$, respectively. As shown in the previous section, these sensitivities may also be obtained directly by solving the adjoint conditions
\begin{gather}
    2c_1\lambda_{\mathrm{am},1}+\eta_\mathrm{da}=0,\label{eq:adjobj}\\
    2c_2\lambda_{\mathrm{am},2}-\eta_\mathrm{da}=0\\
    -2a_1\lambda_{\mathrm{am},1}+(1-\omega_1^2)\lambda_{\mathrm{de},1,\mathrm{re}}-2\zeta\omega_1\lambda_{\mathrm{de},1,\mathrm{im}}=0,\\
    -2a_2\lambda_{\mathrm{am},2}+(1-\omega_2^2)\lambda_{\mathrm{de},2,\mathrm{re}}-2\zeta\omega_2\lambda_{\mathrm{de},2,\mathrm{im}}=0,\\
    -2b_1\lambda_{\mathrm{am},1}+2\zeta\omega_1\lambda_{\mathrm{de},1,\mathrm{re}}+(1-\omega_1^2)\lambda_{\mathrm{de},1,\mathrm{im}}=0,\\
    -2b_2\lambda_{\mathrm{am},2}+2\zeta\omega_2\lambda_{\mathrm{de},2,\mathrm{re}}+(1-\omega_2^2)\lambda_{\mathrm{de},2,\mathrm{im}}=0,\\
    \lambda_\mathrm{fr}+2(\zeta b_1-\omega_1a_1)\lambda_{\mathrm{de},1,\mathrm{re}}-2(\zeta a_1+\omega_1b_1)\lambda_{\mathrm{de},1,\mathrm{im}}+\eta_\mathrm{av}/2=0,\\
    -\lambda_\mathrm{fr}+2(\zeta b_2-\omega_2a_2)\lambda_{\mathrm{de},2,\mathrm{re}}-2(\zeta a_2+\omega_2b_2)\lambda_{\mathrm{de},2,\mathrm{im}}+\eta_\mathrm{av}/2=0,\\
    2\omega_1 b_1\lambda_{\mathrm{de},1,\mathrm{re}}-2\omega_1 a_1\lambda_{\mathrm{de},1,\mathrm{im}}+2\omega_2 b_2\lambda_{\mathrm{de},2,\mathrm{re}}-2\omega_2 a_2\lambda_{\mathrm{de},2,\mathrm{im}}+\eta_\zeta=0,\\
    -\lambda_\mathrm{fr}+\eta_\epsilon=0,
    \label{eq:adjobjlast}
\end{gather}
with $\eta_\mathrm{da}=1$, for the remaining adjoint variables. Notably, the linearized equations \eqref{eq:sec2lin1}-\eqref{eq:sec2lin8} impose eleven linear constraints on the eleven unknowns $\delta_{a_1}$, $\delta_{b_1}$, $\delta_{c_1}$, $\delta_{a_2}$, $\delta_{b_2}$, $\delta_{c_2}$, $\delta_{\omega_1}$, $\delta_{\omega_2}$, $\delta_\zeta$, $\delta_\epsilon$, and $\delta_\Delta$, but the only expression of interest is the solution for $\delta_\Delta$ and this is only obtained as a linear combination of constraint violations and parameter perturbations. In contrast, the adjoint conditions impose ten linear constraints on the ten unknowns $\lambda_{\mathrm{am,1}}$, $\lambda_{\mathrm{am,2}}$, $\lambda_\mathrm{fr}$, $\lambda_{\mathrm{de},1,\mathrm{re}}$, $\lambda_{\mathrm{de},1,\mathrm{im}}$, $\lambda_{\mathrm{de},2,\mathrm{re}}$, $\lambda_{\mathrm{de},2,\mathrm{im}}$, $\eta_\mathrm{av}$, $\eta_\epsilon$, and $\eta_\zeta$, all of which are of interest and immediately represent the sought sensitivities of $\Delta$ with respect to the individual constraint violations and perturbations of the continuation parameters $\bar{\omega}$, $\epsilon_0$, and $\zeta_0$.

We proceed to implement this analysis in \textsc{coco} \jrem{with fixed $\epsilon=0.01$ and $\zeta=0.1$. In this case, the derivative of the right-hand side of \eqref{eq:Deltaexp} with respect to $\bar{\omega}$ vanishes for $\bar{\omega}\approx 0.92043919$ and $1.06205837$ with residuals of $2\times 10^{-7}$ and $4\times 10^{-7}$, respectively.} \textsc{coco}-compatible encodings of the zero function $\Phi$ and its Jacobian\footnote{In the absence of explicit encodings of second derivatives, \textsc{coco} relies on a suitable finite-difference approximation of these derivatives, as necessary.} $D\Phi$ are as follows.
\begin{lstlisting}[language=coco-highlight]
function [data, f] = phi(prob, data, u)

v = num2cell(u);
[a1, b1, c1, a2, b2, c2, o1, o2, ze, ep] = deal(v{:});
f = [c1^2-a1^2-b1^2; c2^2-a2^2-b2^2; o1-o2-ep;
     (1-o1^2)*a1+2*ze*o1*b1-1; (1-o1^2)*b1-2*ze*o1*a1;
     (1-o2^2)*a2+2*ze*o2*b2-1; (1-o2^2)*b2-2*ze*o2*a2];

end
\end{lstlisting}
\begin{lstlisting}[language=coco-highlight]
function [data, J] = dphi(prob, data, u)

v = num2cell(u);
[a1, b1, c1, a2, b2, c2, o1, o2, ze, ep] = deal(v{:});
J = [-2*a1,-2*b1,2*c1,0,0,0,0,0,0,0;
     0,0,0,-2*a2,-2*b2,2*c2,0,0,0,0;
     0,0,0,0,0,0,1,-1,0,-1;
     1-o1^2,2*ze*o1,0,0,0,0,-2*o1*a1+2*ze*b1,0,2*o1*b1,0;
     -2*ze*o1,1-o1^2,0,0,0,0,-2*o1*b1-2*ze*a1,0,2*o1*a1,0;
     0,0,0,1-o2^2,2*ze*o2,0,0,-2*o2*a2+2*ze*b2,2*o2*b2,0;
     0,0,0,-2*ze*o2,1-o2^2,0,0,-2*o2*b2-2*ze*a2,2*o2*a2,0];

end
\end{lstlisting}
Similarly, \textsc{coco}-compatible encodings of the monitor function $\Psi$ and its Jacobian $D\Psi$ are given below.
\begin{lstlisting}[language=coco-highlight]
function [data, f] = psi(prob, data, u)

v = num2cell(u);
[a1, b1, c1, a2, b2, c2, o1, o2, ze, ep] = deal(v{:});
f = [c1-c2; (o1+o2)/2; ep; ze];

end
\end{lstlisting}
\begin{lstlisting}[language=coco-highlight]
function [data, J] = dpsi(prob, data, u)

v = num2cell(u);
[a1, b1, c1, a2, b2, c2, o1, o2, ze, ep] = deal(v{:});
J = [0,0,1,0,0,-1,0,0,0,0;
     0,0,0,0,0,0,1/2,1/2,0,0;
     0,0,0,0,0,0,0,0,0,1;
     0,0,0,0,0,0,0,0,1,0];

end
\end{lstlisting}
We proceed to construct and initialize the corresponding extended continuation problem using the following sequence of commands.
\begin{lstlisting}[language=coco-highlight]
>> prob = coco_prob();
>> prob = coco_add_func(prob, 'phi', @phi, @dphi, [], 'zero', ...
     'u0', [-0.49; 4.9; 4.9; 0; 5; 5; 1.01; 1; .1; .01]);
>> prob = coco_add_func(prob, 'psi', @psi, @dpsi, [], 'inactive', ...
     {'da', 'av', 'ep', 'ze'}, 'uidx', 1:10);
\end{lstlisting}
Here, an initial solution guess for the continuation variables is encapsulated in the vector array following the \mcode{'u0'} flag in the first call to \mcode{coco_add_func}. The \mcode{'inactive'} flag in the second call to \mcode{coco_add_func} implies that the corresponding continuation parameters, here designated by the string labels \mcode{'da'}, \mcode{'av'}, \mcode{'ep'}, and \mcode{'ze'}, are initially inactive.

The commands
\begin{lstlisting}[language=coco-highlight]
>> prob = coco_add_adjt(prob, 'phi');
>> prob = coco_add_adjt(prob, 'psi', {'e.da','e.av','e.ep','e.ze'}, ...
     'aidx', 1:10);
\end{lstlisting}
append the corresponding adjoint conditions with additional \textit{complementary continuation parameters}, designated by the string labels \mcode{'e.da'}, \mcode{'e.av'}, \mcode{'e.ep'}, and \mcode{'e.ze'}, and associated with complementary monitor functions whose values equal $\eta_\mathrm{da}$, $\eta_\mathrm{av}$, $\eta_\epsilon$, and $\eta_\zeta$, respectively, in the notation of this section. These are inactive by default. All the adjoint variables are initialized with their default value $0$ by this construction.

To obtain the sensitivities of $\Delta$ with $\mathbb{I}=\{3,4\}$, $\mathbb{J}_1=\{2\}$, and $\mathbb{J}_2=\{1\}$, we first perform continuation along the one-dimensional solution manifold obtained by allowing \mcode{'da'}, \mcode{'e.da'}, \mcode{'e.av'}, \mcode{'e.ep'}, and \mcode{'e.ze'} to vary, while keeping \mcode{'av'}, \mcode{'ep'}, and \mcode{'ze'} fixed.
\begin{lstlisting}[language=coco-highlight]
>> coco(prob, 'run', [], 1, {'da', 'e.da', 'e.av', 'e.ep', 'e.ze'}, ...
     {[], [0 1]})
\end{lstlisting}
\begin{lstlisting}[language=coco-small]
    STEP   DAMPING               NORMS              COMPUTATION TIMES
  IT SIT     GAMMA     ||d||     ||f||     ||U||   F(x)  DF(x)  SOLVE
   0                          2.40e-01  1.00e+01    0.0    0.0    0.0
   1   1  1.00e+00  3.72e-02  6.80e-04  1.00e+01    0.0    0.0    0.0
   2   1  1.00e+00  9.75e-05  4.76e-09  1.00e+01    0.0    0.0    0.0
   3   1  1.00e+00  6.83e-10  3.56e-15  1.00e+01    0.0    0.0    0.0

 ...  LABEL  TYPE          da        e.da           e.av         e.ep         e.ze
 ...      1  EP     -7.3832e-02   0.0000e+00   0.0000e+00   0.0000e+00   0.0000e+00
 ...      2         -7.3832e-02   2.5752e-01   1.2060e+00   1.8906e+00  -3.1450e-01
 ...      3         -7.3832e-02   5.6770e-01   2.6587e+00   4.1678e+00  -6.9333e-01
 ...      4         -7.3832e-02   8.7788e-01   4.1114e+00   6.4451e+00  -1.0722e+00
 ...      5  EP     -7.3832e-02   1.0000e+00   4.6833e+00   7.3416e+00  -1.2213e+00
\end{lstlisting}
We proceed to extract the solution data from the point obtained with \mcode{'e.da'} equal to $1$ and use this to reconstruct and reinitialize the augmented continuation problem. 
\begin{lstlisting}[language=coco-highlight]
>> prob = coco_prob();
>> chart = coco_read_solution('phi', 'run', 5, 'chart');
>> prob = coco_add_func(prob, 'phi', @phi, @dphi, [], 'zero', ...
     'u0', chart.x);
>> prob = coco_add_func(prob, 'psi', @psi, @dpsi, [], 'inactive', ...
     {'da', 'av', 'ep', 'ze'}, 'uidx', 1:10);
>> chart = coco_read_adjoint('phi', 'run', 5, 'chart');
>> prob = coco_add_adjt(prob, 'phi', 'l0', chart.x);
>> chart = coco_read_adjoint('psi', 'run', 5, 'chart');
>> prob = coco_add_adjt(prob, 'psi', {'e.da','e.av','e.ep','e.ze'}, ...
     'aidx', 1:10, 'l0', chart.x);
\end{lstlisting}
Continuation along the solution manifold obtained for fixed \mcode{'ep'}, \mcode{'ze'}, and \mcode{'e.da'}, and for \mcode{'av'} in the interval $[0.5,2.5]$ is then triggered by the commands
\begin{lstlisting}[language=coco-highlight]
>> prob = coco_set(prob, 'cont', 'ItMX', 500, 'NPR', 100);
>> coco(prob, 'run', [], 1, {'da', 'av', 'e.av', 'e.ep', 'e.ze'}, ...
     {[], [0.5 2.5]})
\end{lstlisting}
where the \mcode{'ItMX'} and \mcode{'NPR'} options regulate the number of continuation steps in each direction from the initial solution and the screen output frequency. The resultant output is shown below
\begin{lstlisting}[language=coco-small]
    STEP   DAMPING               NORMS              COMPUTATION TIMES
  IT SIT     GAMMA     ||d||     ||f||     ||U||   F(x)  DF(x)  SOLVE
   0                          1.06e-14  1.90e+01    0.0    0.0    0.0

 ...  LABEL  TYPE            da           av         e.av         e.ep         e.ze
 ...      1  EP     -7.3832e-02   1.0050e+00   4.6833e+00   7.3416e+00  -1.2213e+00
 ...      2  FP     -2.0587e-01   1.0621e+00  -1.7633e-06   2.0541e+01  -3.9743e+00
 ...      3         -1.2158e-01   1.1556e+00  -9.0966e-01   1.2164e+01  -1.0384e+00
 ...      4  EP     -1.7965e-03   2.5000e+00  -2.6831e-03   1.7965e-01  -3.4810e-04

 ...  LABEL  TYPE            da           av         e.av         e.ep         e.ze
 ...      5  EP     -7.3832e-02   1.0050e+00   4.6833e+00   7.3416e+00  -1.2213e+00
 ...      6          1.7834e-01   9.2050e-01   2.9412e-03  -1.7795e+01   3.7095e+00
 ...      7  FP      1.7834e-01   9.2044e-01   9.2729e-07  -1.7795e+01   3.7069e+00
 ...      8  EP      1.6854e-02   5.0000e-01  -7.5097e-02  -1.6857e+00   1.8071e-02
\end{lstlisting}
Here, the points denoted by \mcode{FP} are fold points (local extrema) in the quantity $\Delta$ and coincide with the loci of sign changes in the sensitivity of $\Delta$ with respect to $\bar{\omega}$, i.e., approximate inflection points. The corresponding values of \mcode{'av'} \jrem{(computed with the default residual tolerance of $10^{-6}$ and rounded off to $0.92044$ and $1.0621$)} agree with those predicted by the closed-form analysis in the first part of this section.

\section{Sensitivities along solutions to ODEs}
\label{sec:Examples}
In this section, we derive several known results about the sensitivity of quantities associated with the behavior of smooth and hybrid dynamical systems to violations of a governing set of differential and algebraic constraints using either linearization of the governing constraints or the adjoint-based approach. Notably, we interpret the adjoint conditions liberally as matrix equalities, enabling simultaneous derivation of the sensitivities of a vector of monitor functions to constraint violations and perturbations in the remaining continuation parameters.

\subsection{A single trajectory segment}
\label{sec:Forward dynamics}
Consider the flow  $F(t,x,p)$ corresponding to the autonomous vector field $f:\mathbb{R}^n\times\mathbb{R}^q\rightarrow\mathbb{R}^n$, such that $x(\tau;T,x_0,p):=F(T\tau,x_0,p)$ is the unique solution to the initial-value problem with rescaled time:
\begin{equation}\label{forw:ode}
    x'=Tf(x,p),\quad x(0)=x_0.
\end{equation}
We apply the formalism of Section~\ref{sec:preliminaries} to
quantities involved in \eqref{forw:ode} in order to evaluate the sensitivities of $x(1;T,x_0,p)$ to its arguments. As in the previous section, we first perform direct differentiation and compare intermediate steps and final results to the predictions of the adjoint analysis.

\subsubsection*{Direct differentiation} 
It follows directly by
differentiation at $\tau=1$ that
\begin{equation}
    \partial_{T}x(1;T,x_0,p)=\partial_tF(T,x_0,p)=f\left(F(T,x_0,p),p\right)=f(x(1;T,x_0,p),p),
\end{equation}
and we obtain the Jacobians $\partial_{x_0}x(1;T,x_0,p)=\partial_xF(T,x_0,p)$ and $\partial_{p}x(1;T,x_0,p)=\partial_pF(T,x_0,p)$ as the solutions at $\tau=1$ to the standard first-order variational initial-value problems
\begin{equation}
\label{eq:vareq}
    X'=T\partial_xf\left(x(\tau;T,x_0,p),p\right)X,\quad X(0)=I_n
\end{equation}
and
\begin{equation}
\label{eq:vareqp}
    P'=T\partial_xf\left(x(\tau;T,x_0,p),p\right)P+T\partial_pf\left(x(\tau;T,x_0,p),p\right),\quad P(0)=0,
\end{equation}
respectively, where $I_n$ denotes the corresponding identity matrix.

Consider, for example, a perturbation of $x_0$ along the vector $f(x_0,p)$. Since
\begin{equation}
    \frac{d}{d\tau}f(x(\tau;T,x_0,p),p)=T\partial_xf\left(x(\tau;T,x_0,p),p\right)f(x(\tau;T,x_0,p),p),
\end{equation}
it follows that
\begin{equation}
\label{eq:vfmap}
    f(x(\tau;T,x_0,p),p)=X(\tau)f(x_0,p),
\end{equation}
i.e., that, to linear order, $x(1;T,x_0,p)$ is perturbed by the vector $f(x(1;T,x_0,p),p)$.
Moreover, by the variation-of-parameters formula, we obtain
\begin{equation}
    P(\tau)=TX(\tau)\int_0^{\tau} X(\sigma)^{-1}\partial_pf\left(x(\sigma;T,x_0,p),p\right)\,\mathrm{d}\sigma.
\end{equation}
This gives the sensitivity $P(1$) of $x(1;T,x_0,p)$ with respect to $p$ directly in terms of a convolution integral.

\subsubsection*{Regularity}
Alternatively, in  notation consistent with that of the previous section, let $\mathcal{U}_\Phi=C^1([0,1];\mathbb{R}^n)\times\mathbb{R}\times\mathbb{R}^q$ be the space of continuation variables, and define
\begin{equation}
    \mathcal{U}\ni u=(x(\cdot),T,p)\mapsto \Phi(u):=x'(\cdot)-Tf(x(\cdot),p)\in C^0([0,1];\mathbb{R}^n)=\mathcal{R}_\Phi
\end{equation}
and
\begin{equation}
\label{eq:monitorforward}
    \Psi(u):=\begin{pmatrix}x(1)\\x(0)\\T\\p\end{pmatrix}\in\mathbb{R}^{2n+1+q}.
\end{equation}
In this case, every solution to the zero problem is regular with dimensional deficit $n+1+q$. Indeed, given a solution $(\tilde{x}(\cdot),T,p)$ and small constraint violation $\delta_\mathrm{ode} (\cdot)$, the linearized equation
\begin{equation}
    \delta_\mathrm{ode}=\delta'_x-T\partial_xf(\tilde{x},p)\delta_x-f(\tilde{x},p)\delta_T-T\partial_pf(\tilde{x},p)\delta_p
\end{equation}
implies that
\begin{align}
\label{eq:linfor}
    \delta_x(\tau)&=X(\tau)\int_0^{\tau} X(\sigma)^{-1}\delta_\mathrm{ode}(\sigma)\,\mathrm{d}\sigma+X(\tau)\delta_x(0)+\tau f(\tilde{x}(\tau),p)\delta_T+P(\tau)\delta_p.
\end{align}
We thus obtain $\delta_x(\cdot)$ in terms of the $n+1+q$ quantities $\delta_x(0)$, $\delta_T$, and $\delta_p$.

\subsubsection*{Adjoint analysis}
Following \eqref{eq:genlag}, we write the Lagrangian
\begin{align}
\label{eq:lagrangetrajseg}
    L&=\int_0^1\lambda^{\transp}(\tau)\left(x'(\tau)-Tf(x(\tau),p)\right)\,\mathrm{d}\tau+\eta_{x(1)}^{\transp}\left(x(1)-\mu_{x(1)}\right)\nonumber\\&\quad+\eta_{x(0)}^{\transp}\left(x(0)-\mu_{x(0)}\right)+\eta_{T}\left(T-\mu_{T}\right)+\eta_{p}^{\transp}\left(p-\mu_{p}\right)
\end{align}
in terms of the adjoint variables\footnote{Here, the dual space $\mathcal{R}^\ast_\Phi$ is the space of functions of bounded variation. We restrict attention to the subspace of continuously differentiable functions $\lambda(\cdot)$ to allow the use of integration of parts when evaluating variations of $L$.} $\lambda(\cdot)\in C^1([0,1];\mathbb{R}^n)$, $\eta_{x(1)},\eta_{x(0)}\in\mathbb{R}^n$, $\eta_T\in\mathbb{R}$, and $\eta_p\in\mathbb{R}^q$. The adjoint conditions \eqref{gen:adjoint} may then be written
\begin{align}
\label{eq:adjforfirst}
\delta x(\cdot)&: &0&=-\lambda^{\prime\,{\transp}}-T\lambda^{\transp}\partial_xf(x,p),\\
\delta x(1)&:&0&=\lambda^{\transp}(1)+\eta_{x(1)}^{\transp},\\
\delta x(0)&:&0&=-\lambda^{\transp}(0)+\eta_{x(0)}^{\transp},\\
\delta T&:&0&=-\int_0^1\lambda^{\transp}(\tau)f(x(\tau),p)\,\mathrm{d}\tau+\eta_T,\\
\delta p&:&0&=-\int_0^1\lambda^{\transp}(\tau)T\partial_pf(x(\tau),p)\,\mathrm{d}\tau+\eta_p^{\transp},
\label{eq:adjforlast}
\end{align}
where we have labeled each condition by the corresponding variation $\delta u$.

We let $\mathbb{I}=\emptyset$, $\mathbb{J}_1=\{n+1,\ldots,2n+1+q\}$, and $\mathbb{J}_2=\{1,\ldots,n\}$. With reference to \eqref{eq:monitorforward}, this choice corresponds to computing the sought sensitivities of $x(1)$ with respect to violations of the governing differential constraint and perturbations to $T$, $x(0)$, and $p$, respectively that drive $(x(\cdot),T,p)$ away from $(x(\cdot),T,p)$. It follows by inspection using \eqref{forw:ode}, \eqref{eq:vareq}, \eqref{eq:vareqp}, and \eqref{eq:adjforfirst}-\eqref{eq:adjforlast} that
\begin{align}
\label{eq:adjfornext1}
    \tilde{\lambda}^{\transp}f(x,p)&\equiv\eta_T=-\eta^{\transp}_{x(1)}f(x(1),p),\\
    \label{eq:adjfornext2}\tilde{\lambda}^{\transp}X&\equiv\eta^{\transp}_{x(0)}=-\eta^{\transp}_{x(1)}X(1),\\
    \label{eq:adjfornext3}\tilde{\lambda}^{\transp}P\big|_{\tau=1}&=\eta_p^{\transp}=-\eta^{\transp}_{x(1)}P(1).
\end{align}
We set $\eta_{x(1)}=I_n$ and obtain the sensitivities of $x(1)$ with respect to variations in $T$, $x(0)$, and $p$ from the quantities $f(x(1),p)$, $X(1)$, and $P(1)$, respectively, as expected from the results obtained using direct differentiation. In this case, \eqref{eq:adjfornext1} and \eqref{eq:adjfornext2} also imply that
\begin{equation}
    X(1)f(x(0),p)=f(x(1),p),
\end{equation}
consistent with \eqref{eq:vfmap}.

\subsection{A single segment with a Poincar\'{e} section}
\label{sec:Poincaresection}
\jrem{We append the zero function $(x(\cdot),T,p)\mapsto h_\mathrm{ps}(x(1),p)$ to the previous construction by adding the term $\lambda_\mathrm{ps}h_\mathrm{ps}(x(1),p)$ to the corresponding Lagrangian \eqref{eq:lagrangetrajseg}. We use the notation $(\cdot)_\mathrm{ps}$ to indicate the association with a Poincar\'{e} section
\begin{equation}
    \label{eq:poincare:seection}
    \{x:h_\mathrm{ps}(x,p)=0\}.
\end{equation}
Then, every solution to the zero problem with nonzero Lie derivative
\begin{equation}
    \mathcal{L}_fh_\mathrm{ps}(x(1),p):=\partial_xh_\mathrm{ps}(x(1),p)f(x(1),p)    
\end{equation}
(i.e., that intersects the Poincar\'{e} section transversally)} is regular with dimensional deficit $n+q$. Indeed, consider the additional small constraint violation $\delta_h$, such that
\begin{equation}
    \partial_xh_\mathrm{ps}(x(1),p)\delta_x(1)+\partial_ph_\mathrm{ps}(x(1),p)\delta_p=\delta_h.
\end{equation}
Then, \eqref{eq:linfor} implies that
\begin{align}
\label{eq:linfor3}
    &\mathcal{L}_fh_\mathrm{ps}(x(1),p)\delta_T=\delta_h-\partial_xh_\mathrm{ps}(x(1),p)X(1)\int_0^{1} X(\sigma)^{-1}\delta_\mathrm{ode}(\sigma)\,\mathrm{d}\sigma\nonumber\\
    &\qquad- \partial_xh_\mathrm{ps}(x(1),p)X(1)\delta_x(0)-\left(\partial_ph_\mathrm{ps}(x(1),p)+ \partial_xh_\mathrm{ps}(x(1),p)P(1)\right)\delta_p,
\end{align}
which may be uniquely solved for $\delta_T$ if $\mathcal{L}_fh_\mathrm{ps}(x(1),p)\ne 0$. Consequently, $\delta_x(\cdot)$ and $\delta_T$ may be obtained in terms of the $n+q$ quantities $\delta_x(0)$ and $\delta_p$. 

With the additional constraint, the adjoint conditions become
\begin{align}
\delta x(\cdot)&: &0&=-\lambda_\mathrm{de}^{\prime\,{\transp}}-T\lambda_\mathrm{de}^{\transp}\partial_xf(x,p),\\
\delta x(1)&: &0&=\lambda_\mathrm{de}^{\transp}(1)+\lambda_\mathrm{ps}\partial_xh_\mathrm{ps}(x(1),p)+\eta_{x(1)}^{\transp},\\
\delta x(0)&: &0&=-\lambda_\mathrm{de}^{\transp}(0)+\eta_{x(0)}^{\transp},\\
\delta T&: &0&=-\int_0^1\lambda_\mathrm{de}^{\transp}(\tau)f(x(\tau),p)\,\mathrm{d}\tau+\eta_T,\\
\delta p&: &0&=-\int_0^1\lambda_\mathrm{de}^{\transp}(\tau)T\partial_pf(x(\tau),p)\,\mathrm{d}\tau+\lambda_\mathrm{ps}\partial_ph_\mathrm{ps}(x(1),p)+\eta_p^{\transp}.
\end{align}
This time, we let $\mathbb{I}=\emptyset$, $\mathbb{J}_1=\{n+1,\ldots,2n,2n+2,\ldots,2n+1+q\}$, and $\mathbb{J}_2=\{1,\ldots,n,2n+1\}$ in order to capture the sensitivities of $x(1)$ and $T$ with respect to constraint violations and perturbation in $x(0)$ and $p$. It follows by the identical steps to the previous section that
\begin{align}
    \lambda_\mathrm{de}^{\transp}f(x,p)&\equiv \eta_T=-\left(\eta^{\transp}_{x(1)}+\lambda_\mathrm{ps}\partial_xh_\mathrm{ps}(x(1),p)\right)f(x(1),p),\\
    \lambda_\mathrm{de}^{\transp}X&\equiv\eta^{\transp}_{x(0)}=-\left(\eta^{\transp}_{x(1)}+\lambda_\mathrm{ps}\partial_xh_\mathrm{ps}(x(1),p)\right)X(1),\\
    \lambda_\mathrm{de}^{\transp}P\big|_{\tau=1}&=\lambda_\mathrm{ps}\partial_ph_\mathrm{ps}(x(1),p)+\eta_p^{\transp}=-\left(\eta^{\transp}_{x(1)}+\lambda_\mathrm{ps}\partial_xh_\mathrm{ps}(x(1),p)\right)P(1)
\end{align}

By considering the case when $\eta_T=0$ and $\eta_{x(1)}=I_n$, we find that the sensitivities of $x(1)$ equal
\begin{equation}
    \frac{f(x(1),p)}{\mathcal{L}_fh_\mathrm{ps}(x(1),p)}
\end{equation}
with respect to $h_\mathrm{ps}(x(1),p)$,
\begin{equation}
    \Pi(x(1),p)X(1)
\end{equation}
with respect to $x(0)$, and 
\begin{equation}
    \Pi(x(1),p)P(1)-\frac{f(x(1),p)\partial_ph_\mathrm{ps}(x(1),p)}{\mathcal{L}_fh_\mathrm{ps}(x(1),p)}
\end{equation}
with respect to $p$, where the nullspaces of the projection matrix
\begin{equation}
    \Pi(x(1),p):= I_n-\frac{f(x(1),p)\partial_xh_\mathrm{ps}(x(1),p)}{\mathcal{L}_fh_\mathrm{ps}(x(1),p)}
\end{equation}
and its transpose
are spanned by $f(x(1),p)$ and $\left(\partial_xh_\mathrm{ps}(x(1),p)\right)^{\transp}$, respectively. Similarly, by considering the case when $\eta_T=1$ and $\eta_{x(1)}=0$, we find that the sensitivities of $T$ equal
\begin{equation}
    \frac{1}{\mathcal{L}_fh_\mathrm{ps}(x(1),p)}
\end{equation}
with respect to $h_\mathrm{ps}(x(1),p)$,
\begin{equation}
    -\frac{\partial_xh_\mathrm{ps}(x(1),p)}{\mathcal{L}_fh_\mathrm{ps}(x(1),p)}X(1)
\end{equation}
with respect to $x(0)$, and 
\begin{equation}
    -\frac{\partial_xh_\mathrm{ps}(x(1),p)P(1)+\partial_ph_\mathrm{ps}(x(1),p)}{\mathcal{L}_fh_\mathrm{ps}(x(1),p)}
\end{equation}
with respect to $p$. These conclusions also follow from implicit differentiation of the constraints
\begin{equation}
    h_\mathrm{ps}(x(1),p)=H,\,x(1)=F(T,x(0),p)
\end{equation}
with respect to the independent variables $H$, $x(0)$, and $p$, or by solving the linearized equations \eqref{eq:linfor} and \eqref{eq:linfor3} for $\delta_x(1)$ and $\delta_T$.

\subsection{Periodic orbits}
\label{sec:Periodic orbits}
\jrem{Next, append the zero function $(x(\cdot),T,p)\mapsto x(0)-x(1)$ to the previous construction by adding the term $\lambda_\mathrm{po}^{\transp}(x(0)-x(1))$ to the corresponding Lagrangian. Here, the notation $(\cdot)_\mathrm{po}$ reflects the corresponding imposition of the periodicity constraint $x(1)=x(0)$.} Then, every solution to the zero problem with $\mathcal{L}_fh_\mathrm{ps}(x(1),p)\ne 0$, and for which the eigenvalue $1$ of the monodromy matrix $X(1)$ is simple, is regular with dimensional deficit $q$. Here, consider the additional small constraint violation $\delta_\mathrm{po}$, such that
\begin{equation}
    \delta_x(0)-\delta_x(1)=\delta_\mathrm{po}.
\end{equation}
Then, \eqref{eq:linfor} and \eqref{eq:linfor3} imply that
\begin{align}
    &\left(I_n-\Pi(x(1),p)X(1)\right)\delta_x(0)=\delta_\mathrm{po}+\Pi(x(1),p)X(1)\int_0^{1} X(\sigma)^{-1}\delta_\mathrm{ode}(\sigma)\,\mathrm{d}\sigma\nonumber\\
    &\qquad+ \frac{f(x(1),p)}{\mathcal{L}_fh_\mathrm{ps}(x(1),p)}\delta_h+ \left(\Pi(x(1),p)P(1)-\frac{f(x(1),p)\partial_ph_\mathrm{ps}(x(1),p)}{\mathcal{L}_fh_\mathrm{ps}(x(1),p)}\right)\delta_p,
\end{align}
which may be uniquely solved for $\delta_x(0)$ provided that $\Pi(x(1),p)X(1)$ has no eigenvalues equal to $1$. Indeed, if there exists a vector $v$ such that $\Pi(x(1),p)X(1)v=v$, then $X(1)v-v$ must be parallel to $f(x(1),p)$ in violation of the assumption on $X(1)$ that the eigenvalue $1$ is simple. When invertibility holds, it follows after substitution in \eqref{eq:linfor3} and some simplification that
\begin{align}
\label{eq:perorbsen}
    \mathcal{L}_fh_\mathrm{ps}(x(1),p)\delta_T&=w^{\transp}\left(\delta_\mathrm{po}+P(1)\delta_p+\int_0^{1} X(\sigma)^{-1}\delta_\mathrm{ode}(\sigma)\,\mathrm{d}\sigma\right),
\end{align}
where the prefactor \begin{equation}
    w^{\transp}:=-\partial_xh_\mathrm{ps}(x(1),p)X(1)\left(I_n-\Pi(x(1),p)X(1)\right)^{-1}
\end{equation}
is the unique left eigenvector of $X(1)$ corresponding to the eigenvalue $1$ such that $w^{\transp}f(x(1),p)=-\mathcal{L}_fh_\mathrm{ps}(x(1),p)$.

With the additional periodicity constraint, the adjoint conditions now become
\begin{align}
\delta x(\cdot)&: &0&=-\lambda_\mathrm{de}^{\prime\,{\transp}}-T\lambda_\mathrm{de}^{\transp}\partial_xf(x,p),\\
\delta x(1)&: &0&=\lambda_\mathrm{de}^{\transp}(1)+\lambda_\mathrm{ps}\partial_xh_\mathrm{ps}(x(1),p)-\lambda_\mathrm{po}^{\transp}+\eta_{x(1)}^{\transp},\\
\delta x(0)&: &0&=-\lambda_\mathrm{de}^{\transp}(0)+\lambda_\mathrm{po}^{\transp}+\eta_{x(0)}^{\transp},\\
\delta T&: &0&=-\int_0^1\lambda_\mathrm{de}^{\transp}(\tau)f(x(\tau),p)\,\mathrm{d}\tau+\eta_T,\\
\delta p&: &0&=-\int_0^1\lambda_\mathrm{de}^{\transp}(\tau)T\partial_pf(x(\tau),p)\,\mathrm{d}\tau+\lambda_\mathrm{ps}\partial_ph_\mathrm{ps}(x(1),p)+\eta_p^{\transp}.
\end{align}
This time, we let $\mathbb{I}=\emptyset$, $\mathbb{J}_1=\{2n+2,\ldots,2n+1+q\}$, and $\mathbb{J}_2=\{1,\ldots,2n+1\}$. It follows again by inspection that
\begin{align}
    \lambda_\mathrm{de}^{\transp}f(x,p)&\equiv \eta_T=-\left(\eta^{\transp}_{x(1)}+\lambda_\mathrm{ps}\partial_xh_\mathrm{ps}(x(1),p)-\lambda_\mathrm{po}^{\transp}\right)f(x(1),p)\nonumber\\
    &=(\lambda_\mathrm{po}^{\transp}+\eta^{\transp}_{x(0)})f(x(1),p),\\
    \lambda_\mathrm{de}^{\transp}X&\equiv\lambda_\mathrm{po}^{\transp}+\eta^{\transp}_{x(0)}\nonumber\\
    &=-\left(\eta^{\transp}_{x(1)}+\lambda_\mathrm{ps}\partial_xh_\mathrm{ps}(x(1),p)-\lambda_\mathrm{po}^{\transp}\right)X(1),\\
    \lambda_\mathrm{de}^{\transp}P\big|_{\tau=1}&=\lambda_\mathrm{ps}\partial_ph_\mathrm{ps}(x(1),p)+\eta_p^{\transp}\nonumber\\&=-\left(\eta^{\transp}_{x(1)}+\lambda_\mathrm{ps}\partial_xh_\mathrm{ps}(x(1),p)-\lambda_\mathrm{po}^{\transp}\right)P(1).
\end{align}
By considering the case when $\eta_T=1$ and $\eta_{x(0)}=\eta_{x(1)}=0$, we obtain
\begin{equation}
    \lambda_\mathrm{ps}=0,\,\lambda_\mathrm{po}^{\transp}f(x(1),p)=1,\,\lambda_\mathrm{po}^{\transp}=\lambda_\mathrm{po}^{\transp}X(1),\text{ and }\eta_p^{\transp}=\lambda_\mathrm{po}^{\transp}P(1),
\end{equation}
and conclude that the sensitivities of $T$ with respect to $h_\mathrm{ps}(x(1),p)$, $x(0)-x(1)$, and $p$ equal $0$, the unique left eigenvector $w^{\transp}$ of the monodromy matrix $X(1)$ corresponding to the eigenvalue $1$ and such that $w^{\transp}f(x(1),p)=-1$, and the product $w^{\transp}P(1)$, respectively, consistent with the result in \eqref{eq:perorbsen}.

\subsection{Segmented trajectories in hybrid systems}
\label{sec:Hybrid dynamics}

We broaden the perspective to hybrid dynamical systems that include discrete jumps and switches between different vector fields. The results in this section generalize to trajectories consisting of any number of consecutive solution segments.

\jrem{For example, consider a time history of a hybrid dynamical system that evolves along a trajectory of a flow $F_1$ until a transversal intersection at $(x_0,p_0)$ with an \emph{event surface}
\begin{equation}
    \{x:h_\mathrm{es}(x,p)=0\},
\end{equation}
followed by evolution along a trajectory of a flow $F_2$ from a point $g(x_0,p_0)$ for some map $g$. In particular, consistent with the discussion in Section~\ref{sec:Poincaresection}, assume that $\partial_xh_\mathrm{es}(x_0,p_0)f_1(x_0,p_0)\ne 0$.} Consider the composition
\begin{equation}
    D(x,p)=F_2(-\sigma(x,p),g(F_1(\sigma(x,p),x,p),p),p),
\end{equation}
where  $h_\mathrm{es}(F_1(\sigma(x,p),x,p),p)=0$ defines $\sigma(x,p)\approx 0$ uniquely for $x$ and $p$ in some neighborhood of $x_0$ and $p_0$, respectively, such that $\sigma(x_0,p_0)=0$. The function $D$ is referred to as the \textit{zero-time discontinuity mapping}~\cite{PWS08} associated with the \textit{jump function} $g$ and the transversal intersection with the event surface $h_\mathrm{es}=0$. As a special case, $D(x_0,p_0)=g(x_0,p_0)$.

We obtain the sensitivities $\partial_xD(x_0,p_0)$ and $\partial_pD(x_0,p_0)$ by analyzing the zero problem
\begin{eqnarray}
    x_1'-\sigma f_1(x_1,p)=0,\,x_2'+\sigma f_2(x_2,p)=0,\\h_\mathrm{es}(x_1(1),p)=0,\,x_2(0)-g(x_1(1),p)=0\label{eq:bcs}
\end{eqnarray}
and monitor functions
\begin{equation}
    \Psi(u):=\begin{pmatrix}x_2(1)\\x_1(0)\\p\end{pmatrix},
\end{equation}
which yield the adjoint conditions
\begin{align}
    \delta x_1(\cdot)&: &0&=-\lambda_{\mathrm{de},1}^{\prime\,{\transp}}-\lambda_{\mathrm{de},1}^{\transp}\sigma f_1(x_1,p),\\
    \delta x_2(\cdot)&: &0&=-\lambda_{\mathrm{de},2}^{\prime\,{\transp}}+\lambda_{\mathrm{de},2}^{\transp}\sigma f_2(x_2,p),\\
    \delta x_1(1)&: &0&=\lambda_{\mathrm{de},1}^{\transp}(1)+\lambda_\mathrm{es}\partial_xh_\mathrm{es}(x_1(1),p)-\lambda_\mathrm{jf}^{\transp}\partial_xg(x_1(1),p),\\
    \delta x_1(0)&: &0&=-\lambda_{\mathrm{de},1}^{\transp}(0)+\eta^{\transp}_{x_1(0)},\\
    \delta x_2(1)&: &0&=\lambda_{\mathrm{de},2}^{\transp}(1)+\eta^{\transp}_{x_2(1)},\\
    \delta x_2(0)&: &0&=-\lambda_{\mathrm{de},2}^{\transp}(0)+\lambda_\mathrm{jf}^{\transp},\\
    \delta \sigma&: &0&=-\int_0^1\lambda_{\mathrm{de},1}^{\transp}f_1(x_1,p)\,\mathrm{d}\tau+\int_0^1\lambda_{\mathrm{de},2}^{\transp}f_2(x_2,p)\,\mathrm{d}\tau,\\
    \delta p&: &0&=-\int_0^1\lambda_{\mathrm{de},1}^{\transp}\sigma\partial_pf_1(x_1,p)\,\mathrm{d}\tau+\int_0^1\lambda_{\mathrm{de},2}^{\transp}\sigma\partial_pf_2(x_2,p)\,\mathrm{d}\tau\nonumber\\
    &&&\qquad+\lambda_\mathrm{es}\partial_ph_\mathrm{es}(x(1),p)-\lambda_\mathrm{jf}^{\transp}\partial_p g(x_1(1),p)+\eta_p^{\transp},
\end{align}
where $\eta=(\eta_{x_2(1)},\eta_{x_1(0)},\eta_p)$, and $\lambda_\mathrm{es}$ and $\lambda_\mathrm{jf}$ are adjoint variables associated with the boundary conditions \eqref{eq:bcs}.

Since $\sigma=0$ for $x_1(0)=x_0$ and $p=p_0$, it follows by inspection that this choice implies that $x_1\equiv x_0$, $x_2\equiv g(x_0,p_0)$,
\begin{align}
\label{eq:adjdiscmapsimpl1}
\lambda_{\mathrm{de},1}^{\transp}&\equiv\eta^{\transp}_{x_1(0)}=\lambda_\mathrm{jf}^{\transp}\partial_xg(x_0,p_0)-\lambda_\mathrm{es}\partial_xh_\mathrm{es}(x_0,p_0),\\
\lambda_{\mathrm{de},2}^{\transp}&\equiv-\eta^{\transp}_{x_2(1)}=\lambda_\mathrm{jf}^{\transp},\\
\eta_p^{\transp}&=\lambda_\mathrm{jf}^{\transp}\partial_pg(x_0,p_0)-\lambda_\mathrm{es}\partial_ph_\mathrm{es}(x_0,p_0),
\end{align}
and
\begin{equation}
\label{eq:adjdiscmapsimpl4}   \eta^{\transp}_{x_1(0)}f_1(x_0,p_0)=-\eta^{\transp}_{x_2(1)}f_2(g(x_0,p_0),p_0).
\end{equation}
Guided by the general theory, we let $\mathbb{I}=\emptyset$, $\mathbb{J}_1=\{n+1,\ldots,2n+q\}$, and $\mathbb{J}_2=\{1,\ldots,n\}$ consistent with computing the sensitivities of $x_2(1)=D(x_1(0),p)$ with respect to $x_1(0)$ and $p$ at $x_1(0)=x_0$ and $p=p_0$. With $\eta_{x_2(1)}=I_n$, we obtain
\begin{align}
    \partial_xD(x_0,p_0)&=-\eta^{\transp}_{x_1(0)}=\partial_xg(x_0,p_0)+\lambda_\mathrm{es}\partial_xh_\mathrm{es}(x_0,p_0),\\
    \partial_pD(x_0,p_0)&=-\eta_p^{\transp}=\partial_pg(x_0,p_0)+\lambda_\mathrm{es}\partial_ph_\mathrm{es}(x_0,p_0),
\end{align}
where
\begin{equation}
    \lambda_\mathrm{es}=\frac{f_2(g(x_0,p_0),p_0)-\partial_xg(x_0,p_0)f_1(x_0,p_0)}{\partial_xh_\mathrm{es}(x_0,p_0)f_1(x_0,p_0)}
\end{equation}
is obtained by multiplying \eqref{eq:adjdiscmapsimpl1} by $f_1(x_0,p_0)$ and simplifying using \eqref{eq:adjdiscmapsimpl4}. The same expressions may be obtained by implicit differentiation of the constraint
\begin{equation}
    h_\mathrm{es}(F_1(\sigma(x,p),x,p),p)=0
\end{equation}
and the defining expression for $D(x,p)$.

As a second example from the theory of hybrid systems, we consider a zero problem for a two-segment periodic orbit of period $T$ obtained by gluing together two individual zero problems of the second form considered in Section~\ref{sec:Forward dynamics} using additional boundary conditions. Specifically, we assume that
\begin{eqnarray}
\label{eq:twoseg1}
    x'_1-\sigma f(x_1,p)=0,\,h_\mathrm{es}(x_1(1),p)=0,\\x'_2-(T-\sigma)f(x_2,p)=0,\,h_\mathrm{ps}(x_2(1),p)=0,\\
    x_2(0)-g(x_1(1),p)=0,\,x_1(0)-x_2(1)=0,\label{eq:twoseg3}
\end{eqnarray}
where the terminal point of the first segment (of total duration $\sigma$) is constrained to the event surface $h_\mathrm{es}=0$ and mapped by the jump function $g$ to the initial point on the second segment (of total duration $T-\sigma$). In this case, any solution that intersects both the event surface and the Poincar\'{e} section at $h_\mathrm{ps}=0$ transversally, i.e., such that
\begin{align}
    \partial_x h_\mathrm{es}(x_1(1),p)f(x_1(1),p)&\ne 0,\label{eq:trans1}\\
    \partial_xh_\mathrm{ps}(x_2(1),p)f(x_2(1),p)&\ne 0,\label{eq:trans2}
\end{align}
is regular with dimensional deficit $q$. If we choose the monitor functions
\begin{equation}
    \Psi(u):=\begin{pmatrix}T\\p\end{pmatrix},
\end{equation}
we obtain the adjoint conditions
\begin{align}
   \delta x_1(\cdot)&: &0&=-\lambda_{\mathrm{de},1}^{\prime\,{\transp}}-\lambda_{\mathrm{de},1}^{\transp}\sigma f(x_1,p),\\
   \delta x_2(\cdot)&: &0&=-\lambda_{\mathrm{de},2}^{\prime\,{\transp}}-\lambda_{\mathrm{de},2}^{\transp}(T-\sigma) f(x_2,p),\\
   \delta x_1(1)&: &0&=\lambda_{\mathrm{de},1}^{\transp}(1)+\lambda_\mathrm{es}\partial_xh_\mathrm{es}(x_1(1),p)-\lambda_\mathrm{jf}^{\transp}\partial_x g(x_1(1),p),\\
   \delta x_1(0)&: &0&=-\lambda_{\mathrm{de},1}^{\transp}(0)+\lambda_\mathrm{po}^{\transp},\\
   \delta x_2(1)&: &0&=\lambda_{\mathrm{de},2}^{\transp}(1)+\lambda_\mathrm{ps}\partial_xh_\mathrm{ps}(x_2(1),p)-\lambda_\mathrm{po}^{\transp},\\
   \delta x_2(0)&: &0&=-\lambda_{\mathrm{de},2}^{\transp}(0)+\lambda_\mathrm{jf}^{\transp},\\
   \delta T&: &0&=-\int_0^1\lambda_{\mathrm{de},2}^{\transp}(\tau)f(x_2(\tau),p)\,\mathrm{d}\tau+\eta_T,\\
   \delta \sigma&: &0&=-\int_0^1\lambda_{\mathrm{de},1}^{\transp}(\tau)f(x_1(\tau),p)\,\mathrm{d}\tau+\int_0^1\lambda_{\mathrm{de},2}^{\transp}(\tau)f(x_2(\tau),p)\,\mathrm{d}\tau,\\
   \delta p&: &0&=-\int_0^1\lambda_{\mathrm{de},1}^{\transp}(\tau)\sigma\partial_pf(x_1(\tau),p)\,\mathrm{d}\tau+\lambda_\mathrm{es}\partial_ph_\mathrm{es}(x_1(1),p)\nonumber\\
   &&&\qquad-\int_0^1\lambda_{\mathrm{de},2}^{\transp}(\tau)(T-\sigma)\partial_pf(x_2(\tau),p)\,\mathrm{d}\tau+\lambda_\mathrm{ps}\partial_ph_\mathrm{ps}(x_2(1),p)\nonumber\\
   &&&\qquad-\lambda_\mathrm{jf}^{\transp}\partial_pg(x_1(1),p)+\eta_p^{\transp},
\end{align}
where $\eta=(\eta_T,\eta_p)$, and $\lambda_\mathrm{es}$, $\lambda_\mathrm{ps}$, $\lambda_\mathrm{jf}$, and $\lambda_\mathrm{po}$ are adjoint variables associated with the four boundary conditions.

As before, the functions $\lambda_{\mathrm{de},1}^{\transp}f(x_1,p)$ and $\lambda_{\mathrm{de},2}^{\transp}f(x_2,p)$ are constant and both equal to $\eta_T$. By transversality \eqref{eq:trans2}, it follows that $\lambda_\mathrm{ps}=0$ and, consequently, that $\lambda_\mathrm{po}=\lambda_{\mathrm{de},1}(0)=\lambda_{\mathrm{de},2}(1)$. Similarly, from \eqref{eq:trans1}, it follows that
\begin{equation}
    \lambda_\mathrm{es}=\lambda_{\mathrm{de},2}^{\transp}(0)\frac{\partial_x g(x_1(1),p)f(x_1(1),p)-f_2(x_2(0),p)}{\partial_x h_\mathrm{es}(x_1(1),p)f(x_1(1),p)}
\end{equation}
and, consequently,
\begin{equation}
    \lambda_{\mathrm{de},1}^{\transp}(1)=\lambda_{\mathrm{de},2}^{\transp}(0)\partial_xD(x_1(1),p).
\end{equation}
Since
\begin{equation}
\lambda_{\mathrm{de},1}^{\transp}(0)=\lambda_{\mathrm{de},1}^{\transp}(1)\partial_x F(\sigma,x_1(0),p)    
\end{equation}
and
\begin{equation}
    \lambda_{\mathrm{de},2}^{\transp}(0)=\lambda_{\mathrm{de},2}^{\transp}(1)\partial_xF(T-\sigma,x_2(0),p),
\end{equation}
it follows that
\begin{equation}
\label{eq:compeigvec}
    \lambda_{\mathrm{de},2}^{\transp}(1)=\lambda_{\mathrm{de},2}^{\transp}(1)\partial_xF(T-\sigma,x_2(0),p)\partial_xD(x_1(1),p)\partial_x F(\sigma,x_1(0),p).
\end{equation}
Further, from the final adjoint condition, after some manipulation we obtain
\begin{align}
    \eta_p^{\transp}&=\lambda_{\mathrm{de},2}^{\transp}(1)\bigg(\partial_xF(T-\sigma,x_2(0),p)\partial_xD(x_1(1),p)\partial_pF(\sigma,x_1(0),p)\nonumber\\
    &\qquad+\partial_xF(T-\sigma,x_2(0),p)\partial_pD(x_1(1),p)+\partial_pF(T-\sigma,x_2(0),p)\bigg).
\end{align}

Inspired by these expressions, we use the discontinuity mapping to define the period-$T$ flow map $G_T$ for fixed $\sigma$ and small changes in $x_1(0)$ and $p$ as follows
\begin{equation}
    G_T(x_1(0),p):=F(T-\sigma,D(F(\sigma,x_1(0),p),p),p).
\end{equation}
It then follows by comparison with \eqref{eq:compeigvec} that $\lambda_\mathrm{po}^{\transp}$ is the unique left eigenvector of the \emph{monodromy matrix}
\begin{equation}
    \partial_xG_T(x_1(0),p)
\end{equation}
corresponding to the eigenvalue $1$ and such that $\lambda_\mathrm{po}^{\transp}f(x_1(0),p)=\eta_T$, and that
\begin{equation}
    \eta_p^{\transp}=\lambda_\mathrm{po}^{\transp}\partial_pG_T(x_1(0),p),
\end{equation}
just like in the case of a smooth periodic orbit. If we let $\mathbb{I}=\emptyset$, $\mathbb{J}_1=\{2,\ldots,q+1\}$, and $\mathbb{J}_2=\{1\}$, we obtain the sensitivities of $T$ with respect to $h_\mathrm{es}$, $h_\mathrm{ps}$, $x_2(0)-g(x_1(1),p)$, $x_1(0)-x_2(1)$, and $p$ from $-\lambda_\mathrm{es}$, $-\lambda_\mathrm{ps}$, $-\lambda_\mathrm{jf}^{\transp}$, $-\lambda_\mathrm{po}^{\transp}$, and $-\eta_p^{\transp}$, respectively.

\section{Quasiperiodic invariant tori}
\label{sec:quasiperiodic invariant tori}
\jrem{In this section, we demonstrate the utility of the adjoint approach for a special class of normally hyperbolic invariant manifolds, namely transversally stable, quasiperiodic invariant tori. We show that this approach may be used to compute the linearization of the corresponding stable fiber projection, as a generalization of the concept of asymptotic phase for periodic orbits. Indeed, in this case, stable fibers span the entire neighborhood of the invariant torus, such that the asymptotic phase can be easily observed, for example, in numerical simulations.} In contrast to the treatment in previous sections, the discussion is here concerned with an infinite-dimensional problem for which determining the regularity of the zero problem is a nontrivial task.

\subsection{Revisiting the periodic orbit}\label{sec:po:revisit}
Before turning to the infinite-dimensional case, however, we return again to the analysis in Section~\ref{sec:Periodic orbits} of a periodic orbit in a smooth vector field and consider its implications to the dynamics of nearby trajectories. To distinguish the periodic orbit from such nearby trajectories, we denote the former by $\tilde{x}(\tau)$. Similarly, let $\tilde{\lambda}_{\mathrm{de},T}(\tau)$ denote the corresponding solution to the adjoint conditions obtained with $\eta_T=1$ and $\eta_{x(0)}=\eta_{x(1)}=0$, such that $\tilde{\lambda}_{\mathrm{de},T}^{\transp}(\tau)=\tilde{\lambda}_{\mathrm{de},T}^{\transp}(\tau+1)$. It follows that the vectors $\tilde{\lambda}_{\mathrm{de},T}^{\transp}(\tau)$ and $f(\tilde{x}(\tau),p)$ are left and right nullvectors of the operator
\begin{equation}
    \Gamma_\tau:= \tilde{X}(\tau+1)\tilde{X}^{-1}(\tau)-I_n,
\end{equation}
where $\tilde{X}(\tau)$ satisfies the variational initial-value problem
\begin{equation}
    X'(\tau)=T\partial_x f(\tilde{x}(\tau),p)X(\tau),\,X(0)=I_n.
\end{equation}
Moreover, since
$\tilde{\lambda}_{\mathrm{de},T}^{\transp}(\tau)f(\tilde{x}(\tau),p)=1$, \jrem{using the notation $(\cdot)_\mathrm{tg}$ for \emph{tangential},}
the definition 
\begin{equation}
    q_\mathrm{tg}(\tau):=f(\tilde{x}(\tau),p)\tilde{\lambda}_{\mathrm{de},T}^{\transp}(\tau)
\end{equation}
defines a projection (since $q_\mathrm{tg}q_\mathrm{tg}=q_\mathrm{tg}$) such that
\begin{equation}
    q_\mathrm{tg}(\tau+1)=q_\mathrm{tg}(\tau),\,\tilde{X}(\tau)q_\mathrm{tg}(0)=q_\mathrm{tg}(\tau)\tilde{X}(\tau),
\end{equation}
and $\Gamma_\tau q_\mathrm{tg}(\tau)=q_\mathrm{tg}(\tau)\Gamma_\tau=0$, i.e., such that $\Gamma_\tau$ maps the image of $q_\mathrm{tr}(\tau):=I_n-q_\mathrm{tg}(\tau)$ to itself and $\Gamma_\tau q_\mathrm{tr}(\tau)=\Gamma_\tau$. \jrem{The notation $(\cdot)_\mathrm{tr}$ indicates that this projection is onto a subspace \emph{transversal} to the tangent of the orbit.} Provided that the right nullspace of $\Gamma_\tau$ is one dimensional, it follows that $\Gamma_\tau$ is invertible on this image. But this is a consequence of the regularity assumption on the monodromy matrix $\tilde{X}(1)$ in Section~\ref{sec:Periodic orbits}, since
\begin{equation}
 \tilde{X}(\tau+1)\tilde{X}^{-1}(\tau)v=v\Leftrightarrow \tilde{X}(1)\tilde{X}^{-1}(\tau)v=\tilde{X}^{-1}(\tau)v
\end{equation}
follows from $\tilde{X}(\tau+1)=\tilde{X}(\tau)\tilde{X}(1)$.

For $x\approx \tilde{x}(\tau)$ for some $\tau$, the scalar equation
\begin{equation}
    \label{eq:localcoords}
    \tilde{\lambda}_{\mathrm{de},T}^{\transp}(\sigma)(x-\tilde{x}(\sigma))=0
\end{equation}
defines $\sigma$ uniquely on a neighborhood of $\tau$, such that $\sigma=\tau$ when $x=\tilde{x}(\tau)$. Indeed, the derivative of the left-hand side of \eqref{eq:localcoords} with respect to $\sigma$, evaluated at $x=\tilde{x}(\tau)$ and $\sigma=\tau$ equals
\begin{equation}
    -T\tilde{\lambda}_{\mathrm{de},T}^{\transp}(\tau)f(\tilde{x}(\tau),p)
\end{equation}
which reduces to the nonzero scalar $-T$. It follows by implicit differentiation of \eqref{eq:localcoords} w.r.t.\ $x$, evaluated at $\sigma=\tau$ and $x=\tilde{x}(\tau)$, that
\begin{equation}
    \partial_x \sigma(\tilde{x}(\tau))=\frac{1}{T}\tilde{\lambda}_{\mathrm{de},T}^{\transp}(\tau).
\end{equation}
Then, for $x=\tilde{x}(\tau)+\delta$ with $\left\|\delta\right\|\ll 1$, it follows that $x=\tilde{x}(\sigma(x))+x_\mathrm{tr}(x)$, where $q_\mathrm{tg}(\sigma(x))x_\mathrm{tr}(x)=0$ \jrem{(i.e., $x_\mathrm{tr}$ is the projection onto the transversal subspace)} and
\begin{align}
    \sigma(x)&=\tau+\frac{1}{T}\tilde{\lambda}_{\mathrm{de},T}^{\transp}(\tau)\delta+O(\left\|\delta\right\|^2),\\
    x_\mathrm{tr}(x)&=q_\mathrm{tr}(\tau)\delta+O(\left\|\delta\right\|^2).
\end{align}

For an arbitrary curve $x(\tau)$ that remains near the limit cycle for all $\tau$, we obtain the unique decomposition
\begin{equation}
\label{eq:assumedform}
    x(\tau)=\tilde{x}(\sigma(\tau))+x_\mathrm{tr}(\tau),
\end{equation}
where $q_\mathrm{tg}(\sigma(\tau))x_\mathrm{tr}(\tau)=0$ for all $\tau$. Suppose, for example, that this is true with $\|x_\mathrm{tr}(\tau)\|=O(\epsilon)$ for a solution $x(\tau)$ of the perturbed differential equation
\begin{equation}
\label{ds:periodicperturbed}
    x'=Tf(x,p)+\delta_\mathrm{ode}(\tau)
\end{equation}
with $\|\delta_\mathrm{ode}(\tau)\|=O(\epsilon)$. Substitution of \eqref{eq:assumedform} into \eqref{ds:periodicperturbed} yields
\begin{equation}
\label{eq:persub}
    Tf(\tilde{x}(\sigma),p)(\sigma'-1)+x'_\mathrm{tr}=T\partial_xf(\tilde{x}(\sigma),p)x_\mathrm{tr}+\delta_\mathrm{ode}+O(\epsilon^2).
\end{equation}
from which it follows that $\sigma'=1+O(\epsilon)$. Multiplication by $q_\mathrm{tr}(\sigma)$ results in
\begin{equation}
\label{eq:subs2}
   q_\mathrm{tr}(\sigma)x_\mathrm{tr}'=Tq_\mathrm{tr}(\sigma)\partial_x f(\tilde{x}(\sigma),p)x_\mathrm{tr}+q_\mathrm{tr}(\sigma)\delta_\mathrm{ode}+O(\epsilon^2).
\end{equation}
But from the form of the rate $\sigma'$, we find
\begin{equation}
    q_\mathrm{tr}(\sigma)x_\mathrm{tr}'=x_\mathrm{tr}'-q'_\mathrm{tr}(\sigma)x_\mathrm{tr}+O(\epsilon^2)
\end{equation}
and, consequently,
\begin{equation}
\label{eq:subs3}
   x_\mathrm{tr}'=Tq_\mathrm{tr}(\sigma)\partial_x f(\tilde{x}(\sigma),p)x_\mathrm{tr}+q'_\mathrm{tr}(\sigma)x_\mathrm{tr}+q_\mathrm{tr}(\sigma)\delta_\mathrm{ode}+O(\epsilon^2).
\end{equation}

Now recall from before that
\begin{equation}
\label{eq:qstVcommuteper}
    \tilde{X}(\sigma)q_\mathrm{tr}(0)=q_\mathrm{tr}(\sigma)\tilde{X}(\sigma).
\end{equation}
Differentiation with respect to $\sigma$ and use of \eqref{eq:qstVcommuteper} then yields
\begin{equation}
    T\partial_x f(\tilde{x}(\sigma),p)q_\mathrm{tr}(\sigma)=q'_\mathrm{tr}(\sigma)+Tq_\mathrm{tr}(\sigma)\partial_x f(\tilde{x}(\sigma),p).
\end{equation}
Multiplication with $x_\mathrm{tr}$ then results in
\begin{equation}
\label{eq:dqstxstcomm}
    T\partial_x f(\tilde{x}(\sigma),p)x_\mathrm{tr}=q'_\mathrm{tr}(\sigma)x_\mathrm{tr}+Tq_\mathrm{tr}(\sigma)\partial_x f(\tilde{x}(\sigma),p)x_\mathrm{tr}%+O(\epsilon^2)
\end{equation}
and, consequently,
\begin{equation}
\label{eq:subs4per}
   x_\mathrm{tr}'=T\partial_x f(\tilde{x}(\sigma),p)x_\mathrm{tr}+q_\mathrm{tr}(\sigma)\delta_\mathrm{ode}+O(\epsilon^2).
\end{equation}
Substitution in \eqref{eq:persub} and multiplication by $\tilde{\lambda}_{\mathrm{de},T}^{\transp}(\sigma)$ finally yields
\begin{equation}\label{eq:phaseper}
    \sigma'=1+\frac{1}{T}\tilde{\lambda}_{\mathrm{de},T}^{\transp}(\sigma)\delta_\mathrm{ode}+O(\epsilon^2).
\end{equation}

We proceed to assume that the eigenvalues of the monodromy matrix $\tilde{X}(1)$ away from $1$ all lie within the unit circle, i.e., that the periodic orbit is normally hyperbolic and orbitally asymptotically stable. Then,
\begin{equation}
    \|\tilde{X}(k)q_\mathrm{tr}(0)\|=\|\tilde{X}(1)-f(\tilde{x}(0),p)\tilde{\lambda}_{\mathrm{de},T}^{\transp}(0)\|^k\le e^{-k/\tau_\mathrm{tr}}
\end{equation}
for some positive constant $\tau_\mathrm{tr}$. It follows from \eqref{eq:subs4per} that if $\|x_\mathrm{tr}(\tau)\|$ is $O(\epsilon)$ at $\tau=0$, then it remains so for all $\tau$ and \eqref{eq:phaseper} remains valid as well for all $\tau$. In this case, for $\delta_\mathrm{ode}(\tau)=0$ and some perturbation $\delta$ to the initial condition $\tilde{x}(\sigma_0)$, it follows that $\sigma(\tau)=\sigma(0)+\tau+O(\delta^2)$, and we find that
\begin{align}
    F\left(T\tau,\tilde{x}(\sigma_0)+\delta,p\right)-F\left(T\tau,\tilde{x}\left(\sigma_0+\frac{1}{T}\tilde{\lambda}_{\mathrm{de},T}^{\transp}(\sigma_0)\delta\right),p\right)
\end{align}
behaves as $O(\|\delta\|^2)+O\left(\exp(-\tau/\tau_\mathrm{tr})\right)$ for large $\tau$. The quantity $\sigma_0+\frac{1}{T}\tilde{\lambda}_{\mathrm{de},T}^{\transp}(\sigma_0)\delta$ is the linear (in $\delta$) approximation to the corresponding \emph{asymptotic phase}.

By the assumption on the eigenvalues of $X(1)$ it further follows that
\begin{equation}
    \tilde{X}(k)=f(\tilde{x}(0),p)\tilde{\lambda}_{\mathrm{de},T}^{\transp}(0)+O\left(\exp(-k/\tau_\mathrm{tr})\right)
\end{equation}
for large $k$. We obtain
\begin{equation}
    \tilde{\lambda}_{\mathrm{de},T}^{\transp}(0)=\lim_{k\rightarrow\infty}\frac{f^{\transp}(\tilde{x}(0),p)}{\|f(\tilde{x}(0),p)\|^2}\tilde{X}(k),
\end{equation}
where the limit has exponential convergence with rate $1/\tau_\mathrm{tr}$.

\subsection{Torus functions}
\label{sec:torus functions}

Let $\mathbb{S}\sim[0,1]$, such that an $\mathbb{R}^n$-valued function on $\mathbb{S}$ is a topological circle. We seek to generalize the treatment in previous sections to the infinite-dimensional boundary-value problem
\begin{equation}
\label{eq:torusbvp}
    \partial_\tau v(\phi,\tau)=Tf(v(\phi,\tau),p),\,v(\phi,1)=v(\phi+\rho,0)
\end{equation}
for the continuously differentiable function $v:\mathbb{S}\times[0,1]\rightarrow\mathbb{R}^n$ in terms of the \emph{a priori} unknown quantities $T$, $\rho$, and $p$. If such a solution exists, then
\begin{equation}
\label{eq:partialvperiodicity}
    \partial_\tau v(\phi,1)=\partial_\tau v(\phi+\rho,0),\,\partial_\phi v(\phi,1)=\partial_\phi v(\phi+\rho,0).
\end{equation}

The boundary condition in \eqref{eq:torusbvp} ensures that the image of $v$ is an invariant topological torus that is covered by a parallel flow described by the rotation number $\rho$. Indeed, given a solution to \eqref{eq:torusbvp}, we may define the \emph{torus function} $u:\mathbb{S}\times\mathbb{S}\rightarrow\mathbb{R}^n$ such that 
\begin{equation}
\label{eq:udef}
u(\theta_1,\theta_2):=v(\theta_1-\rho\theta_2,\theta_2)\mbox{ and  }v(\phi,\tau)=u(\phi+\rho\tau,\tau).
\end{equation}
It follows that $x(\tau)=u(\theta_1(\tau),\theta_2(\tau))$ satisfies $x'=Tf(x,p)$ if $\theta_1'=\rho$ and $\theta_2'=1$. Thus, in terms of the angular coordinates $\theta_1$ and $\theta_2$, the flow on the invariant torus is a rigid rotation. If $\rho$ is irrational, trajectories on the invariant torus are \textit{quasiperiodic} and cover the torus densely.

We may use the torus function $u$ to construct the two-dimensional family
\begin{equation}
    v_{(\theta_1,\theta_2)}:(\phi,\tau)\mapsto u(\theta_1+\phi+\rho\tau,\theta_2+\tau)
    \label{eq:2dimtorusfamily}
\end{equation}
of solutions to \eqref{eq:torusbvp}, defined for $(\phi,\tau)\in\mathbb{S}\times\mathbb{R}$ and
parameterized by the initial condition $u(\theta_1+\phi,\theta_2)$ for arbitrary constants $\theta_1,\theta_2\in\mathbb{S}$. We obtain the original solution in the special case that $\theta_1=\theta_2=0$ and omit subscripts in this case. Using this definition, we obtain $v(\phi,\tau+1)=v(\phi+\rho,\tau)$ for all $(\phi,\tau)\in\mathbb{S}\times\mathbb{R}$, and \eqref{eq:partialvperiodicity} generalizes to
\begin{equation}
\label{eq:partvquasi}
    \partial_\tau v(\phi,\tau+1)=\partial_\tau v(\phi+\rho,\tau),\,\partial_\phi v(\phi,\tau+1)=\partial_\phi v(\phi+\rho,\tau)
\end{equation}
for all $(\phi,\tau)\in\mathbb{S}\times\mathbb{R}$.

Let $V(\phi,\tau)$ be the solution to the variational initial-value problem
\begin{equation}
\label{eq:variationaltorus}
    \partial_\tau V=T\partial_xf(v,p)V,\,V(\cdot,0)=I_n,
\end{equation}
such that 
\begin{equation}
\label{eq:Vprod}
    V(\phi,\tau+1)=V(\phi+\rho,\tau)V(\phi,1)
\end{equation}
follows from \eqref{eq:torusbvp}. Differentiation of both sides of the following equalities with respect to $\tau$ then shows that
\begin{equation}
\label{eq:Xpartphiv}
\partial_\phi v(\phi,\tau)=V(\phi,\tau)\partial_\phi v(\phi,0)
\end{equation}
and 
\begin{equation}
\label{eq:Xparttauv}
    \partial_\tau v(\phi,\tau)=V(\phi,\tau)\partial_\tau v(\phi,0).
\end{equation}
We conclude that $\partial_\tau v(\cdot,\tau)$ and $\partial_\phi v(\cdot,\tau)$ are right nullvectors of the linear operator \begin{equation}\label{def:glin}
    \Gamma_{\rho,\tau}:\delta(\cdot)\mapsto V(\cdot-\rho,\tau+1)V^{-1}(\cdot-\rho,\tau)\delta(\cdot-\rho)-\delta(\cdot).
\end{equation}

The non-uniqueness implied by \eqref{eq:2dimtorusfamily} may be removed by appending two \emph{phase conditions} to \eqref{eq:torusbvp}. In general, these take the form
\begin{align}
h_1(v(\cdot,\cdot),p)=h_2(v(\cdot,\cdot),p)=0\label{eq:torusphase}
\end{align}
in terms of two functionals $h_1(\cdot,p)$ and $h_2(\cdot,p)$ that satisfy a suitable non-degeneracy condition. Here, we assume that these functionals are chosen so that the square matrix
\begin{equation}\label{phase:ndeg:h}
    \int_\mathbb{S}\int_0^1\begin{bmatrix}\partial_vh_1(v(\cdot,\cdot),p)\\\partial_vh_2(v(\cdot,\cdot),p)\end{bmatrix}(\phi,\tau)\begin{bmatrix}\partial_\tau v(\phi,\tau) & \partial_\phi v(\phi,\tau)\end{bmatrix} \mathrm{d}\tau\,\mathrm{d}\phi
\end{equation}
is nonsingular.

\subsection{Normal hyperbolicity}
\label{sec:normal hyperbolicity}
\jrem{We assume henceforth that the solution $v$ of \eqref{eq:torusbvp} is \textit{normally hyperbolic}. For our scenario, this implies the existence of an invariant continuous family of \emph{tangential} projections
\begin{equation}
\label{eq:qtgdef}
    q_\mathrm{tg}(\phi,\tau):=\partial_\tau v(\phi,\tau)q_\tau^{\transp}(\phi,\tau)+\partial_\phi v(\phi,\tau)q_\phi^{\transp}(\phi,\tau)
\end{equation}
onto the tangent spaces of the torus, and a complementary family of \emph{transversal} projections $q_\mathrm{tr}(\phi,\tau):=I_n-q_\mathrm{tg}(\phi,\tau)$, such that
\begin{equation}
\label{tg:qper:inv}
    q_\mathrm{tg}(\phi,\tau+1)=q_\mathrm{tg}(\phi+\rho,\tau),\quad
    q_\mathrm{tg}(\phi,\tau)V(\phi,\tau)=V(\phi,\tau)q_\mathrm{tg}(\phi,0).
\end{equation}
Moreover, the map $\Gamma_{\rho,\tau}$ is a bijection with bounded inverse on the space of transversal perturbations $\phi\mapsto q_\mathrm{tr}(\phi,\tau)\delta(\phi)$ for arbitrary continuous periodic $\delta(\phi)$.} By \eqref{eq:Vprod}, we obtain the conjugacy
\begin{equation}
\label{eq:Gsimilar}
    V^{-1}(\phi,\tau)\Gamma_{\rho,\tau}\left[\delta(\cdot)\right](\phi,\tau)=\Gamma_{\rho,0}\left[V^{-1}(\cdot,\tau)\delta(\cdot)\right](\phi,\tau)
\end{equation}
between the operators $\Gamma_{\rho,\tau}$ and
$\Gamma_{\rho,0}$. Given the continuous families of projections
$q_\mathrm{tg}(\phi,\tau)$ and $q_\mathrm{tr}(\phi,\tau)$, normal
hyperbolicity then follows if and only if the map $\Gamma_{\rho,0}$ is
a bijection with bounded inverse on the space of functions
$\phi\mapsto q_\mathrm{tr}(\phi,0)\delta(\phi)$ for arbitrary
continuous periodic $\delta(\phi)$. In other words, invertibility of $\Gamma_{\rho,\tau}\vert_{\rg q_\mathrm{tr}(\cdot,\tau)}$ for all $\tau$ is equivalent to invertibility of $\Gamma_{\rho,0}\vert_{\rg q_\mathrm{tr}(\cdot,0)}$. General theorems about persistence of normally hyperbolic manifolds \cite{HPS77} imply, for example, that torus functions $v$ with non-rigid rotation (such that $\rho(\phi)$ is non-constant periodic) persist.

Since \jrem{the tangential projection} $q_\mathrm{tg}$ is \jrem{applied pointwise}, $q_\mathrm{tg}q_\mathrm{tg}=q_\mathrm{tg}$ and $q_\mathrm{tg}q_\mathrm{tr}=q_\mathrm{tr}q_\mathrm{tg}=0$ hold everywhere. It follows from \eqref{eq:qtgdef} that
\begin{align}
\label{eq:orthogqtauqphi}
  \begin{bmatrix}
    1 &0\\ 0&1
  \end{bmatrix}=
  &  \begin{bmatrix}
    \smallskip q_\phi^{\transp}(\phi,\tau)\\ q_\tau^{\transp}(\phi,\tau)
  \end{bmatrix}
  \begin{bmatrix}
    \partial_\phi v(\phi,\tau),
    &
    \partial_\tau v(\phi,\tau)
  \end{bmatrix}
\end{align}
for all $(\phi,\tau)\in\mathbb{S}\times\mathbb{R}$. Moreover, by \eqref{eq:partialvperiodicity}, \eqref{eq:Xpartphiv}, \eqref{eq:Xparttauv},  and \eqref{tg:qper:inv}, we obtain
\begin{equation}
\label{eq:perqtauqphi}
    q_\tau(\phi,\tau+1)=q_\tau(\phi+\rho,\tau),\,q_\phi(\phi,\tau+1)=q_\phi(\phi+\rho,\tau),
\end{equation}
and
\begin{equation}
    q_\tau^{\transp}(\phi,\tau)V(\phi,\tau)=q_\tau^{\transp}(\phi,0),\,q_\phi^{\transp}(\phi,\tau)V(\phi,\tau)=q_\phi^{\transp}(\phi,0).
\end{equation}
It follows that
\begin{equation}
\label{eq:shiftqtau}
    0=q_\tau^{\transp}(\phi+\rho,\tau)V(\phi,\tau+1)V^{-1}(\phi,\tau)-q_\tau^{\transp}(\phi,\tau)
\end{equation}
and
\begin{equation}
\label{eq:shiftqphi}
    0=q_\phi^{\transp}(\phi+\rho,\tau)V(\phi,\tau+1)V^{-1}(\phi,\tau)-q_\phi^{\transp}(\phi,\tau).    
\end{equation}
Multiplication of each of these equalities by an arbitrary function $\delta(\phi)$ and integration over $\mathbb{S}$ then yields
\begin{align}
    0&=\int_\mathbb{S}\left(q_\tau^{\transp}(\phi+\rho,\tau)V(\phi,\tau+1)V^{-1}(\phi,\tau)-q_\tau^{\transp}(\phi,\tau)\right)\delta(\phi)\,\mathrm{d}\phi\nonumber\\
    &=\int_\mathbb{S}q_\tau^{\transp}(\phi,\tau)\left(V(\phi-\rho,\tau+1)V^{-1}(\phi-\rho,\tau)\delta(\phi-\rho)-\delta(\phi)\right)\,\mathrm{d}\phi\nonumber\\
    &=\int_\mathbb{S}q_\tau^{\transp}(\phi,\tau)\Gamma_{\rho,\tau}\left[\delta(\cdot)\right](\phi)\,\mathrm{d}\phi
\end{align}
and
\begin{align}
    0&=\int_\mathbb{S}\left(q_\phi^{\transp}(\phi+\rho,0)V(\phi,\tau+1)V^{-1}(\phi,\tau)-q_\phi^{\transp}(\phi,0)\right)\delta(\phi)\,\mathrm{d}\phi\nonumber\\
    &=\int_\mathbb{S}q_\phi^{\transp}(\phi,\tau)\left(V(\phi-\rho,\tau+1)V^{-1}(\phi-\rho,\tau)\delta(\phi-\rho)-\delta(\phi)\right)\,\mathrm{d}\phi\nonumber\\
    &=\int_\mathbb{S}q_\phi^{\transp}(\phi,\tau)\Gamma_{\rho,\tau}\left[\delta(\cdot)\right](\phi)\,\mathrm{d}\phi,
\end{align}
i.e., that the linear functionals $\int_\mathbb{S}q_\tau^{\transp}(\phi,\tau)\left(\cdot\right)\,\mathrm{d}\phi$ and $\int_\mathbb{S}q_\phi^{\transp}(\phi,\tau)\left(\cdot\right)\,\mathrm{d}\phi$ lie in the left nullspace of the operator $\Gamma_{\rho,\tau}$.

With the help of \eqref{eq:partvquasi}, \eqref{eq:Vprod}, \eqref{eq:Xpartphiv}, \eqref{eq:Xparttauv},  and \eqref{def:glin}, it follows from \eqref{eq:qtgdef} and \eqref{tg:qper:inv} that
\begin{align}
    &q_\mathrm{tg}(\phi,\tau)\Gamma_{\rho,\tau}\left[\delta(\cdot)\right](\phi,\tau)=\Gamma_{\rho,\tau}\left[q_\mathrm{tg}(\cdot,\tau)\delta(\cdot)\right](\phi,\tau)=-q_\mathrm{tg}(\phi,\tau)\delta(\phi)\nonumber\\
    &\qquad+\left(\partial_\tau v(\phi,\tau)q_\tau^{\transp}(\phi-\rho,\tau)+\partial_\phi v(\phi,\tau)q_\phi^{\transp}(\phi-\rho,\tau)\right)\delta(\phi-\rho).
\end{align}
Then, \eqref{eq:orthogqtauqphi} implies that
\begin{equation}
    q_\mathrm{tg}(\phi,\tau)\Gamma_{\rho,\tau}\left[q_\mathrm{tr}(\cdot,\tau)\delta(\cdot)\right](\phi,\tau)=0,
\end{equation}
i.e., that the space of functions  $\phi\mapsto q_\mathrm{tr}(\phi,\tau)\delta(\phi)$ for arbitrary continuous periodic functions $\delta(\phi)$ is, in fact, invariant under $\Gamma_{\rho,\tau}$. That $\Gamma_{\rho,\tau}$ is a bijection on this space with a bounded inverse is then equivalent to the existence of a bounded inverse of the map $\hat{\Gamma}_{\rho,\tau}$ given by
\begin{equation}
    \label{nhyp:gst:def}
  \hat{\Gamma}_{\rho,\tau}\left[\delta(\cdot)\right](\phi,\tau):=\Gamma_{\rho,\tau}\left[ q_\mathrm{tr}(\cdot,\tau)\delta(\cdot)\right](\phi,\tau)-q_\mathrm{tg}(\phi,\tau)\delta(\phi)
\end{equation}
on the space of continuous periodic functions $\delta(\phi)$. From \eqref{eq:Gsimilar}, we obtain
\begin{equation}
  V^{-1}(\phi,\tau)\hat{\Gamma}_{\rho,\tau}\left[\delta(\cdot)\right](\phi,\tau)=\hat{\Gamma}_{\rho,0}\left[ V^{-1}(\cdot,\tau)\delta(\cdot)\right](\phi,\tau),
\end{equation}
i.e., that it again suffices to study the restriction to the case $\tau=0$. In particular, normal hyperbolicity implies the existence of a unique continuous periodic solution $\delta=\hat{\Gamma}_{\rho,0}^{-1}\left[\delta_\mathrm{rhs}(\cdot)\right]$ for continuous periodic $\delta_\mathrm{rhs}$ (with the norm of $\delta$ bounded by a fixed multiple of the norm of $\delta_\mathrm{rhs}$) such that
\begin{equation}
\label{torus:transversal:hyp}
    V(\phi-\rho,1)q_\mathrm{tr}(\phi-\rho,0)\delta(\phi-\rho)-\delta(\phi)=\delta_\mathrm{rhs}(\phi).
\end{equation}

Suppose, for example, that $\delta(\cdot)$ is an eigenfunction of $\hat{\Gamma}_{\rho,0}$ with eigenvalue $z$, i.e., such that
\begin{equation}
  V(\phi-\rho,1)q_\mathrm{tr}(\phi-\rho,0)\delta(\phi-\rho)=(1+z)\delta(\phi)
\end{equation}
for all $\phi$. Then, by induction and liberal use of \eqref{eq:Vprod} and \eqref{tg:qper:inv},
\begin{equation}
    V(\phi-k\rho,k)q_\mathrm{tr}(\phi-k\rho,0)\delta(\phi-k\rho)=(1+z)^k\delta(\phi)
\end{equation}
and, consequently,
\begin{equation}
    |1+z|^k\le\|V(\cdot,k)q_\mathrm{tr}(\cdot,0)\|.
\end{equation}
If, in addition to $v$ being normally hyperbolic, the sequence $V(\cdot,k)q_\mathrm{tr}(\cdot,0)$  is bounded by 
\begin{align}
  \label{eq:normdecay}
  \|V(\cdot,k)q_\mathrm{tr}(\cdot,0)\|\leq C_\mathrm{tr}\exp(-k/\tau_\mathrm{tr})
\end{align}
for some positive constants $C_\mathrm{tr}$ and $\tau_\mathrm{tr}$, we say that $v$ is \emph{transversally stable}. In this case
\begin{equation}
    |1+z|\le \exp(-1/\tau_\mathrm{tr})
\end{equation}
and the spectral radius of $\hat{\Gamma}_{\rho,0}^{-1}$ must be bounded by $1/(1-\exp(-1/\tau_\mathrm{tr}))$.

In general,
\begin{equation}
    V(\phi,k)q_\mathrm{tg}(\phi,0)=\partial_\tau v(\phi+k\rho,0)q_\tau^{\transp}(\phi,0)+\partial_\phi v(\phi+k\rho,0)q_\phi^{\transp}(\phi,0).
\end{equation}
For a transversally stable $v$, it follows from $q_\mathrm{tr}=I_n-q_\mathrm{tg}$ that
\begin{equation}
\label{qtg:limk}
    V(\phi,k)=\partial_\tau v(\phi+k\rho,0)q_\tau^{\transp}(\phi,0)+\partial_\phi v(\phi+k\rho,0)q_\phi^{\transp}(\phi,0)+O(\exp(-k/\tau_\mathrm{tr})).
\end{equation}
Provided that the $2\times2$ matrix
\begin{align}
  A(\phi):=
  \begin{bmatrix} 
   \partial_\phi v^{\transp}(\phi,0)\partial_\phi v(\phi,0)&
   \partial_\phi v^{\transp}(\phi,0)\partial_\tau v(\phi,0)\\
   \partial_\tau v^{\transp}(\phi,0)\partial_\phi v(\phi,0)&
   \partial_\tau v^{\transp}(\phi,0)\partial_\tau v(\phi,0)
  \end{bmatrix}
\end{align}
has a uniformly bounded inverse (in $\phi$), we may multiply \eqref{qtg:limk} by $\partial_\phi v^{\transp}(\phi+k\rho,0)$ and $\partial_{\tau} v^{\transp}(\phi+k\rho,0)$, respectively, to obtain
\begin{align}
  \begin{bmatrix}
    \partial_\phi v^{\transp}(\phi+k\rho,0)\\
    \partial_{\tau} v^{\transp}(\phi+k\rho,0)
  \end{bmatrix}
  V(\phi,k)&=A(\phi+k\rho)
  \begin{bmatrix}
    q_\phi^{\transp}(\phi,0)\\[0.5ex]q_\tau^{\transp}(\phi,0)
  \end{bmatrix}+O(\exp(-k/\tau_\mathrm{tr})).
\end{align}
Consequently,
\begin{align}
  q_\mathrm{tg}(\phi,0)=\begin{bmatrix}
    v_\phi(\phi,0)&v_\tau(\phi,0)
  \end{bmatrix}\lim_{k\to\infty}A(\phi+k\rho)^{-1}
  \begin{bmatrix}
    \partial_\phi v^{\transp}(\phi+k\rho,0)\\
    \partial_{\tau} v^{\transp}(\phi+k\rho,0)
  \end{bmatrix}
V(\phi,k),
\end{align}
where the limit has exponential convergence with rate $1/\tau_\mathrm{tr}$.

\subsection{Adjoint conditions}\label{sec:torus:adj}
Inspired by the analysis in previous sections, we next consider the Lagrangian
\begin{align}
    L&=\int_\mathbb{S}\int_0^1\lambda_\mathrm{de}^{\transp}(\phi,\tau)(\partial_\tau v(\phi,\tau)-Tf(v(\phi,\tau),p))\,\mathrm{d}\tau\,\mathrm{d}\phi\nonumber\\
    &\qquad+\int_\mathbb{S} \lambda_\mathrm{bc}^{\transp}(\phi)(v(\phi+\rho,0)-v(\phi,1))\,\mathrm{d}\phi+\lambda_\mathrm{ps}^{\transp}h(v(\cdot,\cdot),p)\nonumber\\
    &\qquad\qquad+\eta_\rho(\rho-\mu_\rho)+\eta_T(T-\mu_T)+\eta^{\transp}_p(p-\mu_p),
\end{align}
from which we derive the adjoint conditions
\begin{align}\label{phase:var:dv}
  &\mbox{$\delta v(\phi,t)$:}&0&=-\partial_\tau \lambda_\mathrm{de}^{\transp}-T\lambda_\mathrm{de}^{\transp}\partial_xf(v,p)+\lambda_\mathrm{ps}^{\transp}\partial_v h(v(\cdot,\cdot),p),\\
\label{phase:var:dv1}
  &\mbox{$\delta v(\phi,1)$:}&    0&=\phantom{-}\lambda_\mathrm{de}^{\transp}(\phi,1)-\lambda_\mathrm{bc}^{\transp}(\phi),\\
\label{phase:var:dv0}
  &\mbox{$\delta v(\phi,0)$:}&    0&=-\lambda_\mathrm{de}^{\transp}(\phi,0)+\lambda_\mathrm{bc}^{\transp}(\phi-\rho),\\
\label{phase:var:dT}
  &\mbox{$\delta T$:}&     0&=-\int_\mathbb{S}\int_0^1\lambda_\mathrm{de}^{\transp}(\phi,\tau)f(v(\phi,\tau),p)\,\mathrm{d}\tau\,\mathrm{d}\phi+\eta_T,\\
\label{phase:var:drho}
  &\mbox{$\delta\rho$:}&    0&=\int_\mathbb{S}\lambda_\mathrm{bc}^{\transp}(\phi)\partial_\phi v(\phi+\rho,0)\,\mathrm{d}\phi+\eta_\rho,\\
\label{phase:var:dp}
  &\mbox{$\delta p$:}&    0&=-\int_\mathbb{S}\int_0^1\lambda_\mathrm{de}^{\transp}(\phi,\tau)T\partial_pf(v(\phi,\tau),p)\,\mathrm{d}\tau\,\mathrm{d}\phi+\eta_p^{\transp}\nonumber\\
  &&&\qquad+\lambda_\mathrm{ps}^{\transp}\partial_ph(v(\cdot,\cdot),p).
\end{align}
For every solution to these conditions, by
elimination of $\lambda_\mathrm{bc}$ from \eqref{phase:var:dv1} and
\eqref{phase:var:dv0}, it must hold that
\begin{equation}\label{phase:lde:periodic}
    \lambda_\mathrm{de}(\phi,1)=\lambda_\mathrm{de}(\phi+\rho,0)
\end{equation}
analogously to \eqref{eq:torusbvp}. Multiplying \eqref{phase:var:dv} by the solution $V(\phi,\tau)$ to
\eqref{eq:variationaltorus} yields
\begin{equation}
    \lambda_\mathrm{ps}^{\transp}\partial_vh(v(\cdot,\cdot),p)V(\phi,\tau)=\partial_\tau\left(\lambda_\mathrm{de}^{\transp}(\phi,\tau)V(\phi,\tau)\right)
\end{equation}
and integrating with respect to $\tau$ from $0$ to $1$ then results in
\begin{equation}
\label{eq:lambdapscond}
    \lambda_\mathrm{ps}^{\transp}\int_0^1\partial_vh(v(\cdot,\cdot),p)V(\phi,\tau)\,\mathrm{d}\tau=\lambda_\mathrm{de}^{\transp}(\phi+\rho,0)V(\phi,1)-\lambda_\mathrm{de}^{\transp}(\phi,0),
\end{equation}
since $V(\phi,0)=I_n$. Multiplication of both sides of \eqref{eq:lambdapscond} by $\partial_\tau v(\phi,0)$ or $\partial_\phi v(\phi,0)$, followed by integration over $\phi$ now yields
\begin{align}
    &\lambda_\mathrm{ps}^{\transp}\int_\mathbb{S}\int_0^1\partial_vh(v(\cdot,\cdot),p)V(\phi,\tau)\partial_\tau v(\phi,0)\,\mathrm{d}\tau\mathrm{d}\phi\nonumber\\
    &\qquad=\int_\mathbb{S}\lambda_\mathrm{de}^{\transp}(\phi,0)\left(V(\phi-\rho,1)\partial_\tau v(\phi-\rho,0)-\partial_\tau v(\phi,0)\right)\,\mathrm{d}\phi\nonumber\\
    &\qquad=\int_\mathbb{S}\lambda_\mathrm{de}^{\transp}(\phi,0)\Gamma_{\rho,0}\left[\partial_\tau v(\cdot,0)\right](\phi)\,\mathrm{d}\phi=0\label{eq:lambdaps1}
\end{align}
    and
\begin{align}
    &\lambda_\mathrm{ps}^{\transp}\int_\mathbb{S}\int_0^1\partial_vh(v(\cdot,\cdot),p)V(\phi,\tau)\partial_\phi v(\phi,0)\,\mathrm{d}\tau\mathrm{d}\phi\nonumber\\
    &\qquad=\int_\mathbb{S}\lambda_\mathrm{de}^{\transp}(\phi,0)\left(V(\phi-\rho,1)\partial_\phi v(\phi-\rho,0)-\partial_\phi v(\phi,0)\right)\,\mathrm{d}\phi\nonumber\\
    &\qquad=\int_\mathbb{S}\lambda_\mathrm{de}^{\transp}(\phi,0)\Gamma_{\rho,0}\left[\partial_\phi v(\cdot,0)\right](\phi)\,\mathrm{d}\phi=0,\label{eq:lambdaps2}
\end{align}
where the final equalities follow from the observation at the end of the previous section that
$\partial_\tau v(\cdot,0)$ and $\partial_\phi v(\cdot,0)$ are right
nullvectors of $\Gamma_{\rho,0}$. 

By the nonsingularity of \eqref{phase:ndeg:h}, the equalities \eqref{eq:lambdaps1} and \eqref{eq:lambdaps2} imply that $\lambda_\mathrm{ps}=0$ and, consequently, that $\lambda_\mathrm{de}^{\transp}(\phi,\tau)V(\phi,\tau)$ must be constant. In particular, by \eqref{eq:Vprod},
\begin{equation}
\label{eq:deXper}
    \lambda_\mathrm{de}^{\transp}(\phi+\rho,\tau)V(\phi,\tau+1)V^{-1}(\phi,\tau)-\lambda_\mathrm{de}^{\transp}(\phi,\tau)=0.
\end{equation}
Multiplication of this equality by an arbitrary function $\delta(\phi)$ and integration over $\mathbb{S}$ then yields the equality
\begin{align}
    \label{phase:lde:left}
    0&=\int_\mathbb{S}\lambda_\mathrm{de}^{\transp}(\phi,\tau)\left(V(\phi-\rho,\tau+1)V^{-1}(\phi-\rho,\tau)\delta(\phi-\rho)-\delta(\phi)\right)\,\mathrm{d}\phi\nonumber\\&\qquad=\int_\mathbb{S}\lambda_\mathrm{de}^{\transp}(\phi,\tau)\Gamma_{\rho,\tau}\left[\delta(\cdot)\right](\phi)\,\mathrm{d}\phi,
\end{align}
i.e., that the linear functional
$\int_\mathbb{S}\lambda_\mathrm{de}^{\transp}(\phi,\tau)\left(\cdot\right)\,\mathrm{d}\phi$
must lie in the left nullspace of the operator $\Gamma_{\rho,\tau}$. Conversely, for every periodic function
$\lambda_\mathrm{bc}(\cdot)$ satisfying
\begin{equation}
\label{eq:lambdabccond}
    \lambda_\mathrm{bc}^{\transp}(\phi)V(\phi,1)-\lambda^{\transp}_\mathrm{bc}(\phi-\rho)=0,
\end{equation}
such that $\int_\mathbb{S}\lambda_\mathrm{bc}^{\transp}(\phi-\rho)\left(\cdot\right)\,\mathrm{d}\phi$
is in the left nullspace of the operator $\Gamma_{\rho,0}$, the function
\begin{equation}
    \lambda_\mathrm{de}^{\transp}(\phi,\tau):=\lambda_\mathrm{bc}^{\transp}(\phi)V(\phi,1)V^{-1}(\phi,\tau)
\end{equation}
satisfies the adjoint boundary value problem
\eqref{phase:var:dv}--\eqref{phase:var:dv0} with
$\lambda_\mathrm{ps}=0$.

It follows from \eqref{eq:shiftqtau} and \eqref{eq:shiftqphi} that the assignment $\lambda_\mathrm{bc}(\phi)=c_\tau q_\tau(\phi+\rho,0)+c_\phi q_\phi(\phi+\rho,0)$ for arbitrary constants $c_\tau$ and $c_\phi$ satisfies \eqref{eq:lambdabccond} and that the corresponding $\lambda_\mathrm{de}^{\transp}(\phi,\tau)$ equals $c_\tau q_\tau^{\transp}(\phi,\tau)+c_\phi q_\phi^{\transp}(\phi,\tau)$. Substitution in \eqref{phase:var:dT} and \eqref{phase:var:drho} then yields $(\eta_T,\eta_\rho)=(c_\tau/T,-c_\phi)$. In the context of the general theory, the solution to the adjoint equations, if it exists, depends on the choice of assignments of $0$'s and $1$'s to the elements of $\eta_{\mathbb{J}_2}$. Here, we focus on two such possibilities, namely that with $\eta_T=1$ and $\eta_\rho=0$, which we denote by an additional subscript $_T$, and that with $\eta_T=0$ and $\eta_\rho=1$, which we denote by the additional subscript $_\rho$. We thus obtain one such pair of solutions $\lambda_{\mathrm{de},T}(\phi,\tau)=Tq_\tau(\phi,\tau)$ and $\lambda_{\mathrm{de},\rho}(\phi,\tau)=-q_\phi(\phi,\tau)$. The corresponding values of $\eta_p^{\transp}$ equal
\begin{equation}
\label{eq: etaT_T}
    T\int_\mathbb{S}q_\tau^{\transp}(\phi,1)P(\phi,1)\,\mathrm{d}\phi
\end{equation}
and
\begin{equation}
\label{eq: etaT_rho}
    -\int_\mathbb{S}q_\phi^{\transp}(\phi,1)P(\phi,1)\,\mathrm{d}\phi,
\end{equation}
respectively, where
\begin{equation}
\label{eq:Pquasi}
    P(\phi,\tau):=TV(\phi,\tau)\int_0^\tau V(\phi,\sigma)^{-1}\partial_pf(v(\phi,\sigma),p)\,\mathrm{d}\sigma\mbox{.}
\end{equation}
In general, since
$\lambda_\mathrm{de}^{\transp}(\phi,\tau)V(\phi,\tau)$ must be
constant in $\tau$, it follows from \eqref{eq:Xpartphiv} and
\eqref{eq:Xparttauv} that the products
\begin{equation}
\lambda_\mathrm{de}^{\transp}(\phi,\tau)\partial_\tau v(\phi,\tau)\mbox{ and }\lambda_\mathrm{de}^{\transp}(\phi,\tau)\partial_\phi v(\phi,\tau)    
\end{equation}
must also be constant in $\tau$. Thus, if the rotation number $\rho$ is irrational (such
that every trajectory on the torus covers it densely) and
$\lambda_\mathrm{bc}$ (and, hence, $\lambda_\mathrm{de}$) is
continuous in $\phi$, the scaling conditions \eqref{phase:var:dT} and
\eqref{phase:var:drho} imply
\begin{align}
\label{eq:orthoggen}
    \lambda_\mathrm{de}^{\transp}(\phi,\tau)\partial_\tau v(\phi,\tau)
  = T\eta_T,\,
  \lambda_\mathrm{de}^{\transp}(\phi,\tau)\partial_\phi v(\phi,\tau)&=-\eta_\rho
\end{align}
for all $(\phi,\tau)$. As an immediate consequence,
\begin{equation}
\label{eq:lambdaqtg}
    \lambda_\mathrm{de}^{\transp}(\phi,\tau)q_\mathrm{tg}(\phi,\tau)=T\eta_Tq_\tau^{\transp}(\phi,\tau)-\eta_\rho q_\phi^{\transp}(\phi,\tau).
\end{equation}
Furthermore, from \eqref{eq:deXper} we find
\begin{align}
\label{eq:lambdadeqst1}
    \lambda_\mathrm{de}^{\transp}(\phi+\rho,\tau)q_\mathrm{tr}(\phi+\rho,\tau)V(\phi,\tau+1)V^{-1}(\phi,\tau)-\lambda_\mathrm{de}^{\transp}(\phi,\tau)q_\mathrm{tr}(\phi,\tau)=0
\end{align}
and, consequently, that the functional $\int_\mathbb{S}\lambda_\mathrm{de}^{\transp}(\phi,\tau)q_\mathrm{tr}(\phi,\tau)(\cdot)\,\mathrm{d}\phi$ must lie in the left nullspace of the operator $\hat{\Gamma}_{\rho,\tau}$. Since $\hat{\Gamma}_{\rho,\tau}$ has a bounded inverse, we conclude that
\begin{equation}
\label{phase:lde:st0}
    \lambda_\mathrm{de}^{\transp}(\phi,\tau)q_\mathrm{tr}(\phi,\tau)=0\mbox{}
\end{equation}
must hold for all $(\phi,\tau)$. 

Considering the
function\footnote{The indicator function
  $\mathbbm{1}_r:\mathbb{S}\rightarrow\{0,1\}$ is nonzero on
  $|\phi|<r/2$ (appropriately defined in the metric on $\mathbb{S}$).}
\begin{equation}
  \delta(\phi)=\mathbbm{1}_r(\phi-\phi_0)\bigl[q_\mathrm{tr}(\phi,\tau)+q_\mathrm{tg}(\phi,\tau)\bigr]\delta_0,
\end{equation}
which equals $\delta_0$ on the ball $B_r$ with diameter $r$ around
some arbitrary $\phi_0$ in $\mathbb{S}$ and $0$ elsewhere,
\eqref{eq:lambdaqtg} and \eqref{phase:lde:st0} imply that
\begin{align}
  \int_\mathbb{S}\lambda_\mathrm{de}^{\transp}(\phi,\tau)\delta(\phi)\,\mathrm{d}\phi
  &=\int_{B_r}\lambda_\mathrm{de}^{\transp}(\phi,\tau)\,\mathrm{d}\phi \,\delta_0\nonumber\\
  &=\left[T\eta_T\int_{B_r}q_\tau^{\transp}(\phi,\tau)\,\mathrm{d}\phi-\eta_\rho\int_{B_r}q_\phi^{\transp}(\phi,\tau)\,\mathrm{d}\phi\right]\delta_0.
\end{align}
Dividing by $r$ and assuming continuity of $\lambda_\mathrm{de}$ in $\phi$, we obtain
\begin{align}
\label{eq:adjointsol}
\lambda_\mathrm{de}(\phi,\tau)=T\eta_Tq_\tau(\phi,\tau)-\eta_\rho q_\phi(\phi,\tau),
\end{align}
since the above integral equality holds for all $\delta_0$, arbitrary
small $r>0$, and arbitrary $\phi_0$. We conclude that for irrational $\rho$, the functions $Tq_\tau(\phi,\tau)$ and $q_\phi(\phi,\tau)$ are, in fact, the unique continuous solutions for $\lambda_{\mathrm{de},T}(\phi,\tau)$ and $\lambda_{\mathrm{de},\rho}(\phi,\tau)$, respectively. It follows from \eqref{eq:orthogqtauqphi} that
\begin{equation}
\label{eq:ortho1}
    \lambda_{\mathrm{de},T}^{\transp}(\phi,\tau)\partial_\tau v(\phi,\tau)=T,\,\lambda_{\mathrm{de},\rho}^{\transp}(\phi,\tau)\partial_\phi v(\phi,\tau)=-1,
\end{equation}
and
\begin{equation}
\label{eq:ortho2}
    \lambda_{\mathrm{de},T}^{\transp}(\phi,\tau)\partial_\phi v(\phi,\tau)=\lambda_{\mathrm{de},\rho}^{\transp}(\phi,\tau)\partial_\tau v(\phi,\tau)=0
\end{equation}
for all $(\phi,\tau)$, consistent with \eqref{eq:orthoggen}.

\subsection{Asymptotic analysis}

As in the analysis of the periodic orbit in Section~\ref{sec:po:revisit}, for $x\approx v(\phi,\tau)$ for some $\phi$ and $\tau$, the scalar equations
\begin{equation}
\label{eq:toruslocalcoords}
    \lambda^{\transp}_{\mathrm{de},T}(\varphi,\sigma)(x-v(\varphi,\sigma))=0,\,\lambda^{\transp}_{\mathrm{de},\rho}(\varphi,\sigma)(x-v(\varphi,\sigma))=0
\end{equation}
define $\varphi$ and $\sigma$ uniquely on a neighborhood of $\phi$ and $\tau$, such that $\varphi=\phi$ and $\sigma=\tau$ when $x=v(\phi,\tau)$. Indeed, the Jacobian of the left-hand sides of \eqref{eq:toruslocalcoords} with respect to $\varphi$ and $\sigma$, evaluated at $x=v(\phi,\tau)$, $\varphi=\phi$, and $\sigma=\tau$ equals
\begin{equation}
    \begin{bmatrix}-\lambda^{\transp}_{\mathrm{de},T}(\phi,\tau)\partial_\phi v(\phi,\tau) & -\lambda^{\transp}_{\mathrm{de},T}(\phi,\tau)\partial_\tau v(\phi,\tau)\\-\lambda^{\transp}_{\mathrm{de},\rho}(\phi,\tau)\partial_\phi v(\phi,\tau) & -\lambda^{\transp}_{\mathrm{de},\rho}(\phi,\tau)\partial_\tau v(\phi,\tau)\end{bmatrix}
\end{equation}
which reduces to the nonsingular matrix
\begin{equation}
    \begin{bmatrix}0 & -T\\1 & 0\end{bmatrix}
\end{equation}
by \eqref{eq:ortho1} and \eqref{eq:ortho2}. It follows by implicit differentiation of \eqref{eq:toruslocalcoords} with respect to $x$, evaluated at $\varphi=\phi$, $\sigma=\tau$, and $x=v(\phi,\tau)$, that
\begin{equation}
    \partial_x\varphi(v(\phi,\tau))=-\lambda^{\transp}_{\mathrm{de},\rho}(\phi,\tau),\,\partial_x\sigma(v(\phi,\tau))=\frac{1}{T}\lambda^{\transp}_{\mathrm{de},T}(\phi,\tau).
\end{equation}
If we define $x_\mathrm{tr}:=x-v(\varphi(x),\sigma(x))$, such that $q_\mathrm{tg}(\varphi(x),\sigma(x))x_\mathrm{tr}(x)=0$, then for $x=v(\phi,\tau)+\delta$ with $\|\delta\|\ll 1$, it follows that $x=v(\varphi(x),\sigma(x))+x_\mathrm{tr}(x)$, where
\begin{align}
    \varphi(x)=\phi-\lambda^{\transp}_{\mathrm{de},\rho}(\phi,\tau)\delta+O(\|\delta\|)^2,\\
    \sigma(x)=\tau+\frac{1}{T}\lambda^{\transp}_{\mathrm{de},T}(\phi,\tau)\delta+O(\|\delta\|)^2.
\end{align}

For an arbitrary curve $x(\tau)$ that remains near the invariant torus for all $\tau$, we obtain the unique decomposition
\begin{align}\label{forcing:decompose}
  x(\tau)=v(\varphi(\tau),\sigma(\tau))+x_\mathrm{tr}(\tau),
\end{align}
where $q_\mathrm{tg}(\varphi(\tau),\sigma(\tau))x_\mathrm{tr}(\tau)=0$ for all $\tau$. Suppose, for example, that this is true with $\|x_\mathrm{tr}(\tau)\|=O(\epsilon)$ for a solution $x(\tau)$ of the perturbed differential equation
\begin{align}
  \label{ds:perturbed}
  x'&=Tf(x,p)+\delta_\mathrm{ode}(\tau),
\end{align}
where $\delta_\mathrm{ode}(\tau)=O(\epsilon)$ for all $\tau$. Substitution of \eqref{forcing:decompose} into \eqref{ds:perturbed} and use of \eqref{eq:torusbvp} yields
\begin{equation}
\label{eq:subs1}
    \partial_\phi v\varphi'+\partial_\sigma v(\sigma'-1)+x_\mathrm{tr}'=T\partial_x f(v,p)x_\mathrm{tr}+\delta_\mathrm{ode}+O(\epsilon^2).
\end{equation}
We conclude that $\varphi'=O(\epsilon)$ and $\sigma'=1+O(\epsilon)$. Following the same arguments for $x_\mathrm{tr}$ as in Section~\ref{sec:po:revisit}, we obtain
\begin{equation}
\label{eq:subs4}
   x_\mathrm{tr}'=T\partial_x f(v,p)q_\mathrm{tr}x_\mathrm{tr}+q_\mathrm{tr}\delta_\mathrm{ode}+O(\epsilon^2).
 \end{equation}
Substitution in \eqref{eq:subs1} and multiplication by
$\lambda_{\mathrm{de},T}^{\transp}(\varphi,\sigma)$ and $\lambda_{\mathrm{de},\rho}^{\transp}(\varphi,\sigma)$,
respectively, then yields the generalization of \eqref{eq:phaseper} for the torus:
\begin{align}
    \varphi'=-\lambda_{\mathrm{de},\rho}^{\transp}(\varphi,\sigma)\delta_\mathrm{ode}+O(\epsilon^2),\label{eq:varphiprime}\\
    \sigma'=1+\lambda_{\mathrm{de},T}^{\transp}(\varphi,\sigma)\frac{\delta_\mathrm{ode}}{T}+O(\epsilon^2).\label{eq:sigmaprime}
\end{align}

Finally, we revisit the assumed uniform $O(\epsilon)$ bound on the magnitude of $x_\mathrm{tr}$. As in Section~\ref{sec:po:revisit}, this follows provided that $v$ is transversally stable, in which case \eqref{eq:varphiprime} and \eqref{eq:sigmaprime} hold for all $\tau$. In contrast to the resultant bounded behavior of $x_\mathrm{tr}$, the variables $\sigma(\tau)-\tau$
and $\phi(\tau)$ may have non-trivial dynamics, for example, determining potential
locking near the invariant torus graph $v(\cdot,\cdot)$. For the special case of $\delta_\mathrm{ode}=0$ and a perturbation $\delta$ to the initial condition $v(\phi_0,\sigma_0)$, transversal stability implies that
\begin{equation}
    \varphi(\tau)=\varphi(0)+O(\|\delta\|^2),\,\sigma(\tau)=\sigma(0)+\tau+O(\|\delta\|^2),
\end{equation}
and we find that
\begin{multline}\label{phase:linear}
F(\tau,v(\phi_0,\sigma_0)+\delta,p)\\-F\left(\tau,v\left(\phi_0-\lambda_{\mathrm{de},\rho}^{\transp}(\phi_0,\sigma_0)\delta,\sigma_0+\frac{1}{T}\lambda_{\mathrm{de},T}^{\transp}(\phi_0,\sigma_0)\delta\right),p\right)
\end{multline}
behaves as $O(\|\delta\|^2)+O(\exp(-\tau/\tau_\mathrm{tr}))$ for large $\tau$. The quantities $\phi_0-\lambda_{\mathrm{de},\rho}^{\transp}(\phi_0,\sigma_0)\delta$ and $\sigma_0+\frac{1}{T}\lambda_{\mathrm{de},T}^{\transp}(\phi_0,\sigma_0)\delta$ are the linear (in $\delta$) approximations to the corresponding asymptotic phases.

\subsection{Regularity}
\label{sec:regularity}
We conclude the general discussion of the quasiperiodic invariant torus by exploring the extent to which the predictions of the adjoint analysis agree with those obtained from the linearization of the zero problem \eqref{eq:torusbvp} and \eqref{eq:torusphase}. We recall from the general theory in Section~\ref{sec:preliminaries} that the values of the adjoint variables $\lambda$ and $\eta_{\mathbb{I}\,\cup\,\mathbb{J}_1}$ capture the sensitivities of the monitor function $\Psi_k(u)$ to violations of the zero problem $\Phi(u)=0$ and variations in $\Psi_{\mathbb{I}\,\cup\,\mathbb{J}_1}(u)$, respectively, provided that the reduced continuation problem
\begin{equation}
    \Phi(u)=0,\,\Psi_{\mathbb{I}\,\cup\,\mathbb{J}_1}(u)-\Psi_{\mathbb{I}\,\cup\,\mathbb{J}_1}(\tilde{u})=0
\end{equation}
is regular at $u=\tilde{u}$ with zero dimensional deficit. In this case,
\begin{equation}
    \begin{pmatrix}\tilde{\lambda} & \tilde{\eta}_{\mathbb{I}\,\cup\,\mathbb{J}_1}\end{pmatrix}=-D\Psi_k(\tilde{u})\begin{pmatrix}D\Phi(\tilde{u})\\D\Psi_{\mathbb{I}\,\cup\,\mathbb{J}_1}(\tilde{u})\end{pmatrix}^{-1},\label{eq:lambdasol}
\end{equation}
where the inverse on the right is assumed to be bounded. As we show below, this assumption fails for the zero problem given by \eqref{eq:torusbvp} and \eqref{eq:torusphase} with irrational $\rho$. Nevertheless, the solutions obtained for $\lambda_\mathrm{de}$, $\lambda_\mathrm{bc}$, $\lambda_\mathrm{ps}$, and $\eta_p$ (whose existence and uniqueness within the space of continuous $\lambda_\mathrm{bc}$ follow from normal hyperbolicity) for $k=1$ ($\eta_\rho=1$ and $\eta_T=0$) and $k=2$ ($\eta_\rho=0$ and $\eta_T=1$), respectively, do represent the sensitivities of $T$ and $\rho$ to constraint violations and variations in $p$.

To show this, given constraint violations $\delta_\mathrm{ode}(\cdot,\cdot)$, $\delta_\mathrm{bc}(\cdot)$, and $\delta_h$, consider the linearized equations
\begin{align}
    \delta_\mathrm{ode}&=\partial_\tau\delta_v-T\partial_xf(v,p)\delta_v-f(v,p)\delta_T-T\partial_pf(v,p)\delta_p\\
    \delta_\mathrm{bc}&=\delta_v(\cdot+\rho,0)-\delta_v(\cdot,1)+\partial_\phi v(\cdot+\rho,0)\delta_\rho\label{eq:linquasi2}\\
    \delta_h&=\int_\mathbb{S}\int_0^1\partial_v h(v(\cdot,\cdot),p)\delta_v(\phi,\tau)\,\mathrm{d}\tau\,\mathrm{d}\phi+\partial_p h(v(\cdot,\cdot),p)\delta_p.\label{eq:phasecondlin}
\end{align}
The first of these implies that
\begin{align}
    \delta_v(\cdot,\tau)&=V(\cdot,\tau)\int_0^\tau V(\cdot,\sigma)^{-1}\delta_\mathrm{ode}(\cdot,\sigma)\,\mathrm{d}\sigma\nonumber\\
    &\qquad+V(\cdot,\tau)\delta_v(\cdot,0)+\tau f(v(\cdot,\tau),p)\delta_T+P(\cdot,\tau)\delta_p,\label{eq:deltavsol}
\end{align}
where $P(\cdot,\tau)$ was given in \eqref{eq:Pquasi}. We substitute this solution with $\tau=1$ into \eqref{eq:linquasi2}, shift
the argument $\phi$ to $\phi-\rho$, replace $f(v(\phi-\rho,1),p)$ by $\frac{1}{T}\partial_\tau v(\phi,0)$ and isolate the
constraint violations on one side of the equation to obtain
\begin{align}
\label{eq:Glineq}
\delta_\mathrm{rhs}(\phi-\rho)&=\partial_\phi v(\phi,0)\delta_\rho-\frac{1}{T}\partial_\tau v(\phi,0)\delta_T\nonumber\\
&\qquad-P(\phi-\rho,1)\delta_p-\Gamma_{\rho,0}[\delta_v(\cdot,0)](\phi),
\end{align}
where
\begin{equation}
    \delta_\mathrm{rhs}(\phi):=\delta_\mathrm{bc}(\phi)+V(\phi,1)\int_0^1 V(\phi,\sigma)^{-1}\delta_\mathrm{ode}(\phi,\sigma)\,\mathrm{d}\sigma.
\end{equation}
Let $\delta_{\tau/\phi}(\phi):=q_{\tau/\phi}^{\transp}(\phi,0)\delta_v(\phi,0)$ and $\delta_{\mathrm{tr}}(\phi):=q_{\mathrm{tr}}(\phi,0)\delta_v(\phi,0)$ such that
\begin{align}
  \delta_v(\phi,0)&=\partial_\tau v(\phi,0)\delta_\tau(\phi)+\partial_\phi v(\phi,0)\delta_\phi(\phi)+ q_\mathrm{tr}(\phi,0)\delta_\mathrm{tr}(\phi)
\end{align}
and \eqref{eq:Glineq} can be written as
\begin{align}
\delta_\mathrm{rhs}(\phi-\rho)&=\partial_\phi v(\phi,0)\left(\delta_\rho+\delta_\phi(\phi)-\delta_\phi(\phi-\rho)\right)\nonumber\\
&\qquad+\partial_\tau v(\phi,0)\left(-\frac{\delta_T}{T}+\delta_\tau(\phi)-\delta_\tau(\phi-\rho)\right)\nonumber\\
&\qquad\qquad-P(\phi-\rho,1)\delta_p-\Gamma_{\rho,0}[q_\mathrm{tr}(\cdot,0)\delta_\mathrm{tr}(\cdot)](\phi).
\end{align}
It follows that
\begin{align}
    q_\mathrm{tr}(\phi,0)\left(\delta_\mathrm{rhs}(\phi-\rho)+P(\phi-\rho,1)\delta_p\right)&=-\hat{\Gamma}_{\rho,0}[q_\mathrm{tr}(\cdot,0)\delta_\mathrm{tr}(\cdot)](\phi),\\
    q_\tau^{\transp}(\phi,0)\left(\delta_\mathrm{rhs}(\phi-\rho)+P(\phi-\rho,1)\delta_p\right)&=-\frac{\delta_T}{T}+\delta_\tau(\phi)-\delta_\tau(\phi-\rho),\label{eq:deltatau}\\
     q_\phi^{\transp}(\phi,0)\left(\delta_\mathrm{rhs}(\phi-\rho)+P(\phi-\rho,1)\delta_p\right)&=\delta_\rho+\delta_\phi(\phi)-\delta_\phi(\phi-\rho).\label{eq:deltaphi}
\end{align}
Since $\hat{\Gamma}_{\rho,0}$ has a bounded inverse, the first of these equations may be solved uniquely for $\delta_\mathrm{tr}$ independently of $\delta_\tau$ and $\delta_\phi$ in terms of $\delta_\mathrm{rhs}$, $\delta_p$, $\delta_T$, and $\delta_\rho$. By integrating the second and third equations over $\mathbb{S}$, we obtain (recall that $q_{\tau/\phi}(\phi,0)=q_{\tau/\phi}(\phi-\rho,1)$)
\begin{align}
\label{eq:taumean}
    \int_\mathbb{S}q_\tau^{\transp}(\phi-\rho,1)\delta_\mathrm{rhs}(\phi-\rho)\,\mathrm{d}\phi+\int_\mathbb{S}q_\tau^{\transp}(\phi-\rho,1)P(\phi-\rho,1)\,\mathrm{d}\phi\,\delta_p=-\frac{\delta_T}{T},
\end{align}
and
\begin{align}
\label{eq:phimean}
    \int_\mathbb{S}q_\phi^{\transp}(\phi-\rho,1)\delta_\mathrm{rhs}(\phi-\rho)\,\mathrm{d}\phi+\int_\mathbb{S}q_\phi^{\transp}(\phi-\rho,1)P(\phi-\rho,1)\,\mathrm{d}\phi\,\delta_p=\delta_\rho,
\end{align}
such that, consistent with \eqref{eq: etaT_T} and \eqref{eq: etaT_rho}, $\eta_p^{\transp}\delta_p=-\delta_T$ (if $\delta_\mathrm{rhs}=0$ and $(\eta_T,\eta_\rho)=(1,0)$) or $\eta_p^{\transp}\delta_p=-\delta_\rho$ (if $\delta_\mathrm{rhs}=0$ and $(\eta_T,\eta_\rho)=(0,1)$). We note, for example, that a necessary condition for the well-posedness of the zero problem \eqref{eq:torusbvp},\,\eqref{eq:torusphase} along an
$m$-dimensional family of quasiperiodic invariant tori with fixed rotation number $\rho$ in $m+1$ parameters is the presence of a nonzero element of the row matrix $\int_\mathbb{S}q_\phi^{\transp}(\phi,0)P(\phi-\rho,1)\,\mathrm{d}\phi$. The tangent space of the family of invariant tori in parameter space is then given by the orthogonal complement to this row matrix (since $\delta_\rho=0$ as $\rho$ is fixed). For $\delta_\mathrm{rhs}=0$ and $(\eta_T,\eta_\rho)=(0,1)$, this row matrix equals $-\eta_p^{\transp}$ by \eqref{eq: etaT_rho}.

We attempt to solve for $\delta_\tau$ and $\delta_\phi$ by expanding all terms in \eqref{eq:deltatau} and \eqref{eq:deltaphi} in Fourier series per the formal ans\"{a}tze
\begin{equation}
    \delta_\tau(\phi)=\sum_{k=-\infty}^\infty e^{2\pi\mathrm{i}k\phi}\delta_{\tau,k},\,\delta_\phi(\phi)=\sum_{k=-\infty}^\infty e^{2\pi\mathrm{i}k\phi}\delta_{\phi,k},
    \label{eq:unknownFourier}
\end{equation}
\begin{equation}
    q_\tau(\phi,0)\left(\delta_\mathrm{rhs}(\phi-\rho)+P(\phi-\rho,1)\delta_p\right)=\sum_{k=-\infty}^\infty e^{2\pi\mathrm{i}k\phi}\mathrm{rhs}_{\tau,k},
\end{equation}
and
\begin{equation}
    q_\phi(\phi,0)\left(\delta_\mathrm{rhs}(\phi-\rho)+P(\phi-\rho,1)\delta_p\right)=\sum_{k=-\infty}^\infty e^{2\pi\mathrm{i}k\phi}\mathrm{rhs}_{\phi,k}.
\end{equation}
Substitution in \eqref{eq:taumean} and \eqref{eq:phimean} then results in the equations
\begin{align}
\label{forward:freq}
  \left(1-\mathrm{e}^{-2\pi\mathrm{i} k\rho}\right)\delta_{\tau/\phi,k}&=P_{\tau/\phi,k}\delta_p+\mathrm{rhs}_{\tau/\phi,k}
\end{align}
for all $k\ne 0$. Notably, $\delta_{\tau,0}$ and $\delta_{\phi,0}$ do not appear in these conditions, but can instead be solved for uniquely in terms of the remaining Fourier coefficients of $\delta_\tau$, $\delta_\phi$, and $\delta_\mathrm{tr}$ from \eqref{eq:phasecondlin} after substitution of \eqref{eq:deltavsol}, since the corresponding coefficient matrix is the nonsingular matrix \eqref{phase:ndeg:h}. 

Although the coefficients $\delta_{\tau,k}$ and $\delta_{\phi,k}$ for $k\ne 0$ may be uniquely determined from \eqref{forward:freq}, arbitrarily small divisors $1-\mathrm{e}^{-2\pi\mathrm{i} k\rho}$ occur for large $k$ and convergence of the Fourier series in \eqref{eq:unknownFourier} cannot be guaranteed. This violates the assumed existence of the bounded inverse on the right-hand side of \eqref{eq:lambdasol} and, further, the assumption of regularity relied upon in the standard implicit-function theorem\footnote{\jrem{See the reviews \cite{delaLlave2001tutorial,moser1966} for the use of generalized implicit function theorems and Diophantine   conditions on $\rho$ (i.e.,  $|\exp(2\pi\mathrm{i} k\rho)-1|\geq C_\mathrm{Diop}|k|^{-\nu}$ for all $k\neq0$ and some constants $C_\mathrm{Diop}>0$ and $\nu>0$) to establish existence of invariant tori with parallel flows and irrational rotation numbers. A general treatment for when formal   expansions (such as \eqref{forward:freq}) permit one to establish the existence of invariant manifolds with a certain degree of regularity is given in \cite{cabre2003:I,cabre2003:II,cabre2005:III}. The more recent monograph \cite{haro2016parameterization} develops numerical algorithms with rigorous error bounds for computing invariant manifolds (such as quasiperiodic tori) in the presence of unbounded inverses in \eqref{eq:lambdasol} and small divisors.}} to imply local solvability of the zero problem. That the subsystem \eqref{eq:taumean} and \eqref{eq:phimean} is well-posed, however, implies that the corresponding rows of the inverse in \eqref{eq:lambdasol} are bounded. In particular, \eqref{eq:taumean} and \eqref{eq:phimean} implies vanishing sensitivity of $T$ or $\rho$ to perturbations $\delta_h$ to the phase conditions, consistent with the observation that $\lambda_\mathrm{ps}=0$.

\subsection{An invariant curve}
\label{sec:invc}
The preceding analysis has considered two-dimensional invariant tori for flows, for which normal hyperbolicity and irrational rotation numbers guarantee the existence of two distinct asymptotic phases associated with the tangent vector fields $\partial_\tau v(\phi,\tau)$ and $\partial_\phi v(\phi,\tau)$. We illustrate the key results of this treatment but consider the simpler problem of an invariant curve of a perturbation of the normal form map for the Neimark-Sacker bifurcation \cite{kuznetsov2013elements} given by
%mrot=[cos(theta),-sin(theta);sin(theta),cos(theta)];
%M=(1+alpha)*mrot*xsym+(xsym.'*xsym)*mrot*[a,-b;b,a]*xsym+r.*xsym;
\begin{align}
  \label{invc:map}
  M(x,p)&=
  \begin{bmatrix}
    \cos\theta
    &-\sin\theta\\
    \sin\theta
    &\phantom{-}\cos\theta
  \end{bmatrix}
      \left((1+\alpha)x+
      \begin{bmatrix}
        a & -b\\
        b &\phantom{-}a
      \end{bmatrix}|x|^2x\right)+
            \begin{bmatrix}
              r_1x_1\\r_2 x_2
            \end{bmatrix}
\end{align}
with $r_1=0$, $\alpha=1/4$, $a=-1/4$, $\theta=\pi\left(\sqrt{\smash[b]{5}}-1\right)$, and $p=(r_2,b)^{\transp}$. For $p=0$, the circle $v:\mathbb{S}\mapsto\mathbb{R}^2$ with radius $\sqrt{-\alpha/a}=1$
and centered on the origin is invariant and attracting. In particular, $M(v(\phi),0)=v(\phi+\theta/2\pi)$, i.e., $\theta/2\pi$ is the corresponding rotation number.

\begin{figure}[ht]
  %\centering
  \hspace*{-3em}\includegraphics[width=1.2\textwidth]{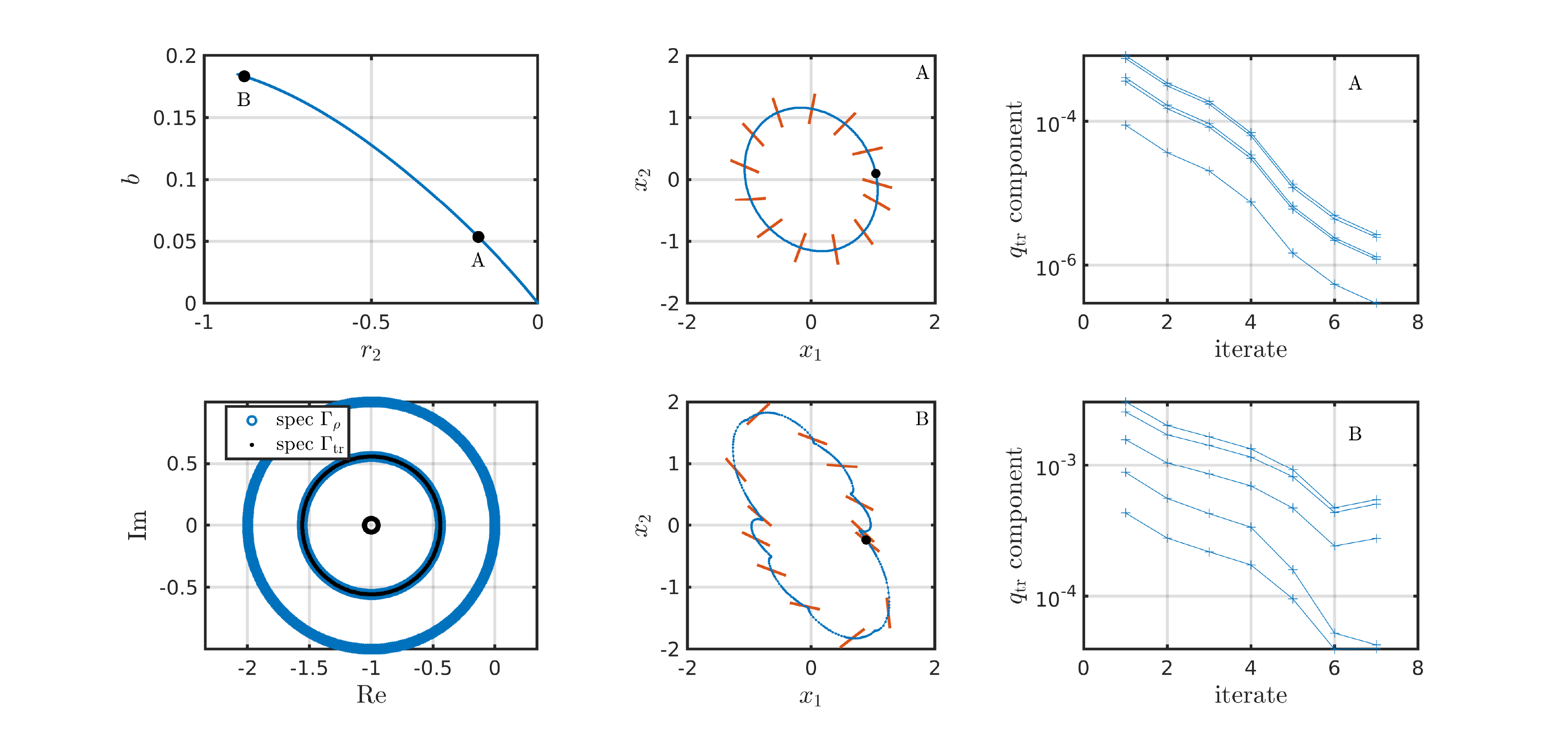}
  \caption[Invariant curve]{\jrem{Illustration of sensitivity analysis along an invariant curve: (top left) Bifurcation diagram, (bottom left) spectrum of $\Gamma_\rho$ and $\hat{\Gamma}_\rho$, (middle column) phase plots of invariant curves at points labeled A and B, and tangents to stable fibers at selected points, (right column) differences between trajectories of initial conditions near a point on the invariant curve (black dot) and forward iterates of the corresponding stable fiber projection.}}
  \label{fig:invc}
\end{figure}

As seen in the top-left panel of Fig.~\ref{fig:invc}, there exists a unique one-dimensional family of invariant curves with rotation number $\theta/2\pi$ and for simultaneous variations in both component of $p$. The panels in the middle column show two such invariant curves and illustrate a loss of smoothness (i.e., increasing higher-order derivatives of the curve function $v$) along this family.

To eliminate the degeneracy associated with arbitrary phase shifts, we append a scalar phase condition to obtain the zero problem
\begin{equation}
\label{invc:zeroproblem}
    M(v(\phi),p)=v(\phi+\rho),\,h(v(\cdot),p)=0
\end{equation}
in terms of the unknown function $v:\mathbb{S}\rightarrow\mathbb{R}^2$ and scalar $\rho$, and assume that
\begin{equation}
\label{invc:nondeg}
    \int_\mathbb{S}\partial_v h(v(\cdot),p)(\phi) v'(\phi)\,\mathrm{d}\phi\ne 0.
\end{equation}
Here, $V(\phi):=\partial_xM(v(\phi),p)$ plays the role of the matrix $V(\phi,1)$ in previous sections. The corresponding operator $\Gamma_\rho$ is then given by
\begin{equation}
  \label{invc:gro}
  \Gamma_\rho:\delta(\cdot)\mapsto V(\cdot-\rho)\delta(\cdot-\rho)-\delta(\cdot)
\end{equation}
and, in particular, $v'(\phi)$ is a right nullvector of $\Gamma_\rho$. Normal hyperbolicity implies the existence of a continuous family of projections
\begin{equation}
    q_\mathrm{tg}(\phi):= v'(\phi)q_\phi^{\transp}(\phi),\,q_\mathrm{tr}(\phi):=I_2-q_\mathrm{tg}(\phi)
\end{equation}
such that
\begin{equation}
    q_\mathrm{tg/tr}(\phi+\rho)V(\phi)=V(\phi)q_\mathrm{tg/tr}(\phi)\Leftrightarrow q_\phi^{\transp}(\phi+\rho)V(\phi)=q_\phi^{\transp}(\phi)
\end{equation}
and $\Gamma_\rho$ is a bijection with bounded inverse on the space of functions $\phi\mapsto q_\mathrm{tr}(\phi)\delta(\phi)$ for arbitrary continuous periodic $\delta(\phi)$. As in previous sections, we conclude that the linear functional $\int_\mathbb{S}q_\phi^{\transp}(\phi)(\cdot)\,\mathrm{d}\phi$ lies in the left nullspace of $\Gamma_\rho$.

By normal hyperbolicity, the map
\begin{equation}
  \label{invc:grohat}
  \hat{\Gamma}_\rho:\delta(\cdot)\mapsto
V(\cdot-\rho)q_\mathrm{tr}(\cdot-\rho)\delta(\cdot-\rho)-\delta(\cdot)
\end{equation}
has a bounded inverse on the space of continuous periodic functions $\delta(\phi)$. For a transversally stable invariant curve, eigenvalues of $\hat{\Gamma}_\rho$ lie inside a circle centered on $-1$ and of radius less than $1$. This prediction is verified by the lower-left panel of Fig.~\ref{fig:invc}, which also shows that $\Gamma_\rho$ has eigenvalues accumulating to $0$. Analogously with the results in previous sections,
\begin{equation}
    q_\phi^{\transp}(\phi)=\lim_{k\rightarrow\infty}\frac{v^{\prime\,{\transp}}(\phi+k\rho)}{\| v'(\phi+k\rho)\|^2}V(\phi)^k
\end{equation}
with exponential convergence. Alternatively, for irrational $\rho$, $q_\phi(\phi)$ may be obtained from the unique continuous solution $\lambda_\mathrm{map}(\phi-\rho)$ of the adjoint conditions
\begin{align}
    0&=\lambda_\mathrm{map}^{\transp}(\phi)V(\phi)-\lambda_\mathrm{map}^{\transp}(\phi-\rho)+\lambda_\mathrm{ps}\partial_v h(v(\cdot),p)(\phi)\label{invc:delv}\\
    0&=-\int_\mathbb{S}\lambda_\mathrm{map}^{\transp}(\phi)v'(\phi+\rho)\,\mathrm{d}\phi+\eta_\rho\label{invc:delrho}\\
    0&=\int_\mathbb{S}\lambda_\mathrm{map}^{\transp}(\phi)\partial_p M(v(\phi),p)\,\mathrm{d}\phi+\lambda_\mathrm{ps}\partial_p h(v(\cdot),p)+\eta_p^{\transp}\label{invc:delp}
\end{align}
with $\eta_\rho=1$. Indeed, multiplication of \eqref{invc:delv} by $v'(\phi)$ and integration over $\mathbb{S}$ shows that $\lambda_\mathrm{ps}=0$ and, consequently, that $\lambda_\mathrm{map}^{\transp}(\phi)V(\phi)=\lambda_\mathrm{map}^{\transp}(\phi-\rho)$ and $\int_\mathbb{S}\lambda_\mathrm{map}^{\transp}(\phi-\rho)(\cdot)\,\mathrm{d}\phi$ lies in the left nullspace of $\Gamma_\rho$.

A forward linear sensitivity analysis requires solving the linear system
\begin{align}
  \label{invc:forward:Gamma}
  \delta_\mathrm{rhs}(\phi)&=V(\phi-\rho)\delta_v(\phi-\rho)-\delta_v(\phi)+v'(\phi)\delta_\rho\\
  \label{invc:forward:ps}
    \delta_h&=\int_\mathbb{S}\partial_xh(v(\cdot),p)(\phi)\delta_v(\phi)\,\mathrm{d}\phi+\partial_p h(v(\cdot),p)\delta_p
\end{align}
with arbitrary $\delta_\mathrm{rhs}(\phi)\in\mathbb{R}^2$ and
$\delta_\mathrm{ps}\in\mathbb{R}$, for $\delta_v(\phi)\in\mathbb{R}^2$, $\delta_\rho\in\mathbb{R}$, and $\delta_p\in\mathbb{R}^2$. Since $\Gamma_\rho$ has a simple
eigenvalue $0$ with eigenvector $v'$, the non-degeneracy condition \eqref{invc:nondeg} ensures that the \eqref{invc:forward:Gamma},\,\eqref{invc:forward:ps}
is solvable. However, as the spectrum of $\Gamma_\rho$ in
Fig.~\ref{fig:invc} indicates, the solution $\delta_v$ is not bounded by $\delta_\mathrm{rhs}$ in the same norm due to the small-divisor problem discussed in Section~\ref{sec:regularity}.

Finally, we illustrate the predicted asymptotic convergence of forward iterates of $v(\phi_0,p)+\delta_0$ toward forward iterates of $v(\phi_0+q_\phi(\phi_0)\delta_0)$ per the theory of asymptotic phase. The panels in the right column of Figure~\ref{fig:invc} show the norm
\begin{multline}
  \label{invc:asymptotic:phase}
  M^k(v(\phi_0,p)+\delta_0,p)-M^k(v(\phi+q_\phi(\phi_0)\delta_0,p),p)\approx\\
  M^k(v(\phi_0,p)+\delta_0,p)-M^k(v(\phi_0,p)+\partial_\phi v(\phi_0,p)q_\phi(\phi_0)\delta_0,p)
\end{multline}
for $20$ initial conditions at distances $|\delta_0|\approx 10^{-4}$
from the points indicated by black dots in the phase portraits in the middle column in Fig.~\ref{fig:invc}.  We observe exponential convergence up to an error of
$|\delta_0|^2\approx 10^{-7}$, consistent with the theoretical prediction \eqref{phase:linear}.

\section{Implementation in \textsc{COCO}}
\label{sec:construction}
We use this section to further highlight the algorithmic nature of the adjoint approach and its implementation in \textsc{coco}. A general discussion that represents a reference for users and toolbox developers is followed by the explicit encoding in \textsc{coco} of the combined sensitivity of $\rho$ with respect to $r_2$ and $b$ along the family of invariant curves in Section~\ref{sec:invc}.

\subsection{Principles of construction}
As alluded to in Section~\ref{sec:preliminaries}, the construction philosophy of the \textsc{coco} software package naturally lends itself to the Lagrangian approach to computing the sensitivities of different monitor functions with respect to violations of constraints and variations of continuation parameters. Specifically, by the additive nature of the problem Lagrangian and the associated adjoint conditions, it is possible to arrive at a complete set of defining equations in terms of the adjoint variables in multiple stages, mirroring the staged addition of zero functions and monitor functions to the extended continuation problem.

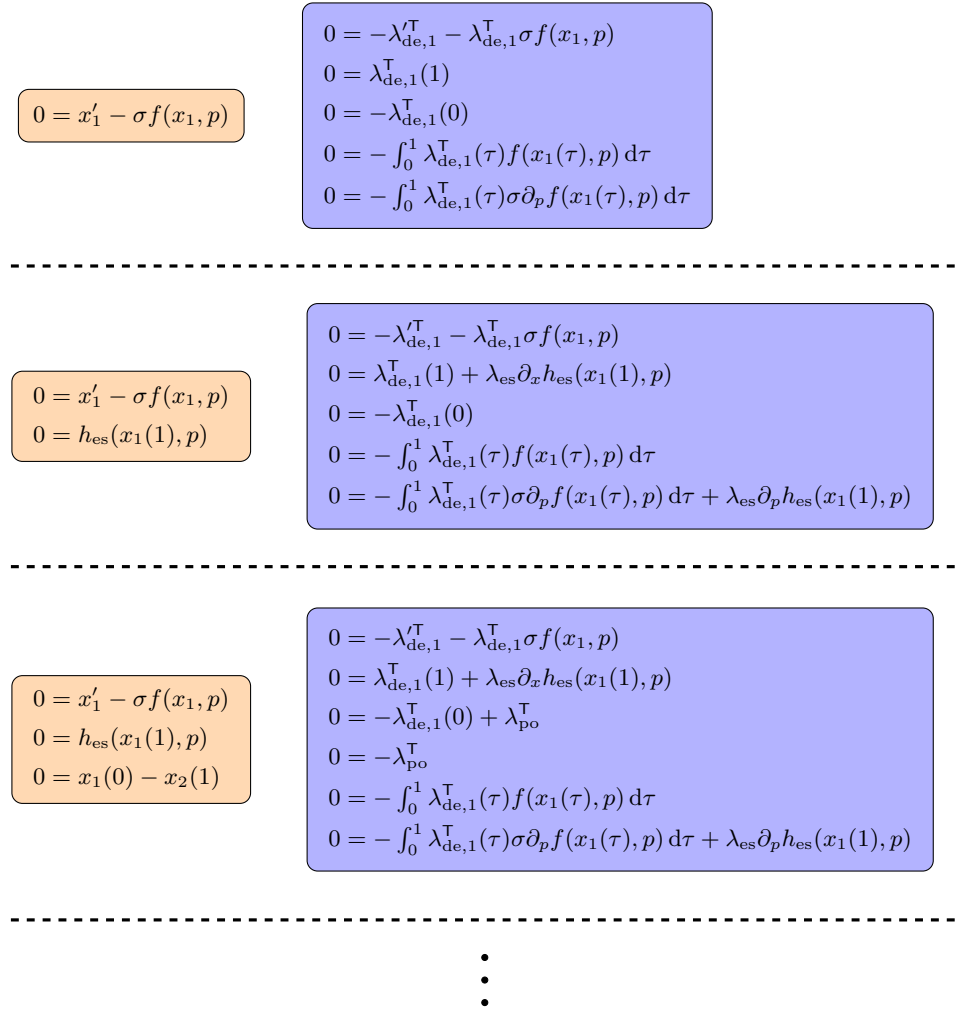
\begin{figure}[h]
    \centering
    {\small
    \begin{tikzpicture}
\node (c1) [rectangle, rounded corners, minimum width=3cm, minimum height=.7cm, text centered, draw=black, fill=orange!30] {$0=x'_1-\sigma f(x_1,p)$};
\node (a1) [rectangle, rounded corners, minimum height=3cm, text centered, draw=black, fill=blue!30,right of=c1,xshift=4cm] {$\begin{array}{l}
    0=-\lambda_{\mathrm{de},1}^{\prime{\transp}}-\lambda_{\mathrm{de},1}^{\transp}\sigma f(x_1,p)\\[.4em]
    0=\lambda_{\mathrm{de},1}^{\transp}(1)\\[.4em]
    0=-\lambda_{\mathrm{de},1}^{\transp}(0)\\[.4em]
    0=-\int_0^1 \lambda_{\mathrm{de},1}^{\transp}(\tau)f(x_1(\tau),p)\,\mathrm{d\tau}\\[.4em]
    0=-\int_0^1 \lambda_{\mathrm{de},1}^{\transp}(\tau)\sigma\partial_p f(x_1(\tau),p)\,\mathrm{d}\tau
\end{array}$};
\node (c2) [rectangle, rounded corners, minimum width=3cm, minimum height=1.2cm, text centered, draw=black, fill=orange!30,below of=c1,yshift=-3cm] {$\begin{array}{l}
0=x'_1-\sigma f(x_1,p)\\[.4em]0=h_\mathrm{es}(x_1(1),p)\end{array}$};
\node (a2) [rectangle, rounded corners, minimum height=3cm, text centered, draw=black, fill=blue!30,below of=a1,yshift=-3cm,xshift=1.5cm] {$\begin{array}{l}
    0=-\lambda_{\mathrm{de},1}^{\prime{\transp}}-\lambda_{\mathrm{de},1}^{\transp}\sigma f(x_1,p)\\[.4em]
    0=\lambda_{\mathrm{de},1}^{\transp}(1)+\lambda_\mathrm{es}\partial_x h_\mathrm{es}(x_1(1),p)\\[.4em]
    0=-\lambda_{\mathrm{de},1}^{\transp}(0)\\[.4em]
    0=-\int_0^1 \lambda_{\mathrm{de},1}^{\transp}(\tau)f(x_1(\tau),p)\,\mathrm{d\tau}\\[.4em]
    0=-\int_0^1 \lambda_{\mathrm{de},1}^{\transp}(\tau)\sigma\partial_p f(x_1(\tau),p)\,\mathrm{d}\tau+\lambda_\mathrm{es}\partial_p h_\mathrm{es}(x_1(1),p)
\end{array}$};
\node (c3) [rectangle, rounded corners, minimum width=3cm, minimum height=1.7cm, text centered, draw=black, fill=orange!30,below of=c2,yshift=-3.3cm] {$\begin{array}{l}
0=x'_1-\sigma f(x_1,p)\\[.4em]0=h_\mathrm{es}(x_1(1),p)\\[.4em]0=x_1(0)-x_2(1)\end{array}$};
\node (a3) [rectangle, rounded corners, minimum height=3.5cm, text centered, draw=black, fill=blue!30,below of=a2,yshift=-3.3cm] {$\begin{array}{l}
    0=-\lambda_{\mathrm{de},1}^{\prime{\transp}}-\lambda_{\mathrm{de},1}^{\transp}\sigma f(x_1,p)\\[.4em]
    0=\lambda_{\mathrm{de},1}^{\transp}(1)+\lambda_\mathrm{es}\partial_x h_\mathrm{es}(x_1(1),p)\\[.4em]
    0=-\lambda_{\mathrm{de},1}^{\transp}(0)+\lambda_\mathrm{po}^{\transp}\\[.4em]
    0=-\lambda_\mathrm{po}^{\transp}\\[.4em]
    0=-\int_0^1 \lambda_{\mathrm{de},1}^{\transp}(\tau)f(x_1(\tau),p)\,\mathrm{d\tau}\\[.4em]
    0=-\int_0^1 \lambda_{\mathrm{de},1}^{\transp}(\tau)\sigma\partial_p f(x_1(\tau),p)\,\mathrm{d}\tau+\lambda_\mathrm{es}\partial_p h_\mathrm{es}(x_1(1),p)
\end{array}$};
\draw [dashed, very thick] (-1.6,-2) -- (11,-2);
\draw [dashed, very thick] (-1.6,-6) -- (11,-6);
\draw [dashed, very thick] (-1.6,-10.7) -- (11,-10.7);
\draw[fill] (4.7,-11.2) circle (0.4mm);
\draw[fill] (4.7,-11.5) circle (0.4mm);
\draw[fill] (4.7,-11.8) circle (0.4mm);
\end{tikzpicture}
}
    \caption{At each stage of construction of the two-segment periodic orbit problem in \eqref{eq:twoseg1}-\eqref{eq:twoseg3}, one obtains an extended continuation problem (left column) in terms of a subset of the complete vector of continuation variables and an associated set of adjoint conditions (right column) in terms of a subset of the complete vector of adjoint variables.}
    \label{fig:staged construction}
\end{figure}

As an illustration of the staged construction paradigm, consider again the analysis of the two-segment periodic orbit of period $T$ in Section~\ref{sec:Hybrid dynamics}. Here, the contributions to the adjoint conditions associated with individual constraints are clearly identified by the corresponding adjoint variables. As constraints are appended individually or in groups, the corresponding contributions to the adjoint conditions may be constructed in terms of linear operations on the corresponding adjoint variables. At any stage of construction, one obtains an extended continuation problem with an associated set of adjoint conditions. We illustrate this principle in Fig.~\ref{fig:staged construction}.

A key concept of the \textsc{coco} construction philosophy is the idea that a constructor that appends a zero or monitor function must be able to operate on an existing extended continuation problem in order to also expand its domains of continuation variables and/or continuation parameters. We say that a constructor is \textit{embeddable} if this is the case. One way this may be visualized is in terms of an operator on the space of logical matrices of arbitrary dimension, whose columns represent a Cartesian decomposition of the domain of continuation variables and whose rows represent subsets of zero functions or monitor functions introduced at individual stages of construction. With each call to an embeddable constructor, the matrix grows by any number of additional columns with all zero entries followed by the addition of rows with a combination of zeros and ones in the previously defined columns and only ones in the newly defined columns. An example is shown in Fig.~\ref{fig:booleanconstraint}, where the interpretation is that the zero or monitor functions added by the constructor in this call depend on the variables identified by ones in the corresponding row, and not on those identified by zeros.
\begin{figure}[h]
    \centering
    \begin{tikzpicture}
    \node (o) {$\begin{bmatrix}
      \cdot & \cdot & \cdot & \cdot \\ \cdot & \cdot & \cdot & \cdot \\ \cdot & \cdot & \cdot & \cdot
    \end{bmatrix}$};
    \node (n) [right of=o, xshift=3cm] {$\begin{bmatrix}
      \cdot & \cdot & \cdot & \cdot & 0 & 0 & 0\\ \cdot & \cdot & \cdot & \cdot & 0 & 0 & 0\\ \cdot & \cdot & \cdot & \cdot & 0 & 0 & 0\\
      1 & 0 & 0 & 0 & 1 & 1 & 1\\
      0 & 1 & 0 & 1 & 1 & 1 & 1
    \end{bmatrix}$};
    \draw [arrow] (o) -- (n);
    \end{tikzpicture}
    \caption{\jrem{A logical matrix representation of the variable dependencies of zero and monitor functions introduced at different stages of problem construction.} At the stage of construction illustrated in the figure, two zero/monitor functions are added that depend on three newly introduced continuation variables. Each function also depends on previously introduced continuation variables identified by ones in the corresponding row.}
    \label{fig:booleanconstraint}
\end{figure}

The new columns added by a call to an embeddable constructor represent continuation variables that are meaningful to the constructor, but not to \textsc{coco}, which is unable to distinguish among subsets of these variables other than through the cardinal number of the corresponding column. The latter depends on the order of construction and on internal details of the constructor, neither of which can be assumed to be known by subsequent constructors, except in the most simple cases. In \textsc{coco}, this challenge is overcome through a subindexing mechanism, whereby constructor-dependent details are stored in a \textit{function data structure} associated with a particular matrix row and used to extract meaningful subsets of continuation variables associated with this row. This mechanism is illustrated in Fig.~\ref{fig:subindexing}.
\begin{figure}[h]
    \centering
    \begin{tikzpicture}
    \node (nr) {$u=\begin{bmatrix}
      \vdots\\[.2em]\bullet\\[.2em]\bullet\\[.2em]\bullet\\[.2em]\bullet\\[.2em]\bullet\\[.2em]\bullet\\\vdots
    \end{bmatrix}$};
    \draw[xshift=.5cm,decoration={calligraphic brace,amplitude=4pt}, decorate, line width=1.5pt] (.2,1.35) -- (.2,0.6) node [midway,xshift=1.5cm] {\mcode{fcndata.x1_idx}};
    \draw[xshift=.5cm,decoration={calligraphic brace,amplitude=4pt}, decorate, line width=1.5pt] (.2,0.4) -- (.2,-1.4) node [midway,xshift=1.5cm] {\mcode{fcndata.x2_idx}};
    \end{tikzpicture}
    \caption{Fields of a function data structure may be used to store context-independent integer indices for subsets of continuation variables associated with a newly constructed function for later reference by subsequent constructor calls. In the figure, the fields \mcode{fcndata.x1_idx} and \mcode{fcndata.x2_idx} reference constructor-dependent details that are unknown to \textsc{coco} and independent of the order of construction.}
    \label{fig:subindexing}
\end{figure}

As shown in Section~\ref{sec:preliminaries}, the collection of adjoint variables naturally decomposes through a Cartesian product of dual spaces to each of the ranges of the zero functions or monitor functions introduced at individual stages of construction. As additional zero or monitor functions are introduced, terms are added to existing adjoint conditions or used to initialize new conditions. The analogy with the dependence on previously defined continuation variables and the introduction of new continuation variables in the construction of zero and monitor functions suggests a similar use of a logical matrix representation to track the target locations of new contributions to the adjoint conditions. In each newly added column, rows with ones imply contributions to a particular (existing or newly constructed) adjoint condition. The analogy is further extended by the realization that it is again necessary to use function data for constructor-specific details that identify the associated adjoint conditions independently of the order of construction. This desired functionality is illustrated in Fig.~\ref{fig:booleanadjoint}. 
\begin{figure}[h]
    \centering
    \begin{tikzpicture}
    \node (o) {$\begin{bmatrix}
      \cdot & \cdot & \cdot & \cdot \\ \cdot & \cdot & \cdot & \cdot \\ \cdot & \cdot & \cdot & \cdot
    \end{bmatrix}$};
    \node (n) [right of=o, xshift=3cm] {$\begin{bmatrix}
      \cdot & \cdot & \cdot & \cdot & 0 & 0 & 0\\ \cdot & \cdot & \cdot & \cdot & 0 & 0 & 0\\ \cdot & \cdot & \cdot & \cdot & 0 & 0 & 0\\
      0 & 1 & 0 & 0 & 1 & 1 & 1\\
      1 & 0 & 0 & 1 & 1 & 1 & 1
    \end{bmatrix}$};
    \draw [arrow] (o) -- (n);
    \end{tikzpicture}
    \caption{\jrem{A logical matrix representation of the dependencies of the adjoint conditions on adjoint variables introduced at different stages of problem construction.} At the stage of construction illustrated in the figure, three adjoint conditions are added that include terms linear in two newly introduced adjoint variables. Each variable also appears in contributions to previously introduced adjoint conditions identified by ones in the corresponding row.}
    \label{fig:booleanadjoint}
\end{figure}

\textsc{coco} toolboxes collect predefined embeddable zero and monitor function constructors that exhibit no dependence on previously defined sets of continuation variables, as well as adjoint constructors that, consequently, do not contribute terms to existing adjoint conditions. Since the examples considered in this paper cover classes of problems of a general, largely problem-independent nature, it is sensible to design toolboxes to enable their immediate analysis. For the forward dynamics problems in Section~\ref{sec:Forward dynamics}, the \textsc{coco}-compatible \mcode{coll} toolbox achieves this objective and may be applied out of the box. For the periodic orbit problem in Section~\ref{sec:Periodic orbits}, the \mcode{bvp} constructors in the \mcode{coll} toolbox achieve the desired result. Alternatively, if the discrete phase condition given by the vanishing of $h(x(1),p)$ is replaced by an integral phase condition, the \textsc{coco}-compatible \mcode{po} toolbox provides the sought support. The first example in Section~\ref{sec:Hybrid dynamics} is best handled with two consecutive calls to \mcode{coll} constructors and the addition of the intersegment boundary conditions using the \textsc{coco} core constructors \mcode{coco_add_func} and \mcode{coco_add_adjt} that were already deployed in Section~\ref{sec:illustration}. For the periodic orbit problem in Section~\ref{sec:Hybrid dynamics} (and, indeed, for arbitrary multi-segment periodic orbit problems), the \mcode{hspo} constructors in the \mcode{po} toolbox implement the desired formalism.

For the case of a quasiperiodic invariant torus, the \mcode{bvp} constructors may be applied to a discretization of the PDE constraints in terms of the coefficients of truncated Fourier series for every $\phi$-dependent function. Since these constructors assume that all unknowns are either discretized state variables, interval durations, or problem parameters, the rotation number $\rho$ must be thought of as a problem parameter, even though the vector field is independent of $\rho$.

\subsection{Demo for invariant curve}
\label{sec:invc:coco}
We allow explicit \textsc{coco} code to illustrate the paradigm of construction and analysis \jrem{for the invariant curves studied in Section~\ref{sec:invc}. The dynamics for this example is given by the  map $M(\cdot,r_2,b):\mathbb{R}^2\to\mathbb{R}^2$, defined in \eqref{invc:map}.} 

To this end, we \jrem{let the rotation number $\theta/2\pi$ equal the golden mean  $(\sqrt{5}-1)/2$ and approximate this} by the truncated continued fraction expansion $p/q$ with $p=233$ and $q=377$ (with error $\approx 3\times10^{-6}$). The map $M$ and its derivatives $\partial_xM$, $\partial_{r_2}M$, and
$\partial_bM$ are then encoded in \textsc{Matlab} using anonymous functions as follows
\begin{lstlisting}[language=coco-highlight]
>> [numer, denom] = deal(233, 377);
>> [theta, alpha, a, r1] = deal(2*pi*numer/denom, 1/4, -1/4, 0);
>> mrot = [cos(theta), -sin(theta); sin(theta), cos(theta)];
>> M    = @(x,r2,b) ((1+alpha)*mrot + ...
     (x'*x)*mrot*[a, -b; b, a] + [r1, 0; 0, r2])*x;
>> dMx  = @(x,r2,b) (1+alpha)*mrot + ...
     (x.'*x)*mrot*[a, -b; b, a] + [r1, 0; 0, r2] + ...
     2*mrot*[a, -b; b, a]*x*x';
>> dMr2 = @(x,r2,b) [0; x(2)];
>> dMb  = @(x,r2,b) (x'*x)*mrot*[0, -1; 1, 0]*x;
\end{lstlisting}
We proceed to discretize the first half of the zero problem \eqref{invc:zeroproblem} in terms of the sequence $\{v_i\}_{i=1}^{q}$ for angles $\phi_i=2\pi(i-1)/q$ such that $v(\phi_i)\approx v_i$  for $i=1,\ldots,q$, and 
\begin{equation}
\label{invc:approxzero}
    M(v_i,p)-v_{\operatorname{mod}(i+p-1,q)+1}-q\left(v_{\operatorname{mod}(i+p,q)+1}-v_{\operatorname{mod}(i+p-1,q)+1}\right)\delta_\rho \approx 0,
\end{equation}
where $\delta_{\rho}=\rho-\theta/2\pi$ is assumed to be small. (Here, the coefficient of $\delta_\rho$ is an approximation for $\partial_\phi v(\phi_{i+p})$, identifying $\phi_{i+q+1}=\phi_i$ for all $i\in\mathbb{Z}$). It follows that when $\delta_\rho=0$, the sequence $\{v_i\}_{i=1}^{q}$ is an approximate periodic orbit of period $q$ on the invariant curve that makes $p$ excursions around the invariant curve before repeating. The commands
\begin{lstlisting}[language=coco-highlight]
>> [r20, b0, drho0] = deal(0, 0, 0);
>> phi = 2*pi*linspace(0, 1-1/denom, denom);
>> v0 = [cos(phi); sin(phi)]*sqrt(-alpha/a);
\end{lstlisting}
generate the corresponding initial solution guess for $(r_2,b)=(0,0)$ with
\begin{equation}
    v_{0,i}=\sqrt{\frac{-\alpha}{a}}\begin{pmatrix}\cos \phi_i \\ \sin \phi_i \end{pmatrix}.
\end{equation}
We use this sequence to construct the discretized phase condition
\begin{align}
  \label{eq:invc:phascond}
  \sum_{i=1}^q[v_{0,\operatorname{mod}(i,q)+1}-v_{0,i}]^{\transp}[v_i-v_{0,i}]=0
\end{align}
as an approximation of the integral condition
\begin{equation}
    \int_\mathbb{S}\partial_\phi v_0^{\transp}(\phi)\left(v(\phi)-v_0(\phi)\right)\,\mathrm{d}\phi=0
\end{equation}
in terms of the function
\begin{equation}
    v_0(\phi)=\sqrt{\frac{-\alpha}{a}}\begin{pmatrix}\cos \phi \\ \sin \phi \end{pmatrix}.
\end{equation}

We construct a zero problem consisting of the vanishing left-hand side of \eqref{invc:approxzero} for $i=1,\ldots,q$ and \eqref{eq:invc:phascond} using the staged construction paradigm of \textsc{coco} as follows. We begin by constructing the function
\begin{equation}
    M_\mathrm{res}(u):=M(u_{1,2},u_5,u_6)-u_{3,4}
\end{equation}
and its Jacobian per the following commands:
\begin{lstlisting}[language=coco-highlight]
>> [mx, my, mr2, mb, mdrho] = deal(1:2, 3:4, 5, 6, 7);
>> Mres  = @(u) M(u(mx), u(mr2), u(mb)) - u(my);
>> dMres = @(u) [dMx(u(mx), u(mr2), u(mb)), -eye(2), ...
     dMr2(u(mx),u(mr2),u(mb)),dMb(u(mx),u(mr2),u(mb)),zeros(2,1)];
\end{lstlisting}
Then, $M_\mathrm{res}(x,y,r_2,b,\delta_\rho)=0$ implies that $y$ is the image of $x$ under the map $M(\cdot,r_2,b)$. Without loss of generality, we carry $\delta_\rho$ along even in the absence of explicit dependence on this variable. The following sequence of commands create $q$ instances of the zero problem $M_\mathrm{res}=0$, initialized with
\begin{align}
u_0=(v_{0,i},v_{0,\operatorname{mod}(i+p-1,q)+1},r_{2,0},b_0,0),
\end{align}
and labeled with the function identifiers \mcode{'M1'}, \mcode{'M2'}, and so on. For each instance, we also construct the corresponding contributions to the adjoint conditions.
\begin{lstlisting}[language=coco-highlight]
>> fcn  = @(f) @(p,d,u) deal(d, f(u)); % convert plain fcn to COCO fcn
>> fid  = @(s,i)[s,num2str(i)];        % enumerated function identifiers
>> prob = coco_prob;
>> for i=1:denom
     irot = mod(i+numer-1, denom)+1;
     prob = coco_add_func(prob, fid('M',i), fcn(Mres), fcn(dMres), ...
       [], 'zero', 'u0', [v0(:,i); v0(:,irot); r20; b0; drho0]);
     prob = coco_add_adjt(prob, fid('M',i));
     muidx(:,i) = coco_get_func_data(prob, fid('M',i), 'uidx');
     maidx(:,i) = coco_get_adjt_data(prob, fid('M',i), 'axidx');
   end
\end{lstlisting}
We store the indices that track the internal ordering of the corresponding continuation variables in the $7\times q$ index array
\mcode{muidx}. The column indices for the associated adjoint conditions are correspondingly kept in the
$7\times q$ index array \mcode{maidx}. 

The $q$ zero problems introduced thus far are uncoupled. We achieve the coupling implied by the vanishing of the left-hand side of \eqref{invc:approxzero} using the following sequence of commands.
\begin{lstlisting}[language=coco-highlight]
>> [bx, bxnext, by, bdrho] = deal(1:2, 3:4, 5:6, 7);
>> fbc = @(u) u(bx) + u(bdrho)*denom*(u(bxnext) - u(bx)) - u(by);
>> dbc = @(u) [eye(2)*(1-denom*u(bdrho)), eye(2)*denom*u(bdrho), ...
     -eye(2), denom*(u(bxnext)-u(bx))];
>> for i=1:denom
     irot  = mod(i+numer-1,denom)+1;
     inext = mod(i+numer,denom)+1;
     prob  = coco_add_func(prob, fid('bc',i), fcn(fbc), fcn(dbc), ...
       [],'zero', 'uidx', [muidx(mx,irot); muidx(mx,inext); ...
       muidx(my,i); muidx(mdrho,i)]);
     prob  = coco_add_adjt(prob, fid('bc',i), ...
       'aidx', [maidx(mx,irot); maidx(mx,inext); ...
       maidx(my,i); maidx(mdrho,i)]);
   end
\end{lstlisting}
We use the index array \mcode{muidx} to refer to the previously initialized continuation variables and the index array \mcode{maidx} to refer to previously initialized adjoint conditions. As each of the $q$ zero problems \mcode{'M1'}, \mcode{'M2'}, \ldots depend on their own instances of the parameters $r_2$, $b$, and $\rho$, we use the \mcode{coco_add_glue} constructor to glue these together across all instances, as shown below.
\begin{lstlisting}[language=coco-highlight]
>> for i=2:denom
     prob = coco_add_glue(prob, fid('pglue',i),...
       muidx([mr2, mb, mdrho], 1), muidx([mr2, mb, mdrho], i));
     prob = coco_add_adjt(prob, fid('pglue',i),...
       'aidx', [maidx([mr2, mb, mdrho], 1); maidx([mr2, mb, mdrho], i)]);
   end
\end{lstlisting}
Finally, we add the phase condition, per the construction
\begin{lstlisting}[language=coco-highlight]
>> dx0       = denom*(v0(:,[2:end,1])-v0);
>> phascond  = @(u)dx0(:)'*(u(:)-v0(:));
>> dphascond = @(u)dx0(:)';
>> prob = coco_add_func(prob, 'phasecond', fcn(phascond), ...
     fcn(dphascond), [], 'zero', 'uidx', muidx(mx,:));
>> prob = coco_add_adjt(prob, 'phasecond', 'aidx', maidx(mx,:));
\end{lstlisting}

Following the analysis in previous sections, we append monitor functions that evaluate to the instances of $r_2$, $b$, and $\delta_\rho$ associated with \mcode{'M1'} and label the corresponding continuation parameters by \mcode{'r2'}, \mcode{'b'}, and \mcode{'drho'}. By default these are initially inactive. 
\begin{lstlisting}[language=coco-highlight]
>> prob = coco_add_pars(prob, 'pars', muidx([mr2, mb, mdrho], 1), ...
     {'r2','b','drho'});
\end{lstlisting}
The call to \mcode{coco_add_adjt} below then appends complementary monitor functions whose values equal $\eta_{r_2}$, $\eta_b$, and $\eta_\rho$, respectively. We associate these with additional complementary continuation parameters, designated by the string labels \mcode{'e.r2'}, \mcode{'e.b'}, and \mcode{'e.drho'}. These are also initially inactive.   
\begin{lstlisting}[language=coco-highlight]
>> prob = coco_add_adjt(prob, 'pars', {'e.r2','e.b','e.drho'}, ...
     'aidx', maidx([mr2, mb, mdrho], 1), 'l0', [0; 0; 1]);
\end{lstlisting}
Here, we set $\eta_\rho=1$, as we intend to determine the sensitivity of $\rho$ with respect to all other variables. 

We obtain results compatible with Fig.~\ref{fig:invc} by allowing \mcode{'r2'}, \mcode{'b'}, \mcode{'e.r2'}, and \mcode{'e.b'} to vary, while holding \mcode{'rho'} and \mcode{'e.rho'} fixed. \jrem{The commands below define a user-defined solution point at $r_2=-0.16$ (the point labeled A in Fig.~\ref{fig:invc}) and perform continuation of the corresponding augmented continuation problem on the computational domain defined by $r_2\in[-0.9,0]$.}
\begin{lstlisting}[language=coco-highlight]
>> prob = coco_add_event(prob, 'A', 'r2', -0.16);
>> coco(prob, 'run', [], 1, {'r2', 'b', 'e.r2', 'e.b'}, [-0.9,0]);
\end{lstlisting}
\begin{lstlisting}[language=coco-small]
     STEP   DAMPING               NORMS              COMPUTATION TIMES
   IT SIT     GAMMA     ||d||     ||f||     ||U||   F(x)  DF(x)  SOLVE
   0                          1.00e+00  2.75e+01    0.1    0.0    0.0
   1   1  1.00e+00  2.44e-01  5.15e-15  2.75e+01    0.2    0.6    0.0
   2   1  1.00e+00  1.11e-14  3.04e-15  2.75e+01    0.3    1.1    0.1

 ...  LABEL  TYPE            r2            b         e.r2          e.b
 ...      1  EP      0.0000e+00   2.5649e-18  -5.3757e-02  -1.5916e-01
 ...      2         -1.3955e-01   4.2507e-02  -5.3742e-02  -1.9276e-01
 ...      3  A      -1.6000e-01   4.8149e-02  -5.4033e-02  -1.9790e-01
 ...      4         -3.2179e-01   8.9201e-02  -5.7266e-02  -2.4255e-01
 ...      5         -5.1513e-01   1.3088e-01  -6.0551e-02  -3.1210e-01
 ...      6         -7.0947e-01   1.6371e-01  -6.0449e-02  -4.2622e-01
 ...      7         -8.8006e-01   1.8328e-01  -6.3439e-02  -7.2503e-01
 ...      8  EP     -9.0000e-01   1.8498e-01  -7.5930e-02  -9.2511e-01
\end{lstlisting}
\jrem{The analysis locates the point A at $(r_2,b)\approx (-0.16,0.04815)$, and finds that $\eta_{r_2}\approx -0.054033$ and $\eta_b\approx -0.19790$ are the corresponding sensitivites of $\delta_\rho$ to variations in $\delta_{r_2}$ and $\delta_b$, respectively, i.e.,
\begin{equation}
    \delta_\rho=-\eta_p^{\transp}\delta_p=-(\eta_{r_2},\eta_b)\begin{pmatrix}
    \delta_{r_2}\\\delta_b
    \end{pmatrix}\approx 0.0540\delta_{r_2}+0.198\delta_b\mbox{,}
\end{equation}
assuming no violations of the governing constraints. Since $\delta_\rho=0$ along the family of invariant tori, it follows that $\delta_b/\delta_{r_2}\approx -0.273$ at this location. This ratio equals the slope at point A of the curve of invariant tori in the $(r_2,b)$-plane in Fig.~\ref{fig:invc}.}

\jrem{The tangents to stable fibers shown in Fig.~\ref{fig:invc} are given by the adjoint variables $\lambda_{\mathrm{map},i}$. These} can now be extracted from the data stored with the solution label \mcode{lab} using the commands
\begin{lstlisting}[language=coco-highlight]
>> chart = coco_read_adjoint(fid('M',i), 'run', lab, 'chart');
>> getfield(chart, 'x')
\end{lstlisting}
\jrem{Finally, for the illustration of the spectrum of the linear map $\Gamma_\rho$ in Fig.~\ref{fig:invc}, we evaluate the map $\partial_xM$ on the solutions generated by \textsc{coco}. The map $V(\cdot-\rho)$ is then approximated by  
\begin{align}\label{invc:Gamma}
    \mathrm{diag}\left(\partial_xM(x_{\sigma_{-p,q}(1)},r_2,b),\ldots,\partial_xM(x_{\sigma_{-p,q}(q)},r_2,b)\right)\cdot (\sigma_{-p,q}\otimes I_2),
\end{align}
where $\sigma_{p,q}$ denotes the permutation (and permutation matrix) of the first $q\in\mathbb{Z}$ integers corresponding to a rotation by $p\in\mathbb{Z}$, and $\mathrm{diag}(A_1,\ldots,A_q)$ denotes the blockdiagonal matrix with $A_j$ on the diagonal. This is computed using the following sequence of commands. First, we read the solution data from the file for the label of point A, extracting the $2\times q$ array for the solution $x$, and the parameters $r_2$ and $b$.}
\begin{lstlisting}[language=coco-highlight]
>> chart_p = coco_read_solution('run', 3, 'chart');
>> x = chart_p.x(muidx(mx,:));
>> [r2, b] = deal(chart_p.x(5), chart_p.x(6));
\end{lstlisting}
\jrem{Then we construct the rotation, corresponding to $x\mapsto x(\cdot-\rho)$.}
\begin{lstlisting}[language=coco-highlight]
>> perm = diag(ones(numer,1),denom-numer) +...
     diag(ones(denom-numer,1),-numer);
>> sigma = kron(perm, eye(2));
\end{lstlisting}
\jrem{Finally, we apply this rotation to $x$, implementing expression \eqref{invc:Gamma}, by applying $\partial_xM$ and multiplying with the rotation matrix again. }
\begin{lstlisting}[language=coco-highlight]
>> xrot = num2cell(reshape(sigma*x(:), 2, denom), 1);
>> dMxval = cellfun(@(x)dMx(x,r2,b), xrot, 'uniformoutput', false);
>> blkdiag(dMxval{:})*sigma;
\end{lstlisting}

\section{Concluding discussion}
\label{sec:conclusions}

A goal of the discussion in the preceding sections has been to demonstrate the relationship between several known results about the asymptotic phase dynamics near transversally stable limit cycles and quasiperiodic invariant tori, on the one hand, and the adjoint conditions obtained from the analysis of a problem Lagrangian, on the other hand. As an example, we have found that commonly adopted normalization conditions result directly from the choice of monitor functions and sought sensitivities. Indeed, our analysis extends the treatment beyond the orbitally stable case to the more general instance of normal hyperbolicity. In this case, while it may no longer be possible to define a concept of asymptotic phase, the association persists between the solution to the adjoint conditions and a continuous family of complementary projections onto and transversally to the periodic orbit or quasiperiodic invariant torus.

Although the analysis in Section~\ref{sec:quasiperiodic invariant tori} was restricted to the case of a two-dimensional torus, it generalizes with minimal effort to the case of an $m$-dimensional torus simply by considering a rotation vector $\rho\in\mathbb{R}^{m-1}$ and a vector-valued angle $\phi\in\mathbb{S}^{m-1}$, integrals over $\mathbb{S}^{m-1}$ rather than $\mathbb{S}$, and additional phase conditions. For the transversally stable case, the results reduce to those of~\cite{demir2010}, which are presented there without derivation. A further generalization to the case of a piecewise-smooth dynamical system is open for investigation (cf.~\cite{Szalai09}).

For the case of dynamical systems with time delay, a preliminary analysis in~\cite{Ahsan21} showed the relationship between an adjoint variable and the sensitivity of the orbital duration to violations of the periodic orbit constraint. This same analysis can, of course, be generalized according to the treatment in this paper without significant additional effort. More interesting is to consider the adjoint-based sensitivity analysis for a normally hyperbolic quasiperiodic invariant torus of a smooth dynamical system with delay. While the adjoint-based approach was used in~\cite{ahsan2020optimization} to search for optimal quasiperiodic invariant tori in an example system with delay, the analysis lacked the theoretical rigor of the present treatment, in particular as relates to the solvability of the adjoint conditions. This appears to be a worthwhile problem to pursue.

Finally, we briefly discuss opportunities associated with further development of \textsc{coco} to support analyses of the form considered here. For example, while the example in Section~\ref{sec:invc:coco} relied on low-level functionality (\mcode{coco_add_func}, \mcode{coco_add_pars}, \mcode{coco_add_glue}, and \mcode{coco_add_adjt}), it is straightforward to build a composite constructor that encapsulates the individual steps of constraint construction and a separate composite constructor that encapsulates the individual steps of constructing contributions to the adjoint conditions. Such development would then benefit from the encoding of utility functions for generating an initial solution guess, extracting solution data, and graphing families of solutions and sensitivities.

Indeed, a similar construction could apply to the case where the map $M$ is only implicitly defined by the solution to a differential equation. In this case, the zero problem $M_\mathrm{res}=0$ is replaced with the collocation problem for the corresponding trajectory segment. Rather than \mcode{coco_add_func} and \mcode{coco_add_adjt}, here calls to the constructors \mcode{ode_isol2coll} and \mcode{adjt_isol2coll} would instantiate individual trajectory segments and the corresponding contributions to the adjoint conditions. Finally, the index arrays \mcode{muidx} and \mcode{maidx} would be constructed from a subset of the indices associated with each trajectory instance, specifically those corresponding to the initial and final points on the trajectory segment and the corresponding parameter values. The implementation in \textsc{coco} of the higher-dimensional case should follow the same pattern, with necessary modifications to the discretization, now over $\mathbb{S}^{m-1}$, and the addition of a required number of phase conditions.

%For acknowledgements section, please don't number the section, please begin it with \section*{Acknowledgements}
\section*{\jrem{Code availability}}
\jrem{The code included in this paper constitutes fully executable scripts. Complete code, including that used to generate the results in Fig.~\ref{fig:invc}, is available at \url{https://github.com/jansieber/adjoint-sensitivity2022-supp}.}
\section*{Acknowledgments}
The second author's research is supported by the UK Engineering and Physical Sciences Research Council (EPSRC) grants EP/N023544/1 and EP/V04687X/1.

\bibliographystyle{plainnat}
\bibliography{references}

\medskip
% The data information below will be filled by AIMS editorial staff
Received xxxx 20xx; revised xxxx 20xx.
\medskip

\end{document}